\documentclass[BCOR8mm,DIV14]{scrartcl}
\usepackage{bbm}
\usepackage{epic}        
\usepackage[english]{babel}
\usepackage{amsmath,amssymb,amsthm}
\usepackage[all]{xy}
\newcommand{\labto}[1]{\xrightarrow{\labelstyle\textstyle #1}}
\newtheorem{theorem}{Theorem}
\newtheorem{corollary}{Corollary}
\newtheorem{proposition}{Proposition}
\newtheorem{lemma}{Lemma}
\newtheorem{definition}{Definition}
\newtheorem{remark}{Remark}
\begin{document}
\title{A Boundary Operator for Simplices}
\author{VOLKER  W. TH\"UREY  \\  Rheinstr. 91  \\   28199 Bremen, Germany
     \thanks{volker@thuerey.de \ \ T: Germany, (0)421/591777}   }
    \date{\today}
\maketitle
    \centerline{MSC-2010: 51M04, 55N20}
    \centerline{Keywords: standard simplex, boundary operator}
  \begin{abstract}
         \centerline{Abstract}
         We generalize the very well known boundary operator of the  ordinary singular homology theory,
         defined in many books about algebraic topology.
         We describe a variant of this ordinary simplicial boundary operator where the usual
         boundary $(n-1)$-simplices of each $n$-simplex are replaced by combinations of internal $(n-1)$-simplices
         parallel to the faces.  This construction may lead to an infinite class of extraordinary
         non-isomorphic homology theories.  We show further some interesting constructions on the standard simplex.
   \end{abstract}
  \tableofcontents
	    \section{Introduction}
	
	         \qquad In their famous book  { \it `Foundations of Algebraic Topology' } {\bf \cite{Eil-Stee}}
	         Samuel Eilenberg and Norman Steenrod presented a new method to distinguish topological
	         spaces.  Since this  time  the  { \it `singular homology theory' } is a very useful and successful method
	         in mathematics and  in other fields of science.
	         This `theory' is in fact a sequence of functors from pairs of topological spaces to the category
	         of abelian groups.  It begins with a definition of a  { \it `boundary operator' }  $\partial_{n}$,
	         that means  if  $\Delta_n$ is the standard simplex and $T$ is any continuous map from  $\Delta_n$
	         into a  topological space $X$,  then  $\partial_{n}(T) $ can be considered as the generated map if
	         we restrict  $T$ to the topological boundary of  $\Delta_n$ as the domain of $T$ instead of the entire
	         $\Delta_n$.
	
	    This construction is basicly done by elementary calculations in the $n$-dimensional real space
	         $ \mathbbm{R}^{n}$. In this paper we generalize this construction.
	         Here our  `boundary operator'  $\partial_{n}$ is determined  not only by the topological boundary
	         but also by parts of the interior of the standard simplex.
	         The author took the idea from a similar work which deals with cubical homology, see  {\bf \cite{thuerey}}.
	         This was the natural way because cubes are easier to handle than simplices.
	
	    In the ordinary  singular homology theory  the boundary operator is constructed by taking the topological
	        boundary of an  $n$-dimensional standard simplex as a linear combination of   $n + 1$  simplices of
	        dimension $(n-1)$, (the  {\it faces}), provided with alternating signs.
          We  generalize this by taking a linear combination of a fixed number $L+1$
          of   $(n-1)$-dimensional simplices parallel to each of its $n+1$ faces, provided  with a coefficient tuple
          $ \vec{m} := \left( m_{0}, m_{1}, m_{2},   \dots   ,m_{L} \right)$.
          Note that for a fixed  $L > 0$, `our' boundary operator   $_{\vec{m}}\partial_{n}$   maps not only the
          topological boundary but also parts of the interior of the standard simplex, in contrast to the
          classical singular homology theory.
	
      In the paper we use the customary notations
          $ \mathbbm{N} :=  \{ 1, 2, 3, \dots  .\} ,   \mathbbm{N}_0 := \mathbbm{N} \cup \{ 0 \} , \
          \mathbbm{Z} :=  \{ \dots , -3, -2, -1, 0, 1, 2, 3 , \dots   \}$,  and $\mathbbm{R}$ for the real numbers.
          We shall use  the brackets    $\left( \ldots \right)$ and   $\left[ \ldots \right]$  for tuples and to
          structure   text and formulas,    $\left[ a,b \right]$     also for intervals.
          The brackets   $\left\langle \ldots  \right\rangle$    will be needed for the  boundary operator. \\
         	Let for all   $ n \in  \mathbbm{N}_0$   and all   $ j \in  \{ 0,1,2,  \ldots , n \}:
	         \ e_j := (0, \ldots , 0, 1 , 0 , \ldots , 0 ) \in \mathbbm{R}^{n+1}$
		      (with a single $1$ at the $j^{th}$ place)  be the $j^{th}$-standard unit vector of the
		       $\mathbbm{R}^{n+1}$.   \   Let
	  \begin{align*}       	
	     \Delta_{n} \ := \
	         \left\{ \left( x_{0}, x_{1}, x_{2}, \;  \dots , x_{n} \right) \in  \mathbbm{R}^{n+1} \; | \;
     	     0 \leq x_j \leq 1 \,  \text{for }  \ j \in  \left\{0,1,2,  \dots , n \right\}  \text{ and  }
		       \sum_{j=0}^{n}   x_j  = 1 \:  \right\} .
	  \end{align*}   	
		  That means that
	        $\Delta_{n}$ is the {\it  oriented  $n$-dimensional  standard simplex }
	        with the usual euclidian topology, i.e.  the convex hull of the   $n+1$   standard unit vectors
	        $e_0, e_1, e_2,  \dots  ,  e_n$     of  the real vectorspace  $\mathbbm{R}^{n+1}$.
	        The elements of   $ \{ e_0, e_1, e_2,  \dots  ,  e_n  \}$   are  called the { \it vertices } of $\Delta_{n}$.
	        Note that   $\Delta_{n} \subset \mathbbm{R}^{n+1}$.
	        The space $\Delta_{1}$  is homeomorphic to   $ { \bf I } :=  [ 0, 1 ]$, the unit interval,  and
	        $\Delta_{0}  =  \{ 1 \}$.   \    	
    	    Let  for all  $n \in \mathbbm{N}_{0}$  and for all topological spaces $X$
	 	  $${\cal SS}_{n}(X) := \left\{ \; T : \Delta_{n} \rightarrow  X \; | \;  T \: \text{is continuous} \; \right\}. $$
      ${\cal SS}$    means { \it `Singular Simplices'}.  Moreover all maps we'll use are continuous.  \
        Let  $ {\cal R}$  be a commutative ring with unit $1_{\cal R}$. We define
	  $${\cal F}\left({\cal R}\right)_{n}(X) :=  \text{The free } {\cal R}\text{-module  with the basis} \
	               {\cal SS}_{n}(X) , $$
	      and let   ${\cal F}\left({\cal R}\right)_{-1}(X)   := \{ 0 \}$.
		    Every  {\tt u} $\in {\cal F}\left({\cal R}\right)_{n}(X) $  is called a  { \it `chain'}.
		    Let $\mathsf{TOP}$  be the  category  of topological spaces and continuous maps as
	      morphisms. That means  $ ( f: X \rightarrow Y ) \in \mathsf{TOP}$  if and only if  $X$ and $Y$ are
	      topological spaces and  $f$  is continuous.
   	    Let ${\cal R}$-$\mathsf{MOD}$   be the  category of   ${\cal R}$-Modules.

	      To describe  a homology theory  in one sentence, we say that it is a sequence
	      $ ({\cal H}_n)_{ n \geq 0 }$ of functors,
	                $$  {\cal H}_n: \mathsf{TOP} \longrightarrow  {\cal R} \text{-}\mathsf{MOD}  \, $$
	      with some additional properties, the  `axioms'.	
	      For more detailed information about singular homology theory see, for instance,  {\bf \cite{Hatcher}},
	      {\bf \cite{Rotman}}, {\bf \cite{Spanier}}, {\bf \cite{StoeZie}},  {\bf \cite{Vick}},
	      or other books about this topic. For cubical singular homology theory see {\bf \cite{Massey}}.
	
	  In this paper we create an infinite set of different boundary operators, i.e. for $ L = 1 $ and for each
	     coefficient pair $ \vec{m} := \left( m_{0}, m_{1} \right)$ we construct a chain complex, and then we can take
	     in each  dimension the quotient module   kernel/image.
	     In this way  we generate a sequence  $ (_{\vec{m}}{\cal H}_n)_{ n \geq 0 }$  of functors.  But to get an
	     `extraordinary homology theory' the {\it Excision Axiom }  and the
	     {\it Homotopy Axiom } are  missing. The proofs of  both axioms seem to be difficult. Even if someone  will
	     have success to prove this,  it is uncertain whether there is any  application  of this extraordinary
	     homology theory. An old paper { \bf \cite{BuCoFl}} shows that the  homology groups of an  extraordinary
	     homology which is created by a chain complex  can  be expressed by a  product of singular homology groups.
	     That means that all  informations we can  get  about a topological space from this new homology theory
	     we already are able to get from the known singular homology. 
	     Hence, the meaning of the present paper lies more in the considerations about the  standard simplex
	     which are made in the   fourth  and fifth section.
	
	  Very briefly we describe our work as follows.  For a fixed $L \in \mathbbm{N}_0$
	    and for  all fixed tuples $\vec{m} : = ( m_{0}, m_{1}, m_{2}, \dots  ,m_{L} ) \in {\cal R}^{L+1}$
  	  we try to construct a functor  $_{\vec{m}}{\cal H}_{n}: \mathsf{TOP} \longrightarrow {\cal R}$-$\mathsf{MOD}$
  	  for  all   $ n \in \mathbbm{N}_0$.  We shall have a complete success  for  $L= 1$, while for $L > 1$ we are
  	  missing a   set of homeomorphisms on the standard simplex $ \Delta_n$. Our construction is a generalization of
  	  the ordinary boundary operator in the usual singular homology theory, i.e. the case $ L = 0$.  \\
	\section{The Boundary Operator}
	   Now we shall define  for fixed  $  L  \in  \mathbbm{N}_0$  and a fixed coefficient tuple
         	  $\vec{m}  = \left( m_{0}, m_{1}, \dots  ,m_{L} \right) \:  \in \: {\cal R}^{L+1}$
	          for all  $ n \in \mathbbm{N}_{0}$ the  `boundary operators'
	         	$_{\vec{m}}\partial_{n}: {\cal F}\left({\cal R}\right)_{n}(X)
	         	         \stackrel{ }{\rightarrow}{\cal F}\left({\cal R}\right)_{n-1}(X) $, i.e. they will have
	          the property  $_{\vec{m}}\partial_{n} \circ \, _{\vec{m}}\partial_{n+1}= 0$.
	
	    Assume that we  have fixed a number  $ L \in \mathbbm{N}_0$  and  an  $( L+1)$-tuple
	          of ring    elements
            $\vec{m} = \left( m_{0}, m_{1}, m_{2}, \dots  ,m_{L} \right) \in  {\cal R}^{L+1}$.
	          For all $n \in \mathbbm{N}$, for  all basis  elements
	          $T  \in  {\cal SS}_{n}(X)$,  for all  $i \in  \{ 0, 1, 2, \dots ,  L \}$  and  all
	          $  j  \in \left\{ 0, 1, 2,  \dots  ,n \right\}$  let the map
	          $\left\langle T\right\rangle_{ L,\: n,\;  i,\; j \;}$ be an element of  ${\cal SS}_{n-1}(X)$
	          by defining  for all vectors
	          $( x_{0} ,\  x_{1} ,\  \dots \; , x_{n - 1} )  \in  \Delta_{n-1} \subset \mathbbm{R}^{n}$
	     \begin{eqnarray*}
	          \left\langle T\right\rangle_{\: L, \:  n,\;  i,\; j \;}
	          (\: x_{0},  x_{1} ,  \ldots , x_{j - 1},  x_{j}, \ldots ,  x_{n - 1} )\  := \
	          T ( y_{0} ,  y_{1} , \ldots , y_{j - 1} , v , y_{j},  \ldots ,  y_{n - 1} ) , \\
	         \text{with}  \quad v  :=  \frac{i}{(L+1) \cdot (n+1)} \  , \quad \text{and  for all}  \quad
	            k  \in  \{ 0, 1, 2, \dots ,  n-1 \}   \quad \text{let}  \quad   y_{k} := (1- v)  \cdot x_k .
	     \end{eqnarray*}
	         $\left\langle T\right\rangle_{\: L, \:  n,\;  i,\; j \;}$ is an element of
	             ${\cal SS}_{n-1,X}$.
	
	        Assume for all $L \in \mathbbm{N}_0$,  for all   $n \in  \mathbbm{N}_{0}$, and
	             for all $i \in  \{ 0, 1, 2, \dots , L \}$ the existence of a  special
	             homeomorphism  $\Theta_{L,n,i}$ on $\Delta_{n}$  which  we  shall construct later.
	             (These  homeomorphisms will be necessary for the proof of
	              $_{\vec{m}}\partial_{n} \; \circ \; _{\vec{m}}\partial_{n+1} =  0$).
	             We define for  all  $n  \in \mathbbm{N}$ and for an arbitrary  $ T \in {\cal SS}_{n}(X)$:
 	    \begin{equation}
	          _{\vec{m}}\partial_{n}	 (T)
  	        :=      \sum_{j =  0}^{n}   (-1)^{j}  \cdot \sum_{i = 0}^{L}
  	        m_{i} \cdot  \left[ \left\langle T\right\rangle_{\: L,\: n,\;  i,\; j \;} \circ \Theta_{L,n-1,i} \right]
	    \end{equation}
	          and let $_{\vec{m}}\partial_{0}(T) := 0$ be the only possible map.  \\
	          Define \ $_{\vec{m}}\partial_{n}$:
	          ${\cal F}\left({\cal R}\right)_{n}(X) \stackrel{ }{\rightarrow}
	                              {\cal F}\left({\cal R}\right)_{n-1}(X)$  by linearity.
	          See the following Figure \ref{picture1}, which illustrates the case
	          ${\cal R} := \mathbbm{Z}$,  $n := 2$ , $L := 1$ , ${\vec{m}} := (9,4)$ , and  $T := id(\Delta_2)$.  \\
	          On the left hand side you see the two-dimensional standard simplex
	           $\Delta_2$, the right hand side shows $ _{\vec{m}}\partial_2 (T)$, i.e. the images of
	           six one-dimensional standard simplices
	           $  \left\langle T\right\rangle_{1,\: 2,\;  i,\; j \;}, i \in \{0,1\}$ and $j \in \{0,1,2\}$,
	           multiplied by  coefficients  9 and 4, elements of the ring ${\cal R} = \mathbbm{Z}$.  We have   \\
            $  _{(9,4)}\partial_2 (T) \ = \
            9 \cdot  [   \left\langle T\right\rangle_{\: 1,\: 2,\;  0,\; 0 \;}   \circ    \Theta_{1,1,0}   ]  +
            4 \cdot  [   \left\langle T\right\rangle_{\: 1,\: 2,\;  1,\; 0 \;}   \circ    \Theta_{1,1,1}   ]  -
            9 \cdot  [   \left\langle T\right\rangle_{\: 1,\: 2,\;  0,\; 1 \;}   \circ    \Theta_{1,1,0}   ]   \\
             {  }     \qquad  \qquad \  \
            - 4 \cdot  [   \left\langle T\right\rangle_{\: 1,\: 2,\;  1,\; 1 \;}   \circ    \Theta_{1,1,1}   ]  +
            9 \cdot  [   \left\langle T\right\rangle_{\: 1,\: 2,\;  0,\; 2 \;}   \circ    \Theta_{1,1,0}   ]  +
            4 \cdot  [   \left\langle T\right\rangle_{\: 1,\: 2,\;  1,\; 2 \;}   \circ    \Theta_{1,1,1}   ] $ .
         \newpage  $ $  \\
	\begin{figure}[h]
    \setlength{\unitlength}{1cm}
    \begin{picture}(8,3)
         \put(0,0){\line(1,0){6.0}}      \put(-0.4,-0.3){$e_0$}    \put(6.1,-0.3){$e_1$}   \put(2.8,4.7){$e_2$}
         \put(0,0){\line(2,3){3.0}}
         \put(6,0){\line(-2,3){3.0}}      
         \put(8,0){\line(1,0){6.0}}      \put(7.6,-0.3){$e_0$}    \put(14.1,-0.3){$e_1$}   \put(10.8,4.7){$e_2$}
         \put(8,0){\line(2,3){3.0}}
         \put(14,0){\line(-2,3){3.0}}       \put(11.4,1.8){4}      \put(12.0,2.4){9}
         \put(8.6,0.9){\line(1,0){4.8}}
         \put(9,0){\line(2,3){2.5}}         \put(9.7,1.8){$-4$}    \put(9.1,2.4){$-9$}
         \put(13,0){\line(-2,3){2.5}}       \put(11.0,0.96){$4$}   \put(11.0,0.06){$9$}
   \end{picture}  \\
     \caption{}\label{picture1} 
  \end{figure}	  
     \\  \\
      	Now we define a set of important maps which will play a central part in the proof that
           	$_{\vec{m}}\partial_{n}$ is a `boundary operator', i.e.  \
	          $_{\vec{m}}\partial_{n} \, \circ \,	_{\vec{m}}\partial_{n+1} =  0$.
	\begin{definition}   \rm
            For every  $n \in \mathbbm{N}_0$ let $id$ be the identical map on  $\Delta_n$.
	          For all $n \in  \mathbbm{N}$ and for all  $i, k \in \{ 0, 1, \ldots , L \}$,
            we look at some injective maps,  (to  make it  better  readable, we omit the parameter $L$  in  such
            expressions  as  $\left\langle T \right\rangle_{ L, n, i, j } \ \text{and} \ \Theta_{L,n,i} $),
      $$     \left\langle id \right\rangle_{ n+1, i, j } \ \circ \ \Theta_{n,i} \ \circ \
  	            \left\langle id  \right\rangle_{ n, k, p } \ \circ \ \Theta_{n-1,k}  \quad \text{for}  \quad 	
  	             j \in \{ 0, 1, \ldots , n+1 \} ,   p \in \{ 0, 1, \ldots , n \}   , $$
  	        they are injective  maps from  $ \Delta_{n-1} \ \text{to}  \ \Delta_{n+1}$.  (Note that still the maps
  	        $ \Theta_{n,i}$ are not defined, we only are describing some of its  properties.)   \\
  	       Now let $ j, p \in \{ 0, 1, \ldots , n \}$ with  $j \leq p$.
           If we have  the equality of the following two maps,
	     \begin{equation}
	            \left\langle id \right\rangle_{\: n+1,\;  i,\; j \; } \ \circ \ \Theta_{n,i}  \ \circ \
  	            \left\langle id \right\rangle_{\: n,\;  k,\; p \; } \ \circ \ \Theta_{n-1,k}  \quad  =  \quad 	
  	         \left\langle id \right\rangle_{\: n+1,\;  k,\; p+1 \; } \ \circ \  \Theta_{n,k} \ \circ \
  	          \left\langle  id  \right\rangle_{\: n,\;  i,\; j \; } \ \circ  \  \Theta_{n-1,i}   \   ,
       \end{equation}  	
             then we abbreviate this important equation by  {\tt   EQUATION$_{n,j \leq p,i,k}$},
             for every fixed $n \in \mathbbm{N},  j, p \in \{ 0, 1, \ldots, n \}$
             with  $j \leq p$, and $i, k \in \{ 0, 1, \ldots , L \}$.     \hfill  $\Box$
	\end{definition}
	\begin{remark}   \rm     \label{remark eins}
	          You will find a corresponding equation in every book about simplicial homology theory,
	          e.g. in  {\bf \cite[p.65]{Rotman}} it appears as  `If $ j < p $,   the face maps satisfy
	          $ \varepsilon^{n+1}_j \circ \varepsilon^{n}_{p-1} = \varepsilon^{n+1}_p \circ \varepsilon^{n}_j$'.
	\end{remark}
     Now we formulate the theorem that if the { \tt   EQUATION$_{n,j \leq p,i,k}$ } holds,
            the above construction leads  to a `boundary operator', that means
            $_{\vec{m}}\partial_{n} \; \circ \;	_{\vec{m}}\partial_{n+1} = 0$.
   \begin{theorem}  \sf   \quad   \label{theorem eins}
		        Let $L$  be a fixed element  of  $ \mathbbm{N}_0$, and   let
		        $\vec{m} := \left(m_{0}, m_{1}, m_{2}, \dots  ,m_{L} \right)$  be a  fixed tuple from $ {\cal R}^{L+1}$.
            In addition  we assume  the  following property:  For every  $n \in  \mathbbm{N}$, for all
            $  i , k \in \{ 0, 1, \ldots , L \}$  the equation  { \tt   EQUATION$_{n,j \leq p,i,k}$ } holds
            for all  $j , p \in \{ 0, 1, \ldots , n \}$  with  $j \leq p$.
	          Then we  have for all $n \in  \mathbbm{N}_0$  for all  $T \in {\cal SS}_{n+1}(X)$
	          (i.e. $ T : \Delta_{n+1} \rightarrow X $ is continuous):
	         $$_{\vec{m}}\partial_{n} \;	 \circ \;    	 _{\vec{m}}\partial_{n+1} (T) \ = \ 0 \ . $$
	 \end{theorem}
	 \begin{proof} \quad
	       The statement is trivial for $n = 0$,  so let $n \in  \mathbbm{N}$.
	            Note that $\left\langle T \right\rangle_{ n+1, i, j }$ is a map with the domain $ \Delta_n $,
	            and note that   $\left\langle T \right\rangle_{ n+1, i, j } =
	                             T \circ \left\langle id \right\rangle_{n+1, i, j }$ . \  We have
	    \begin{align*}
	         &  _{\vec{m}}\partial_{n} \;	 \circ \; 	 _{\vec{m}}\partial_{n+1} (T)
	                \ = \  _{\vec{m}}\partial_{n} \;	\left( \ \sum_{j =  0}^{n+1} (-1)^{j} \cdot \sum_{i = 0}^{L}
  	              m_{i} \cdot \left[\left\langle T\right\rangle_{\: n+1,\;  i,\; j \;} \ \circ \ \Theta_{n,i} \right]
  	              \right) \\
        =  & \ \  \sum_{p =  0}^{n}   (-1)^{p} \cdot \sum_{k = 0}^{L}    m_{k} \cdot  	
	             \sum_{j =  0}^{n+1}   (-1)^{j} \cdot \sum_{i = 0}^{L} m_{i} \cdot
	             \left[ \left\langle \left\langle T \right\rangle_{\: n+1,\;  i,\; j \;} \ \circ \
	             \Theta_{n,i}  \right\rangle_{\: n,\;  k,\; p \; } \ \circ \ \Theta_{n-1,k} \right]    \\
        =  & \ \  \sum_{p =  0}^{n}   (-1)^{p} \cdot \sum_{k = 0}^{L}    m_{k} \cdot  	
	             \sum_{j =  0}^{n+1}   (-1)^{j} \cdot \sum_{i = 0}^{L}
               m_{i} \cdot  \left[ \left\langle T \right\rangle_{\: n+1,\;  i,\; j \; } \ \circ \ \Theta_{n,i} \
               \ \circ \ \left\langle id \right\rangle_{\: n,\;  k,\; p \; } \ \circ \ \Theta_{n-1,k} \right]
     \end{align*}
    \begin{align}      \label{equation drei}
               =   \ & \ \sum_{p =  0}^{n}  \sum_{j =  0}^{n+1}   \sum_{i , k = 0}^{L}   	
	                (-1)^{j+p}  \cdot  m_{i} \cdot m_k  \cdot  \left[
	                T \ \circ \ \left\langle id\right\rangle_{\: n+1,\;  i,\; j \; } \ \circ \ \Theta_{n,i} \ \circ \
  	              \left\langle  id  \right\rangle_{\: n,\;  k,\; p \; } \ \circ \ \Theta_{n-1,k} \right]   \  .
     \end{align}
     The sign depends only on  $ j $ and  $p$.     \\
  	           The set $ M :=  \left\{ 0, 1, 2, \dots  ,n, n+1 \right\} \times \left\{ 0, 1, 2, \dots  ,n \right\}$
  	           contains    $(n+2)\cdot(n+1) $ elements. With  $M_{small}  :=  \{ (j,p) \in M  |  j \leq p \}$
  	           and  $M_{big}  :=  \{ (j,p) \in M  |  j > p \}$ ,   we have  $ M =  M_{small} \cup  M_{big}$,
  	           and  $M_{small} \cap M_{big} = \emptyset$.  The map
  	           ${ \cal B }: \ M_{small} {\rightarrow}   M_{big} , \ (j,p) \mapsto (p+1,j)$  is  bijective.
  	           We have
    	\begin{align*}
    	     _{\vec{m}}\partial_{n} \; \circ \; _{\vec{m}}\partial_{n+1} (T) \
         & =   \sum_{ (j,p) \in M_{small}}   \ \sum_{i , k = 0}^{L}   	
	         (-1)^{j+p} \cdot  m_{i} \cdot m_k  \cdot
	         \left[ T \circ \left\langle id\right\rangle_{ n+1, i, j } \circ \Theta_{n,i} \circ
  	              \left\langle  id  \right\rangle_{n, k, p } \circ \Theta_{n-1,k} \ \right]   \\
          & +  \sum_{ (j,p) \in M_{big}}   \ \sum_{i , k = 0}^{L}   	
	         (-1)^{j+p} \cdot  m_{i} \cdot m_k  \cdot
	         \left[ T \circ \left\langle id\right\rangle_{ n+1, i, j } \circ \Theta_{n,i} \circ
  	              \left\langle  id  \right\rangle_{ n, k, p } \circ \Theta_{n-1,k} \ \right] .
  	  \end{align*}
  	  We rename the elements   $(j,p) \in M_{big}$   in   $(p+1,j)$.   Further, because of
  	  $ \sum_{i , k = 0}^{L}$, we can exchange the parts of  $i$ and $k$ in the second sum.
  	  Hence we get that \ $_{\vec{m}}\partial_{n} \; \circ \; _{\vec{m}}\partial_{n+1} (T) \ = $
  	  \begin{eqnarray}   	     	 \label{equation vier}
	        &   &    \sum_{ (j,p) \in M_{small}}   \ \sum_{i , k = 0}^{L}   	 
	              (-1)^{j+p} \cdot  m_{i} \cdot m_k  \cdot
	         \left[ T \circ \left\langle id\right\rangle_{ n+1, i, j } \circ \Theta_{n,i} \circ
  	              \left\langle  id  \right\rangle_{n, k, p } \circ \Theta_{n-1,k} \  \right]   \\
  	               \label{equation fuenf}
          &  + & \sum_{ (p+1,j) \in M_{big}}    \sum_{i , k = 0}^{L}   	
	            (-1)^{(p+1)+j}    \cdot     m_{k} \cdot m_i  \cdot
	         \left[   T \circ \left\langle id\right\rangle_{ n+1, k, p+1 } \circ \Theta_{n,k} \circ
  	              \left\langle  id  \right\rangle_{n, i, j } \circ \Theta_{n-1,i} \  \right]     
  	 \end{eqnarray}
           Because of the bijection  ${ \cal B }$ of  $M_{small}$  and  $M_{big}$,
  	       every summand in  (\ref{equation vier}) corresponds to another in  (\ref{equation fuenf}).
  	       Because of  {\tt EQUATION$_{n,j \leq p,i,k}$}  and   because of different signs   the
  	       $ (n+1)\cdot(n+2)\cdot(L+1)^{2}$  summands in equation(\ref{equation drei}) cancel pairwise.
  	       Thus it follows  that $_{\vec{m}}\partial_{n - 1} \;	\circ \; _{\vec{m}}\partial_{n }(T) = 0$.
    \end{proof}
    \begin{remark}     \rm     \label{remark zwei}	
	       If you read the above  definition of the boundary operator for the first time,
	       you  may  miss the idea behind it. Here is an attempt to explain it: For any  $(n+1)$-simplex
	       $\Delta \subset \mathbbm{R}^{n}$, the `boundary' $_{\vec{m}}\partial_{n+1}(\Delta)$ is a linear
	       combination of $n$-simplices,  regarded as subsets of $ \Delta $.
	       The set \ $_{\vec{m}}\partial_{n} \circ \, _{\vec{m}}\partial_{n+1}(\Delta)$ \ is a union of
	       $(n-1)$-simplices.  In fact it is the union of the intersections of  two  at a time of the  $n$-simplices
	       of   $_{\vec{m}}\partial_{n+1}(\Delta)$. Every $(n-1)$-simplex of
	       $_{\vec{m}}\partial_{n} \circ \, _{\vec{m}}\partial_{n+1}(\Delta)$   occurs twice.
	       And, by  factors $m_i$  and different signs, they cancel each other. Hence \
	       $_{\vec{m}}\partial_{n} \circ \, _{\vec{m}}\partial_{n+1}(\Delta) = 0 $.    See again the first 
	       Figure \ref{picture1}, where  \ $_{\vec{m}}\partial_{1} \circ \, _{\vec{m}}\partial_{2}(\Delta_2)$ \ 
	       is the union of  0-simplices, i.e.  of dots.   
	    \end{remark}
	   The  shortest  way  to complete  this new boundary operator  is by writing down an  order: \\
	   `The construction of the homeomorphisms  $\Theta_{L,n,i}$  is left to the reader as an easy exercise!'  \\
	   Unfortunately  this is only a joke; indeed, the construction of the $\Theta_{L,n,i}$'s is the most
         difficult part in the  proof of Theorem  \ref{theorem eins}. (In the theorem we assumed the existence of
         these homeomorphisms).  We  mention a solution  for  $L = 0$
      	 because it   is  well known, and we  present a solution  for  the special case  $L = 1$,
      	 which also may show an idea of a  general construction for the case of an arbitrary
         $ L   \in  \mathbbm{N}$, but this is still an open problem and is left to the reader
       `as an easy exercise'.
   \newpage
                                        	
   \section{Beginning of the Induction}   	
      	 We formulate a trivial but important lemma.	
      \begin{lemma}       \label{lemma eins}
	          The map   $\left\langle T\right\rangle_{\: L, \:  n,\;  i,\; j \;} $   has the important property
	          that it respects permutations.  As before, let us fix  $ L \in \mathbbm{N}_0, \ n \in \mathbbm{N}$, and
	          $\vec{m}  = \left( m_{0}, m_{1}, m_{2}, \dots  ,m_{L} \right) \in {\cal R}^{L+1}$.
	          Let   $i \in \{ 0,1,2, \ldots , L \}$, and $j, \widetilde{j} \in  \{ 0,1,2, \ldots , n \}$,
	          with (for instance) $ j < \widetilde{j}$.    \  \
	          We take an arbitrary  $T \in {\cal SS}_{n}(X)$, and an  $n$-tuple
	          $(\: x_{0}, x_{1}, x_2, \dots , x_{n - 1}) \in  \Delta_{n-1}$. \  If   \\
	      \centerline{$ \left\langle T\right\rangle_{\: L, \:  n,\;  i,\; j \;}
	                 \left( \: x_{0},\, x_{1}, \; \dots \; , x_{j-1}, \, x_j, \, \dots \, , x_{n - 1}\;  \right)
	                 \  =  \  T \left( \: y_{0} ,\  y_{1} ,\  y_{2} , \dots , \;y_{j - 1} ,\; v ,\; y_{j} , \; \dots  \;
	                    , \; y_{n - 1 \:}\;  \right) $,  }
	        as it is defined above,  we have   \\
	          \centerline{ $\left\langle T\right\rangle_{\: L, \:  n,\;  i,\; \widetilde{j} \;}
	                  \left( \: x_{0} ,\  x_{1} ,  \;   \ldots \; ,   x_{n - 1}\;  \right) \
	                 =   T  \left( \: y_{0} ,\  y_{1} ,    \ldots   ,   \;y_{j - 1} , \; y_{j} , \ldots ,
	                 \; y_{\widetilde{j} - 1} ,\; v ,\; y_{\widetilde{j}} , \; \ldots  \; , \; y_{n - 1 \:}\; \right) $.}
	                 \\     And if   $\vartheta$  is an arbitrary permutation of $ \{ 0,1,2, \ldots , n-1 \}$, \
	                  then it holds that  \\
	          $ \left\langle T\right\rangle_{\: L, \:  n,\;  i,\; j \;}   \left( \: x_{\vartheta(0)} ,\
	                   x_{\vartheta(1)} ,  \; \ldots \; , x_{\vartheta(n-1)}\; \right) \
	                   =    T   \left( \: y_{\vartheta(0)} ,\ \ldots , \;y_{\vartheta(j-1)} ,\; v , \
	                    y_{\vartheta(j)} , \; \dots   , \; y_{\vartheta(n-1) \:} \;  \right) \ $   \\     and
	            \centerline{$ \left\langle T\right\rangle_{\: L, \:  n,\;  i,\; \widetilde{j} \;}
	                  \left( \: x_{\vartheta(0)},\, x_{\vartheta(1)}, \; \ldots \;, x_{\vartheta(j-1)}, \;
	                  x_{\vartheta(j)}, \; \dots , \; x_{\vartheta(\widetilde{j}-1)},\; x_{\vartheta(\widetilde{j})},
	                   \; \ldots , \, x_{\vartheta(n-1)}\; \right) \ $ } \\
	            \centerline{ $ =   T  \left( \:   y_{\vartheta(0)} ,\  y_{\vartheta(1)} , \,  \ldots , \,
	                      \;y_{\vartheta(j-1)} ,\; y_{\vartheta(j)} , \; \dots   ,
	                      \;y_{\vartheta(\widetilde{j}-1)} ,\; v \; , y_{\vartheta(\widetilde{j})} , \;
	                     \ldots  \; ,    \; y_{\vartheta(n-1) \:}\; \;  \right) $.  }
	\end{lemma}
	\begin{proof}
	         These facts are trivial, but it is necessary to mention them.
	\end{proof}
	    Now we look at special cases of \ $ L $.  \\    \\	
		  { \bf  The case L := 0}.
	         Let  for all  $ n  \in  \mathbbm{N}_{0}$:    \
	       $\Theta_{0,n,0} :=  id({\Delta_{n}})$, the identity map on  $\Delta_{n}$, and  with
	       $\vec{m} : = ( 1 )$ we get the well known ordinary  singular simplicial
	       homology theory  with the coefficient module  ${\cal R}$, which was  introduced in {\bf \cite{Eil-Stee}},
	       and it  is also described for  instance in  {\bf \cite{Hatcher}}, or {\bf \cite{Rotman}}.
	       \\   \\
	 	  { \bf The case L := 1}. \
	       This case is a little bit more complicated, and we shall need the rest of the paper to explain  it.
	
	   We must construct two homeomorphisms $\Theta_{1,n,0},  \Theta_{1,n,1}$  on $\Delta_{n}$,  for each
	       $n \in  \mathbbm{N}_0$. We have the singleton $\Delta_{0} = \{1\}$, hence
	       $\Theta_{1,0,0} = \Theta_{1,0,1} := id(\{1\})$,    of course.
	
	   The homeomorphisms  $ \Theta_{1,n,0}$ will  be described by a general construction, and the maps
	       $ \Theta_{1,n,1}$  will  be defined  by  induction on $n$ to  make the equations
	       { \tt   EQUATION$_{n,j \leq p,i,k}$ }  true.  \\  \\
	   \underline{Beginning  of   the   induction.} \
	         For  $ n=1$  let us construct two homeomorphisms  $ \Theta_{1,1,0} , \Theta_{1,1,1}$   on
	         $ \Delta_{1}$.  Use the auxiliary homeomorphisms
	         $ \eta, \kappa: [0,1] \stackrel{\cong}{\longrightarrow}  [0,1]$,
	        let $ \eta$   be the polygon through	four points
	        $ \left\{(0,0), \left(\frac{1}{4},\frac{1}{6}\right), \left(\frac{3}{4},\frac{5}{6}\right), (1,1) \right\}$,
	        and  let  $\kappa$  be the polygon through	four points
	        $ \left\{(0,0), \left(\frac{1}{4},\frac{1}{5}\right), \left(\frac{3}{4},\frac{4}{5}\right), (1,1) \right\}$.
	        Note that both maps $ \eta , \kappa$   are symmetrical at the point
	        $\left(\frac{1}{2},\frac{1}{2}\right)$.  We  have  \
	        $ \eta(x) + \eta(1-x) = 1 =  \kappa(x) + \kappa(1-x) $,   for all   $ x \in  [0,1]$.
	        Further note \ $\eta\left(\frac{1}{4}\right) = \frac{1}{6}$, \ and
	        $\kappa\left(\frac{1}{4}\right) = \frac{1}{5}$.
	  \begin{definition}    \label{zweite Definition}  \quad
	         Define  for  all pairs  $(x, 1-x) \in \Delta_1$  \\
	    \centerline{ $ \Theta_{1,1,0}(x,1-x)   :=   (\eta(x),\eta(1-x)) \quad \text{and} \quad
	                   \Theta_{1,1,1}(x,1-x)   :=   (\kappa(x),\kappa(1-x)) \ .$ }   ${  }$ \hfill $\Box$
	  \end{definition}
	   Now we have to consider $ 4 \cdot 3$ equations  {\tt EQUATION$_{n=1, j \leq p,i,k}$},  for
	        $ i,k \in \{0,1\}, \ \text{and} \ j,p \in \{0,1\}$ with   $j \leq p$.
	        But fortunately, because of  Lemma \ref{lemma eins}, we can fix   $ j = p = 0 $.  The other pairs
	        $(j,p) \in \{ (0,1) , (1,1) \}$  work in the  same manner.  \\
	        We want to prove the equality of (in the following equations we  omit the parameter $L=1$)
	     \begin{equation*}
	        \left\langle id\right\rangle_{ 2, i, j=0 } \ \circ \ \Theta_{1,i} \ \circ \
  	              \left\langle  id  \right\rangle_{1, k, p=0 } \ \circ \ \Theta_{0,k} \  = \
  	       \left\langle id\right\rangle_{ 2, k, p+1=1 } \ \circ \ \Theta_{1,k} \ \circ \
  	              \left\langle  id  \right\rangle_{1, i, j=0 } \ \circ \ \Theta_{0,i} \ \ .
  	    \end{equation*}
	         To begin with we fix $i = 0, k = 1 $, i.e. we consider the {\tt EQUATION$_{n=1, j=0 \leq p=0,i=0,k=1}$}.
	         We have to show that     
    \begin{equation*}
	        \left\langle id\right\rangle_{ 2, 0, 0 } \ \circ \ \Theta_{1,0} \ \circ \
  	              \left\langle  id  \right\rangle_{1, 1, 0 } \ \circ \ \Theta_{0,1} \  = \
  	       \left\langle id\right\rangle_{ 2, 1, 1 } \ \circ \ \Theta_{1,1} \ \circ \
  	              \left\langle  id  \right\rangle_{1, 0, 0 } \ \circ \ \Theta_{0,0} \ \ .
     \end{equation*}    	
	   We have to map the set $ \Delta_0 = \{1\}$,  hence the left hand side of the equation is
	   \begin{center}
	          $\left\langle id\right\rangle_{ 2, 0, 0 } \ \circ \ \Theta_{1,0} \ \circ \
  	        \left\langle  id  \right\rangle_{1, 1, 0 } \ \circ \ \Theta_{0,1}$ \    $ (1) $   \ = \
  	        $\left\langle id\right\rangle_{ 2, 0, 0 } \ \circ \ \Theta_{1,0} \ \circ \
  	        \left\langle  id  \right\rangle_{1, 1, 0 }$   $ (1) $  \\
  	 =    $\left\langle id\right\rangle_{ 2, 0, 0 } \ \circ \ \Theta_{1,0} \ \left(\frac{1}{4},\frac{3}{4} \right) $ \
  	       =   $\left\langle id\right\rangle_{ 2, 0, 0 } \ \left(\frac{1}{6},\frac{5}{6} \right) $  \
        	 =   $ \left(0,\frac{1}{6},\frac{5}{6} \right) \in  \Delta_2 $ .
      \end{center}  	
  	 The right hand side of the equation is
  	    \begin{center}
	          $\left\langle id\right\rangle_{ 2, 1, 1 } \ \circ \ \Theta_{1,1} \ \circ \
  	        \left\langle  id  \right\rangle_{1, 0, 0 } \ \circ \ \Theta_{0,0} \   (1) $   \ = \
  	        $\left\langle id\right\rangle_{ 2, 1, 1 } \ \circ \ \Theta_{1,1} \ \circ \
  	        \left\langle  id  \right\rangle_{1, 0, 0 }$   $ (1) $  \\
  	 =    $\left\langle id\right\rangle_{ 2, 1, 1 } \ \circ \ \Theta_{1,1} \ \left( 0,1 \right) $ \
  	       =   $\left\langle id\right\rangle_{ 2, 1, 1 } \ \left( 0,1 \right) $  \
        	 =   $ \left(0,\frac{1}{6},\frac{5}{6} \right) $, \ \
               hence $ {\tt EQUATION_{1, 0 \leq 0,0,1}}$   holds.
      \end{center}    $ $ 	      \\
  	 More  beautiful is a comutative diagram (Figure \ref{picture2}):  \\  \\  \\  \\ %
  	\begin{figure}[ht]
  \centering
  \setlength{\unitlength}{1cm}
  \begin{picture}(8,3)         
     \put(-3.5,2.0){$\Delta_0$ }   \put(-2.5,2.9){$\vector(2,1){2}$}   \put(-2.0,2.7){$\Theta_{0,0} = id$}  
     \put(-0.0,3.9){$\Delta_0$ }   \put(1.0,4.0){$\vector(1,0){2}$} 
                                   \put(1.5,3.5){$\left\langle id \right\rangle_{\: 1,\;  0,\; 0 \; }$}  
                                   \put(3.8,3.9){$\Delta_1$ }  
     \put(5.0,4.0){$\vector(1,0){2}$}   \put(5.7,3.55){$\Theta_{1,1}$}  \put(7.8,3.9){$\Delta_1$ }
     \put(10.0,3.4){$\left\langle id \right\rangle_{\: 2,\;  1,\; 1 \; }$}    \put(9.0,3.8){$\vector(2,-1){2}$}    
     \put(-2.5,1.2){$\vector(2,-1){2}$}    \put(-2.0,1.0){$\Theta_{0,1} = id$} 
     \put(0.0,-0.1){$\Delta_0$ }       \put(1.0,0.0){$\vector(1,0){2}$} 
                                       \put(1.5,0.3){$\left\langle id \right\rangle_{\: 1,\;  1,\; 0 \; }$}  
     \put(3.8,-0.1){$\Delta_1$ }        \put(5.0,0.0){$\vector(1,0){2}$}
     \put(5.7,0.2){$\Theta_{1,0}$}   \put(7.8,-0.1){$\Delta_1$ } 
     \put(9.0,0.2){$\vector(2,1){2}$}   \put(9.9,0.3){$\left\langle id \right\rangle_{\: 2,\;  0,\; 0 \; }$}   
     \put(11.0,2.0){$\Delta_2$ }  
  \end{picture}   
      \\  \caption{}\label{picture2}   
  \end{figure}  
   \\   \\  
        We get a corresponding  diagram if we replace  the set   $ \Delta_0$   by its single element $(1)$, 
  	  	and we show again the equation {\tt EQUATION$_{n=1, j=0 \leq p=0,i=0,k=1}$} .  \\   \\  \\
    \begin{figure}[ht]
  \centering
  \setlength{\unitlength}{1cm}
  \begin{picture}(9,4)
     \put(-3.5,2.0){$(1)$ }   \put(-2.5,2.9){$\vector(2,1){2}$}   \put(-2.0,2.7){$\Theta_{0,0} = id$}  
     \put(-0.0,3.9){$(1)$ }   \put(1.0,4.0){$\vector(1,0){2}$} 
                              \put(1.5,3.5){$\left\langle id \right\rangle_{\: 1,\;  0,\; 0 \; }$}  
     \put(3.8,3.9){$(0,1)$ }  
     \put(5.0,4.0){$\vector(1,0){2}$}   \put(5.7,3.55){$\Theta_{1,1}$}  \put(7.8,3.9){$(0,1)$ }
     \put(10.0,3.4){$\left\langle id \right\rangle_{\: 2,\;  1,\; 1 \; }$}    \put(9.0,3.8){$\vector(2,-1){2}$}    
     \put(-2.5,1.2){$\vector(2,-1){2}$}    \put(-2.0,1.0){$\Theta_{0,1} = id$} 
     \put(-0.0,-0.1){$(1)$ }       \put(1.0,0.0){$\vector(1,0){2}$} 
                                  \put(1.5,0.3){$\left\langle id \right\rangle_{\: 1,\;  1,\; 0 \; }$}  
     \put(3.8,-0.1){$\left(\frac{1}{4},\frac{3}{4}\right)$ }  
     \put(5.0,0.0){$\vector(1,0){2}$}      \put(5.7,0.2){$\Theta_{1,0}$}                                                     \put(7.8,-0.1){$\left(\frac{1}{6},\frac{5}{6}\right)$ } 
     \put(9.0,0.2){$\vector(2,1){2}$}   \put(9.9,0.3){$\left\langle id \right\rangle_{\: 2,\;  0,\; 0 \; }$}   
     \put(11.0,2.0){$\left(0,\frac{1}{6},\frac{5}{6}\right)$ }    
  \end{picture}   
     \\    \caption{}\label{picture3} 
  \end{figure}          
  	     { $  $   }  \\  \\  
  	      If we exchange $i$ and $k$ we get a symetrical diagram, both sides of the
  	      {\tt EQUATION$_{1, 0 \leq 0,i=1,k=0}$} map the single element $(1)$ to the point
  	      $\left(\frac{1}{6},0,\frac{5}{6}\right)$.  Hence this equation also holds.
  	      Further, the reader may convince himself that the case $ i = k = 0$, i.e. the
  	      {\tt EQUATION$_{1, 0 \leq 0,i=0,k=0}$}  is trivial. Only  the
  	      {\tt EQUATION$_{1, 0 \leq 0,i=1,k=1}$} is missing. \ We map the single element   $ (1) \in \Delta_0 $,
  	      and we prove  {\tt EQUATION$_{1, 0 \leq 0,i=1,k=1}$}, see the following Figure \ref{picture4}.
  	      \\  \\  \\  \\
 \begin{figure}[ht]
  \centering
  \setlength{\unitlength}{1cm}
  \begin{picture}(9,4)
     \put(-3.5,2.0){$(1)$ }   \put(-2.5,2.9){$\vector(2,1){2}$}   \put(-2.0,2.7){$\Theta_{0,1} = id$}  
     \put(-0.0,3.9){$(1)$ }   \put(1.0,4.0){$\vector(1,0){2}$} 
                              \put(1.5,3.5){$\left\langle id \right\rangle_{\: 1,\;  1,\; 0 \; }$}  
     \put(3.8,3.9){$\left(\frac{1}{4},\frac{3}{4}\right)$ }  
     \put(5.0,4.0){$\vector(1,0){2}$}       \put(5.7,3.55){$\Theta_{1,1}$} 
                                            \put(7.8,3.9){$\left(\frac{1}{5},\frac{4}{5}\right)$ }
     \put(10.0,3.4){$\left\langle id \right\rangle_{\: 2,\;  1,\; 1 \; }$}    \put(9.0,3.8){$\vector(2,-1){2}$}    
     \put(-2.5,1.2){$\vector(2,-1){2}$}    \put(-2.0,1.0){$\Theta_{0,1} = id$} 
     \put(-0.0,-0.1){$(1)$ }       \put(1.0,0.0){$\vector(1,0){2}$} 
                                  \put(1.5,0.3){$\left\langle id \right\rangle_{\: 1,\;  1,\; 0 \; }$}  
     \put(3.8,-0.1){$\left(\frac{1}{4},\frac{3}{4}\right)$ }  
     \put(5.0,0.0){$\vector(1,0){2}$}      \put(5.7,0.2){$\Theta_{1,1}$}                                                     \put(7.8,-0.1){$\left(\frac{1}{5},\frac{4}{5}\right)$ } 
     \put(9.0,0.2){$\vector(2,1){2}$}   \put(9.9,0.3){$\left\langle id \right\rangle_{\: 2,\;  1,\; 0 \; }$}   
     \put(11.0,2.0){$\left(\frac{1}{6},\frac{1}{6},\frac{4}{6}\right)$ }    
  \end{picture}         \\     \caption{}\label{picture4}    
  \end{figure}    
  	     { $  $   }  \\    \\ \\  
     As we  mentioned  above, the other pairs
         $(j,p) \in \{ (0,1) , (1,1) \}$ with  $(i,k) \in \{ (0,0) , (1,1), (0,1), (1,0) \}$
         work in the  same manner. Hence we have done the beginning of the induction  for  $ n=1$.  \\
                           
     Before we can continue the induction with
         the construction  of  $ \Theta_{n,0} $  and  $ \Theta_{n,1}$ for greater  $n$,
	       we have to  make some  technical considerations about the standard simplices $ \Delta_n$.  \\   \\
	  \newpage
	  \section{The Geometry of the Standard Simplex}
	        The constructions of this section will be needed in the next one.   Some of the following definitions
	        and propositions can be skipped over during the first reading.  The reader may take a look on the next
	        section first.
	 \begin{definition}   \rm  For all  $ n \in \mathbbm{N}_0$,  let     
       ${ \bf Center}_n  := \left( \frac{1}{n+1},\frac{1}{n+1}, \ \ldots , \ \frac{1}{n+1}\right) \in \Delta_n$, 
            the  center of the $n$-dimensional standard simplex \ $\Delta_n \subset  \mathbbm{R}^{n+1}$. 
             \hfill $\Box$ 
       \end {definition}
   \begin{definition}  \rm  \label{vierte definition}
            Let for fixed  $ n \in \mathbbm{N}_0$ { \rm and } $\alpha  \in  [ 0, 1 ]$  the set  
            ${ \clubsuit}_{n,\alpha}$  be a subset of  $\Delta_n$,     \\
            ${ \clubsuit}_{n,\alpha}$ :=  \{ $(x_0,  \ldots , x_n ) \ \in  \Delta_n  \ |$ 
            there is  at least one  $j \in  \{0, \ldots , n\}$  
            such that $ x_j \ = \ \alpha $ \} .   { $ \  $ }    \hfill $\Box$   
     \end{definition}
      ${ \clubsuit}_{n,\alpha}$ is called the  `$ \alpha-Cross \ of \ \Delta_n$' .      
      \\  \\    Examples : 
         \begin{figure}[ht]
                 \centering
                 \setlength{\unitlength}{1.0cm}
             \begin{picture}(6,5)
                  \put(-4.0,3.0){${ \clubsuit}_{2,\frac{1}{6}}$ : } 
                  \dashline{0.2}(-4,1)(2,1)     \dashline{0.2}(-4,1)(-1,5.5)   \dashline{0.2}(-1,5.5)(2,1)  
                  \put(-4.4,0.7){$e_0$}  \put(2.1,0.7){$e_1$}  \put(-1.2,5.7){$e_2$}
                  \thicklines  \put(-3.4,1.9){\line(1,0){4.8}}
                  \put(-3,1){\line(2,3){2.5}}     \put(1,1){\line(-2,3){2.5}}   
                  \thinlines
                  \put(4.0,3.0){${ \clubsuit}_{2,\frac{5}{6}}$ : } 
                  \dashline{0.2}(4,1)(10,1)     \dashline{0.2}(4,1)(7,5.5)   \dashline{0.2}(7,5.5)(10,1)  
                  \put(3.6,0.7){$e_0$} \put(10.1,0.7){$e_1$}  \put(6.8,5.7){$e_2$}
                  \thicklines       \put(9.4,1.9){\line(-2,-3){0.6}}           
                  \put(4.6,1.9){\line(2,-3) {0.6}}        \put(6.4,4.6){\line(1,0){1.2}}       
            \end{picture}     \\     \caption{}\label{picture5}   
         \end{figure}	 \\  
                                  
    Note that  $ { \clubsuit}_{n,1} = \{  e_0, e_1, \ldots , e_n \}$, the set of the vertices of $ \Delta_n$,
            and for positive $n$   that  
            $ { \clubsuit}_{n,0} $  is the topological boundary of $ \Delta_n$.  
    \begin{lemma}    \label{ein weiteres lemma}  
            We have for \ $ 0 < \alpha , \beta <  \frac{1}{n+1} $ \ that the  $ \alpha $-cross 
            $ { \clubsuit}_{n,\alpha} $ is   homeomorphic to the $ \beta $-cross $ { \clubsuit}_{n,\beta}$, i.e.
            \qquad \qquad \qquad   \qquad \qquad \qquad  \qquad
            $ {{\clubsuit}_{n, \alpha}} \stackrel{\cong}{\longrightarrow}  {{ \clubsuit}_{n, \beta }} $ .      
    \end{lemma} 
    \begin{proof}
            To prove the lemma we need any increasing
            homeomorphism  $ f $  on the  interval  $  \left[ 0, \frac{1}{n+1} \right]$  with  $ f(\alpha) = \beta $.  
            For instance take the polygon through three points  
            $ \left\{ (0,0) , ( \alpha, \beta ) ,  \left(\frac{1}{n+1},\frac{1}{n+1} \right) \right\}$. \ Then we
            use  Proposition  \ref{proposition vier}, with which we shall deal a few pages later,  and  the 
            restriction  of $ { \bf \Lambda_n}(f)$ \ to the subset \ $ { \clubsuit}_{n,\alpha} $  of   $ \Delta_n$ \ 
            yields the desired homeomorphism.  Note Remark  \ref{remark fuenf}.	 
    \end{proof}
    \begin{definition}
          {  \rm  Let for all } $ n \in   \mathbbm{N}$ \\ 
          $ {\bf BOU}_{n} := \ \{  (x_0,  \ldots , x_n ) \in  \Delta_n \ |$ \rm there is at least one 
          $ j \in \{ 0,1,2,  \ldots , n \} $  such that $ x_{j} = 0 $ \}    \hfill $\Box$ 
    \end{definition}   
          $ { \bf BOU}_n $  is called  the  `$ Boundary \ of \ \Delta_n$',  because 
          $ { \bf BOU}_n $  is the topological  boundary of  $ \Delta_n  \subset  \mathbbm{R}^{n+1}$.   
          Note that  $ { \bf BOU}_n $  is homeomorphic to the  $(n-1)$-sphere, and that  
          $ { \bf BOU}_n = { \clubsuit}_{n,0} \, $.      
   \begin{definition}   \quad   \rm 
          {\rm  Let for all} $n \in \mathbbm{N}$ {\rm and for an} \ $\alpha \in \left[ 0, \frac{1}{n+1} \right]$,
           $  {\bf Layer}_{n,\alpha} \ \subset { \clubsuit}_{n,\alpha}$, \\
          $ { \bf Layer}_{n,\alpha} := \ \{ (x_0, \  \ldots \ , x_n ) \ \in  \Delta_n \ |$ \
          $ \alpha = \min \{ x_0, x_1 , \ \ldots \ x_n \}  \} $   \hfill $\Box$ 
    \end{definition}	
    $  { \bf Layer}_{n,\alpha} $ is called the `$ \alpha-Layer \ of \ \Delta_n  $'.   \\
          Note that   $ \ { \bf Layer}_{n,\frac{1}{n+1}} =  { \{ { \bf Center}_n \} }$  
          and   $ { \bf Layer}_{n,0}  = { \bf BOU}_{n} = { \clubsuit}_{n,0} $ ,  and that for  
          $ 0 \leq \alpha < \frac{1}{n+1} $  we have  that  $ { \bf Layer}_{n,\alpha} $  is homeomorphic to the 
          $(n-1)$-sphere.   There is a homeomorphism 
          $$ \vec{x} \  \stackrel{}{\longmapsto} \ \vec{x} - { \bf Center}_n \ \stackrel{}{\longmapsto} \ 
                      \frac{1}{ \parallel\vec{x} - { \bf Center}_n \parallel}
                      \cdot ( \vec{x} - { \bf Center}_n ) , \ \quad 
                      {\rm for \ all \ } \  \vec{x} \in {\bf Layer}_{n,\alpha} \, . $$    
          We get a disjoint union 
          $ \Delta_n =  \bigcup \ \{ { \bf Layer}_{n,\alpha} \ | \ 0 \leq \alpha \leq \frac{1}{n+1} \}$. 
    \begin{definition}   \quad   \rm 
          Let  for  $n \in  \mathbbm{N}$  and for fixed 
          $j \in \{ 0, 1, 2, \ldots, n \} \ {\bf Section}_{n,j}$ be a subset of $\Delta_{n}$, \\
      \centerline{     ${\bf Section}_{n,j} := \{ (x_0, \ldots , x_n ) \in  \Delta_n \ |$ 
                          $\min\{x_0, \ldots , x_n \} = x_j \}$. }
   \end{definition}
       $ { }  $     $ { \bf Section}_{n,j}$  is called the  {`\it $j$-Section of $ \Delta_n$'}.   \\
      A subset of  $ \Delta_n $ of the form  $ { \bf BOU}_{n} \cap { \bf Section}_{n,j} $ is called a 
       { \it face } of   $ \Delta_n $,  for $ j = 0,1, \ldots, n$ \  we have the face  \quad  
       $ { \bf BOU}_{n} \cap  { \bf Section}_{n,j}  =  \{ (x_0, \ldots , x_n ) \in  \Delta_n \ | 
          \min\{x_0, \ldots , x_n \} = x_j = 0 \}$.      \\  	 
    
       Note  the union $ \Delta_n  = \bigcup \ \{ { \bf Section}_{n,j} \ | \  j = 0, 1, 2, \ \ldots \ , n \} $
          which is not disjoint, e.g. we have for all 
          $ j \in \{0,1, \ldots, n \}$ that ${\bf Center}_n $  is an element of ${ \bf Section}_{n,j}$ .    \\
                                         
       Now, for  natural numbers  $ n > 0$  we shall project  every 
           $\vec{x} \in \Delta_n, \ \vec{x} \neq { \bf Center}_n$ onto the $\alpha$-Layer  of $ \Delta_n$, 
            for all $0 \leq \alpha \leq \frac{1}{n+1}$ .    \\
       We shall define $ \pi := \pi_0 : \ \Delta_n\backslash \{ { \bf Center}_n \}  \longrightarrow  
          { \bf BOU}_n = { \bf Layer}_{n,0}$,  and  for $0 < \alpha < \frac{1}{n+1}$  we shall define maps
          $\pi_{\alpha}: \Delta_n\backslash \{{ \bf Center}_n \}  \longrightarrow  {\bf Layer}_{n,\alpha}$, 
          and let $\pi_{\frac{1}{n+1}}: \Delta_n  \longrightarrow \{{ \bf Center}_n \}$  be the constant map. 
                                     
      In details:   Let  $ n \in  \mathbbm{N}$. 
       For an arbitrary  $\vec{x} \in \Delta_n\backslash \{{ \bf Center}_n \} ,  \vec{x} =
        (x_0, x_1, x_2, \ \ldots , x_n)$,  let  $ \vartheta$  be a permutation on   $ \{0,1,2, \ldots , n \}$ 
         such that  \
         $ 0 \leq x_{\vartheta(0)} \leq  x_{\vartheta(1)} \leq x_{\vartheta(2)} \leq \ldots  \leq x_{\vartheta(n)}$.
            
      Then define  $\pi(\vec{x}) := \vec{b} := ( b_0, b_1, b_2, \ \ldots \ , b_n )$,  
      we set for all \ $i \in \{ 0,1,2, \ldots ,n\}$ 
      $$ b_{i} := \frac{1}{1 - (n+1) \cdot x_{\vartheta(0)}} \cdot  \left( x_i - x_{\vartheta(0)}\right) \, ,$$ 
            hence  $ b_{\vartheta(0)} = 0 $  and  $\pi(\vec{x}) \in   { \bf BOU_n}$. \\ 
            For  $ 0 \leq \alpha \leq \frac{1}{n+1} $ ,   let  $\pi_{\alpha}(\vec{x}) 
            := (y_0, y_1, y_2, \ \ldots \ , y_n)$,   for  $ i \in \{ 0, 1, 2, \ \ldots , \ n \}$    we define 
      $$ y_{i} \ := \ \alpha + ( 1 - (n+1)\cdot \alpha ) \cdot b_i \ = \  \alpha +  \frac{ 1 - (n+1) \cdot 
      \alpha }{1 - (n+1) \cdot x_{\vartheta(0)}} \cdot  \left( x_i - x_{\vartheta(0)}\right) \, . $$ 
      We get  $  y_{\vartheta(0)} = \alpha $,  and for all  $ i \in \{0, \ldots , n\}$ we have  
            $ y_i \ \geq \ \alpha $, hence \  $\pi_{\alpha}(\vec{x}) \in   {\bf Layer}_{n,\alpha}$ .  \\
   Consider also the continuous surjective map ${\cal A}: \Delta_n \longrightarrow \left[ 0, \frac{1}{n+1} \right]$  
            for all  $ \vec{x} \in  \Delta_n  $. \ There is an unique number  
            ${\cal A}(\vec{x}) := x_{\vartheta(0)} = \min \{ x_0, x_1 , \ \ldots , \ x_n \} $, such that for 
            $  \vec{x} \neq { \bf Center}_n $ we have
   $$ \vec{x} \ = \ \pi(\vec{x}) +  {\cal A}(\vec{x}) \cdot (n+1) \cdot ({ \bf Center}_n - \pi(\vec{x}) ) \ , \ \  
            { \rm  and } \ \ \ {\cal A}({ \bf Center}_n ) = \frac{1}{n+1} \ . $$
   For $\vec{x} \in \Delta_n $,  we have  
            $\vec{x} = \pi_{{\cal A}(\vec{x})}(\vec{x}) \in { \bf Layer}_{n,{\cal A}(\vec{x})} $,  
            and for  ${\cal A}(\vec{x}) \neq \alpha$, and \ $0 < \alpha , {\cal A}(\vec{x})  < \frac{1}{n+1}$, 
            we have four collinear points 
            \{$ {\bf Center}_n , \, \vec{x}, \, \pi_{\alpha}(\vec{x}), \, \pi(\vec{x}) $\}.   \\
   Note that if we restrict  ${\cal A}$  to an  $ \alpha $-layer, for all $ \ 0 \leq \alpha \leq \frac{1}{n+1}$, 
            ${\cal A}$ \ will be a constant map, 
            ${\cal A}( \vec{x}) = \alpha $  for all  $ \vec{x} \in { \bf Layer}_{n,\alpha}$, \ because \ 
            $ \min \{x_0, \ldots x_n\} = \alpha$.    
  	\begin{definition}  \rm \quad  Let  $n \in  \mathbbm{N}$ and $ \vec{a}, \ \vec{b} \in \mathbbm{R}^{n}$.                      Then we define  $ \left[\vec{a}, \vec{b} \right] \subset \mathbbm{R}^{n}$. Let \\
	       $ \left[\vec{a}, \vec{b} \right] \ := \ \left\{ \ t \cdot \vec{a} + (1-t) \cdot \vec{b} \ | \
	        t \in [0,1] \ \right\}$ be the straight line confined by  $\vec{a}$  and  $\vec{b}$.  \hfill $\Box$ 
  	\end{definition}
	Note that for any standard \ $n$-simplex $\Delta_n$ \ we have the (nearly disjoint) union \\
	       $ \Delta_n \ = \  \bigcup \ \left\{ \ \left[ { \bf Center}_n  , \vec{b} \right] \ | \ 
	       \vec{b} \in  { \bf BOU}_n \ \right\} $,   all the lines intersect only in \ $ { \bf Center}_n$.   \\  
	Further note that for all  $ \vec{b} \in  { \bf BOU}_n$,   we get a constant projection  
      $$ \pi|_{ \left[ { \bf Center}_n, \vec{b} \right]\backslash \{{ \bf Center}_n \}} \, :  \ \ 
        \left[ { \bf Center}_n, \vec{b} \right]\backslash \{{ \bf Center}_n \} \ \longrightarrow  \ \vec{b} \, .$$
    For a point  $ \vec{x} = ( x_0 , \ldots , x_n ) \in  \left[ { \bf Center}_n , \vec{b} \right] $, 
     $ \vec{x} \neq { \bf Center}_n $,
    i.e.  $ \pi( \vec{x}) =  \vec{b}$,   we have the unique number $ t := {\cal A}(\vec{x}) \cdot (n+1) \in [0,1]$ 
    such that  $ \vec{x} = t \cdot { \bf Center}_n + (1-t) \cdot \vec{b}$. \
    Let  $ \vec{b} = ( b_0,  \ldots , b_n )$. \ For a component   $ x_j$  we get the notation \ 
    $ x_j =  t \cdot \frac{1}{n+1} + \left(1-t\right) \cdot b_j  = t \cdot \left( \frac{1}{n+1} - b_j \right) + b_j$.     \begin{definition}  \rm \quad  Let  $n \in  \mathbbm{N}$. \ For any subset $ M \subset  \mathbbm{R}^{n}$  let \\ 
        $ { \sf Sponge}(M) := \{ ( x_1, x_2, \ldots , x_n ) \in M \, | \, x_i = x_j \text{ if and only if } i=j \}$.  
                                              \hfill  $\Box$   
  \end{definition} 
	      That means  that ${\sf Sponge}(M)$  contains only $n$-tuples $ (x_1,x_2, \ldots, x_n) $ with 
	      pairwise different components.  
	  \begin{remark}   \rm    \label{remark drei}
	      A continuous map $ f: \Delta_n \rightarrow \Delta_n$ is uniquely determined by the induced map 
	      $ \widetilde{f}:  {\sf Sponge}(\Delta_n) \rightarrow \Delta_n$, \ with \ 
	      $  \widetilde{f} := f_{| {\sf Sponge}(\Delta_n)} $.     
	  \end{remark}
	  \begin{proof} This is trivial since \ closure$({\sf Sponge}(\Delta_n)) = \Delta_n $ .  
	  \end{proof}  
	   { $ $   }   
	  Now we define a nice subgroup  of the  group of all homeomorphisms on   $ \Delta_n$ .  		
	  \begin{definition} 	\quad \rm  \label{irgendeine Definition von Abbildungen}
        	For a fixed    $ n \in  \mathbbm{N}_0$  let \  ${ \cal COMFORT}(\Delta_n)  \subset
  	      \{ F: \Delta_n \rightarrow \Delta_n \} 	$.
        	 A map $ F $ is an element of    ${ \cal COMFORT}(\Delta_n) $  if an only if   $F$
        	 fulfils the following conditions   ${ \bf \widehat{[1]}, \widehat{[2]}, \widehat{[3]}}$. \\
   	       ${ \bf  \widehat{[1]}}$: \quad $F$  is  a homeomorphism on $\Delta_n$ . \\
  	       ${ \bf  \widehat{[2]}}$: \quad $F$   respects permutations, that means if
  	       $ \vec{x} =  ( x_{0} , x_{1} ,   \ldots   \ldots ,    x_{n} )  \in  \Delta_n $    and if
  	       $F(\vec{x}) = \vec{y} =  ( y_{0} , y_{1} ,   \ldots    \ldots , y_{n} )  \in  \Delta_n$,   and if
  	       $ \vartheta $   is a permutation  on   $ \{ 0,1,2,   \ldots ,   n \} $, then  \\
  	       \centerline{ $ F\left( x_{\vartheta(0)} , x_{\vartheta(1)} , x_{\vartheta(2)} ,\ \ldots   \ldots ,  \
  	                    x_{\vartheta(n)} \right) =  \left( y_{\vartheta(0)} , y_{\vartheta(1)} , y_{\vartheta(2)} ,\
  	                    \ldots    \ldots ,  \  y_{\vartheta(n)} \right) $.  } \\
           ${ \bf  \widehat{[3]}}$: \quad  $F$   keeps the order. Trivially,  for every
           $ \vec{x} =  ( x_{0} , x_{1} ,   \ldots   \ldots ,    x_{n} )  \in  \Delta_n $
           exists a permutation   $\vartheta$ on   $ \{ 0,1,2,   \ldots ,   n \} $ and a number
           ${ \tt r} \in \{ 0,1,2,   \ldots ,   n \}$ \ 
           (We introduce here the number ${ \tt r}$. It will play a major part not until 
           Proposition \ref{proposition vier}) such that  \\
           \centerline{  $ 0 \leq  x_{\vartheta(0)} \leq  x_{\vartheta(1)} \leq    \ \ldots   \
	                  \leq  x_{\vartheta({ \tt r})} \leq  \frac{1}{n+1} <
	                   x_{\vartheta({ \tt r}+1)} \leq     \  \ldots    \  \leq  x_{\vartheta(n)} \leq  1 $. } \\
		       If     $F(\vec{x}) = \vec{y} = ( y_{0} , y_{1} ,   \ldots    \ldots , y_{n} )$, then we demand that \\
		       \centerline{    $ 0 \leq  y_{\vartheta(0)} \leq  y_{\vartheta(1)} \leq    \ \ldots   \
	                    \leq  y_{\vartheta(j)} \leq   y_{\vartheta(j+1)} \leq     \  \ldots    \
		                  \leq  y_{\vartheta(n)} \leq  1 $  } \\
		       holds for all $ j \in \{ 0, 1, 2, \ldots , n-1 \} $.  
		                          
		  For any subset   $ { \bf S } \subset \Delta_n $   we say that a homeomorphism
		       $ F : {\bf S} \rightarrow {\bf S}$   is an element of  ${ \cal COMFORT}({\bf S}) $
		       if and only if  $F$ fulfils  ${\bf \widehat{[2]}}, {\bf \widehat{[3]}}$.    \hfill $\Box$
		  \end{definition}
		       It follows that for an   $ F \in { \cal COMFORT}(\Delta_n)$  the homeomorphism $F$ yields a homeomorphism
		       on  each   ${ \bf Section}_{n,k}, \ F|_{{ \bf Section}_{n,k}} \in { \cal COMFORT}({ \bf Section}_{n,k})$
		       for   $ k \in \{ 0,1,2, \ldots , n \}$, \ and $F|_{\bf BOU_n}$ is an element of
		       ${ \cal COMFORT}({\bf BOU_n})$.
		
		 For   $ F,G \in   { \cal COMFORT}(\Delta_n)$  it holds that $ F^{-1} \in { \cal COMFORT}(\Delta_n)$ and
		       furthermore   $F \circ G \in { \cal COMFORT}(\Delta_n)$, \  hence 
		       $\left( {\cal COMFORT}(\Delta_n), \circ \right) $   is a subgroup of the group 
  	       of all homeomorphisms on  $\Delta_n$.  
  	                        
     We remark that a homeomorphism \  $ f \in {\cal COMFORT}({\sf Sponge}(\Delta_n))$  can be uniquely
  	       extended to a map  $ F \in {\cal COMFORT}(\Delta_n)$ \ with \ $ f = F_{| {\sf Sponge}(\Delta_n)}$.
  	
  	  Note that each element of the  set
  	       $ \left\{ \pi_{\alpha} \, | \, \alpha \in \left[0, \frac{1}{n+1} \right] \right\}$  of projections fulfils
  	       the  conditions \  $ {\bf \widehat{[2]}} \text{ and }  {\bf \widehat{[3]}} \text{ of the above Definition }
  	       \ref{irgendeine Definition von Abbildungen} $, and all $ \pi_{\alpha} $  except
  	       $ \pi_{\frac{1}{n+1}}$ are elements    of the following monoid of endomorphisms  \\
  	      \centerline {$ \left( \left\{ f:
  	      \Delta_n\backslash \{{\bf Center}_n \} \rightarrow  \Delta_n\backslash \{{\bf Center}_n \}   \ | \
  	            f \text{ fulfils the conditions } {\bf \widehat{[2]}} \text{ and } {\bf \widehat{[3]}}
  	            \text{ of Definition } \ref{irgendeine Definition von Abbildungen} \right\}, \circ  \right) $.} \\ \\
  	         	
    The next few lemmas deal with the behaviour of homeomorphisms  $\Phi \in {\cal COMFORT}(\Delta_n)$.
    \begin{lemma} \quad  \label{Lemma zwei}
        Let  $ \Phi $   be any homeomorphism on   $ \Delta_n$ . Then  \
        $  \Phi|_{{ \bf BOU}_n} : { \bf BOU}_n  \stackrel{\cong}{\rightarrow} { \bf BOU}_n $ .
    \end{lemma}
    \begin{proof} \quad This is trivial, because ${ \bf BOU}_n$ is the topological boundary of  $ \Delta_n$ .
    \end{proof}
    The following three lemmas describe the fact that a map $ \Phi \in { \cal COMFORT}(\Delta_n)$ preserves
    equalities and inequalities of the components of an   $ \vec{x} \in  \Delta_n $   .
   \begin{lemma} \label{Lemma drei}  \quad
        Let   $ \Phi $   be an element of    $ { \cal COMFORT}(\Delta_n)$.  Trivially,
        for every  $ \vec{b} \in { \bf BOU}_n $ with components $ \vec{b} = ( b_0, b_1,   \ldots   , b_n ) $ \
        there is a permutation   $ \vartheta $   on   $ \{ 0,1,   \ldots   , n \}$  and there are two natural
         numbers \    ${ \tt q} , { \tt r}$,   with   $ 0 \leq { \tt q} \leq  { \tt r} \leq n-1 $ \ such  that  \
        $ \Phi (\vec{b}) =: \vec{c} =: ( c_0, c_1,   \ldots   , c_n ) $,  and
        $$ 0 = b_{\vartheta(0)} = b_{\vartheta(1)} =     \ldots    =  b_{\vartheta({\tt q})}  = 0
                        < b_{\vartheta({\tt q}+1)} \leq    \ldots   \leq  b_{\vartheta({\tt r})} \leq \frac{1}{n+1}
                        <  b_{\vartheta({ \tt r}+1)} \leq   \ldots   \leq b_{\vartheta(n)} \leq 1   \ . $$
        (The case   $ { \tt q} = { \tt r} $   is possible).
          Because \   $ \Phi  \in   { \cal COMFORT}(\Delta_n)$ we have  $ \vec{c} \in { \bf BOU}_n$  and \\
         \centerline{  $ 0 = c_{\vartheta(0)} \leq c_{\vartheta(1)} \leq   \ldots   \leq
                       c_{\vartheta({\tt q})}  \leq  c_{\vartheta({\tt q}+1)} \leq    \ldots   \leq
                       c_{\vartheta({\tt r})} \leq  c_{\vartheta({ \tt r}+1)} \leq   \ \ldots   \
                        \leq c_{\vartheta(n)} \leq 1  \ . $   }  \\
        Then the claim of this lemma is
           $$    c_{\vartheta(0)} = c_{\vartheta(1)} =  \ldots  = c_{\vartheta({\tt q})}  = 0
                <  c_{\vartheta({\tt q}+1)}  \ .  $$
   \end{lemma}
   \begin{proof}  \quad  
         $ \Phi \in { \cal COMFORT}(\Delta_n)$ means that   $ \Phi$   keeps the order. \
         Because of   $ \Phi(\vec{b}) =  \vec{c} $, and because of
         $  0 =  b_{\vartheta({\tt q})} \leq b_{\vartheta(0)} = 0 $, it  follows that
         $  c_{\vartheta({\tt q})} \leq c_{\vartheta(0)} = 0 $,   hence   $  c_{\vartheta({\tt q})} =0 $.  \\
         Now assume that  $ c_{\vartheta({\tt q}+1)} = 0 $.  As we described above,  $\Phi^{-1}$   is an
         element of   $ { \cal COMFORT}(\Delta_n). $
         We have    $ \Phi^{-1}( \vec{c} ) =  \vec{b}$,    and \
         $ c_{\vartheta({\tt q}+1)} = 0 \leq c_{\vartheta(0)} = 0 $,   and we get a contradiction to
         $ b_{\vartheta({\tt q}+1)} >  b_{\vartheta(0)} = 0 $.   Hence the only possibility is
         $ c_{\vartheta({\tt q}+1)} > 0 $.
   \end{proof}
    \begin{lemma} \quad  \label{wahrscheinlich lemma vier}
        Let   $ \Phi $   be an element of    $ { \cal COMFORT}(\Delta_n)$.  Trivially,
        for  $ \vec{x} = ( x_0, x_1,   \ldots   , x_n ) \in  \Delta_n $
        there is a permutation   $ \vartheta $   on   $ \{ 0,1, \ldots , n \}$, and there are two natural numbers
        ${ \tt q} , { \tt r}$,   with   $ 0 \leq { \tt q} \leq  { \tt r} \leq n $ \ such  that
        $ \Phi (\vec{x}) =: \vec{y} =: ( y_0, y_1,   \ldots   , y_n ) \in \Delta_n$,  and
        $$ x_{\vartheta(0)} = x_{\vartheta(1)} =  \ldots   =
         x_{\vartheta({\tt q})}   <   x_{\vartheta({\tt q}+1)} \leq    \ldots   \leq  x_{\vartheta({\tt r})}
        \leq \frac{1}{n+1} <  x_{\vartheta({ \tt r}+1)} \leq   \ldots   \leq x_{\vartheta(n)}   \ . $$
          $({ \tt q} = { \tt r} $   is possible.   If   $ \vec{x}  = { \bf Center}_n $   then \
         $ { \tt q} = { \tt r} = n )$.   It holds \  $ \Phi( { \bf Center}_n) =  { \bf Center}_n \, .$
         If    $ \vec{x} \neq { \bf Center}_n $, then  because
         $ \Phi  \in   { \cal COMFORT}(\Delta_n) $  we get \\
         \centerline{  $ y_{\vartheta(0)} \leq y_{\vartheta(1)} \leq   \ldots   \leq y_{\vartheta({\tt q})}
                       \leq  y_{\vartheta({\tt q}+1)} \leq    \ldots   \leq y_{\vartheta({\tt r})}
                       \leq   y_{\vartheta({ \tt r}+1)} \leq   \ \ldots   \ \leq y_{\vartheta(n)}   \ . $ } \\
         Then we claim that
     $$ y_{\vartheta(0)} = y_{\vartheta(1)} =   \ldots   = y_{\vartheta({\tt q})} <  y_{\vartheta({\tt q}+1)} \ . $$
   \end{lemma}
   \begin{proof} \quad Follow the lines of the previous proof.
   \end{proof}
   \begin{lemma} \quad   \label{Lemma fuenf}
        Let   $ \Phi $   be an element of    $ { \cal COMFORT}(\Delta_n)$.   Let
        $ \vec{x} \in \Delta_n\backslash \{ { \bf Center}_n \}  , \\
        \vec{x} = ( x_0, x_1,   \ldots   , x_n ). $   Let \
        $ \Phi (\vec{x}) =: \vec{y} =: ( y_0, y_1,   \ldots   , y_n ) \in \Delta_n$.
        There is a permutation   $ \vartheta $   on   $ \{ 0,1,   \ldots   , n \}$,  and there are
        two natural  numbers \ ${ \tt q}, { \tt t}$, with  $0 \leq { \tt q}  \leq n , \ 1 \leq
        { \tt t} \leq n - {\tt q}$, \ and
    $$ x_{\vartheta(0)} \leq x_{\vartheta(1)} \leq    \ldots    \leq
        x_{\vartheta({\tt q})} < x_{\vartheta({\tt q}+1)} =  x_{\vartheta({\tt q}+2)}\ \ldots   \ldots
        = x_{\vartheta({\tt q}+  { \tt t} )} <  x_{\vartheta({\tt q}+  { \tt t} +1  )} \leq \ldots   \ldots   \leq
                                                     x_{\vartheta(n)} \, . $$
        $( \text{The \ case \ }  { \tt t}= n -{ \tt q} \text{ is possible} )$.   \;
         Because  $ \Phi  \in {\cal COMFORT}(\Delta_n) $   we have \\
         $ y_{\vartheta(0)} \leq y_{\vartheta(1)} \leq    \ldots    \leq
         y_{\vartheta({\tt q})} \leq y_{\vartheta({\tt q}+1)} \leq  y_{\vartheta({\tt q}+2)}\ \leq   \ldots   \leq
         y_{\vartheta({\tt q}+  { \tt t} )} \leq  y_{\vartheta({\tt q}+  { \tt t} +1  )}   \leq   \ldots    \leq \
         y_{\vartheta(n)}   \ . $  \\
         We claim that
         $$ y_{\vartheta({\tt q})}   <   y_{\vartheta({\tt q}+1)} =  y_{\vartheta({\tt q}+2)}\   =    \ldots   \ldots \
         =   y_{\vartheta({\tt q}+  { \tt t} )}   <   y_{\vartheta({\tt q}+  { \tt t} +1  )} \, .  $$
   \end{lemma}
   \begin{proof} \quad Essentially it is the same proof as before twice.
   \end{proof}
   \begin{corollary} \label{corollary eins}
     It follows from the previous lemma that an element $ F \in {\cal COMFORT}(\Delta_n)$
  	       yields an   $ F_{| {\sf Sponge}(\Delta_n)}  \in {\cal COMFORT}({\sf Sponge}(\Delta_n))$.
   \end{corollary}	     { \ $ $  }   \\
   We continue the investigations of maps $ \Phi  \in   { \cal COMFORT}(\Delta_n)$ with an important statement.
   Note that in the  following  Lemma  \ref{propositionlemma eins}  we assume a homeomorphism  \
         $ { \clubsuit}_{n,\alpha}    \stackrel{\cong}{\longrightarrow} {\clubsuit}_{n,\beta}$ \ for some
         $ \alpha $ and $ \beta $, which is not always possible, e.g. for
         $ \alpha > 0 $ and $ \beta = 0 $.
         We believe, but we have no proof that a homeomorphism
          $ { \clubsuit}_{n,\alpha}    \stackrel{\cong}{\longrightarrow} {\clubsuit}_{n,\beta}$ \
	       is   possible  in the following cases of  $\alpha$ and $\beta$: \\
	       {  $ {  }  $  }  \qquad   \qquad    \ $ \frac{1}{k+1} < \alpha , \beta  < \frac{1}{k} $, \
                                        for all natural numbers \ $ 1 \leq k \leq n $.  \\
         Further there exists a homeomorphism   
         $ { \clubsuit}_{n,\alpha}    \stackrel{\cong}{\longrightarrow} {\clubsuit}_{n,\beta}$ \  if 
         $\alpha = \beta $, of course, and also  if  \\        
         {  $ {  }  $  }   \qquad   \qquad    \ $ 0 < \alpha , \beta  < \frac{1}{n+1} \, $. \ 
         This fact is proved in Lemma  \ref{ein weiteres lemma} .  
   \begin{lemma}  \label{propositionlemma eins}
           	\quad Let   $n \in  \mathbbm{N}$. Let $0 \leq \alpha , \beta < \frac{1}{n+1} $, and let
            $ \Phi  \in   { \cal COMFORT}(\Delta_n) $, and we assume  that $\Phi $  induces a homeomorphism \
	          $ \Phi |_{ {\clubsuit}_{n,\alpha} }:  \
	          { \clubsuit}_{n,\alpha}    \stackrel{\cong}{\longrightarrow} {\clubsuit}_{n,\beta}$.  \\
	          Let  $ \vec{x} = ( x_0, x_1, \ldots , x_n ) \in  {\clubsuit}_{n,\alpha} \cap {\sf Sponge}(\Delta_n)$,
	          i.e. we have  a single $\alpha $  in the components \
	          $ \{ x_0, x_1, \ldots , x_n \} \ \text{of} \  \vec{x}$. \ Trivially, there is a permutation
	          $ \vartheta $  on  $ \{ 0,1, \ldots, n \}$, and  a natural  number  ${ \tt q}$, with
	          $0 \leq { \tt q} \leq n-1 $ \ such that the  $ \vartheta ({\tt q})^{th}$ component of
	          $ \vec{x} \ \text{is  the single } \alpha $,  i.e. $ \alpha = x_{\vartheta({\tt q})}$, and
    $$ x_{\vartheta(0)} < x_{\vartheta(1)}  < \ldots <   x_{\vartheta({\tt q}-1)} <  \alpha  < x_{\vartheta({\tt q}+1)}
                        < x_{\vartheta({\tt q}+2)}  <  \ldots \ldots < x_{\vartheta(n)} \ . $$
      Let    $ \Phi (\vec{x}) =: \vec{y} =: ( y_0, y_1, \ldots , y_n)$.
         We claim that the  $ \vartheta ({\tt q})^{th}$ component of $ \vec{y} $ is  $\beta$,
         i.e.   $ \beta  = y_{\vartheta({\tt q})}$, and this is the only $\beta$  in the components \
         $ \{ y_0, y_1, \ldots , y_n \} \ \text{of} \  \vec{y} $, \ i.e. we claim
     $$ y_{\vartheta(0)} < y_{\vartheta(1)} < \ldots  <  y_{\vartheta({\tt q}-1)} <  \beta <
          y_{\vartheta({\tt q}+1)} <  \ldots \ldots <  y_{\vartheta(n)}   \ . $$
   \end{lemma}
   \begin{proof} \quad    The case $n=1$ is trivial. We have \
       ${\clubsuit}_{1,\alpha} = \{ (\alpha, 1- \alpha), (1-\alpha, \alpha) \} $, \  hence \
       $ \Phi(\alpha, 1 - \alpha) =  (\beta, 1 - \beta). $  \\ 
       Let $n \geq 2$. \\ 
       Since  \ $ \Phi |_{ {\clubsuit}_{n,\alpha} }:
	          { \clubsuit}_{n,\alpha}    \stackrel{\cong}{\longrightarrow} {\clubsuit}_{n,\beta}$ \
       we can  assume that either \ $\alpha = \beta = 0$ \ or \ $ 0 < \alpha, \beta < \frac{1}{n+1} \, . $  \\
      For  $\alpha = \beta = 0$  we apply Lemma \ref{Lemma drei}, with  ${ \tt q} = 0 $.  \ \
      We assume    $ \alpha > 0$ and $ \beta > 0$.

      Because  $\vec{x} \in {\clubsuit}_{n,\alpha} \cap {\sf Sponge}(\Delta_n)$
       and Corollary \ref{corollary eins}  the   components  $\{ y_0, y_1,   \ldots   , y_n \} $ of
       $ \vec{y}$ contain a single  $\beta$.
       Since  $ \Phi \in {\cal COMFORT}(\Delta_n)$, we can deduce from the inequality \\
   \centerline {   $ x_{\vartheta({\tt q}-1)} < \alpha =  x_{\vartheta({\tt q})} < x_{\vartheta({\tt q}+1)}$  } \\
       and Corollary \ref{corollary eins}   the  inequality \\
   \centerline{  $ y_{\vartheta({\tt q}-1)} <  y_{\vartheta({\tt q})}  < y_{\vartheta({\tt q}+1)}$.  }  \\
       (The case $ {\tt q} = 0 $ is posssible.)  We want to show that $ y_{\vartheta({\tt q})} = \beta  $.
       We call   `$\vartheta({ \tt k})$' the index of $\beta$, i.e. $\beta = y_{\vartheta({\tt k})},$
        for a suitable   ${ \tt k} \in  \{ 0,1, \ldots , n-1 \}$,  and  we have  \\
        \centerline{  $ y_{\vartheta(0)} < y_{\vartheta(1)} < \ldots  < y_{\vartheta({\tt k}-1)} <
                        \beta  < y_{\vartheta({\tt k}+1)} < \ldots <  y_{\vartheta(n)} \ . $ }   \\
       We want to show  ${\tt k} = {\tt q}$.  To prove this we consider the two other cases
        ${\tt k} < {\tt q}$ and   ${\tt k} > {\tt q}$ and we seek contradictions.  We shall  find a contradiction in
        the case of  ${\tt k} < {\tt q}$. The  case ${\tt k} > {\tt q}$ can be treated  in the same way.
        (Because $ \Phi^{-1}  \in   { \cal COMFORT}(\Delta_n)$  we can exchange the parts of
       $ \alpha$   and   $ \beta$, note \ $ \Phi^{-1} |_{ {\clubsuit}_{n,\beta}}:
	          { \clubsuit}_{n,\beta}  \stackrel{\cong}{\longrightarrow} { \clubsuit}_{n,\alpha} $).
	
    The case  \  ${\tt k} < {\tt q}$: \\
    It is  $ 0 < {\tt q}$ in this case. Further note that we can exclude the case  \ $ x_{\vartheta(0)} = 0,
    x_{\vartheta(1)} = \alpha$, since then from Lemma \ref{Lemma drei} would follow that \ $ y_{\vartheta(0)} = 0$,
    this contradicts the assumption  ${\tt k} < {\tt q}$.  (Note $ 0 < \beta = y_{\vartheta({\tt k})}$).

   We define an infinite connected subset \ $ Subset[{\clubsuit}_{n,\alpha}] \subset {\clubsuit}_{n,\alpha}$
        to use a topological argument.    Let  for all $ \varepsilon \in [0,1]$ the element
       $\vec{a}_{\varepsilon} := (a_0, a_1, \ldots, a_n) \in {\clubsuit}_{n,\alpha} $ \ by setting
      \[    a_j   :=
                   \begin{cases}   x_j \cdot (1-\varepsilon)   &
                     \text{ for } \ \   j \in  \{ \vartheta(0) , \vartheta(1) , \ldots , \vartheta({\tt q}-1) \} \  \\
                      x_j    &   \text{ for } \ \  j \in
                      \{\vartheta({\tt q}),\vartheta({\tt q}+1),  \ldots , \vartheta(n-1) \}  \\
                    x_{\vartheta(n)} +   \varepsilon \cdot \sum_{i=0}^{{\tt q}-1}   x_{\vartheta(i)}
                     &   \text{ for }  \ \   j =  \vartheta(n)
                   \end{cases}
       \]
   We have $ a_{\vartheta({\tt q})} = \alpha$, hence   $  \vec{a}_{\varepsilon} \in  { \clubsuit}_{n,\alpha}$
            for all $ \varepsilon \in [0,1]$. For $ \varepsilon  \neq 1 $ we have
             $ \vec{a}_{\varepsilon} \in  {\sf Sponge}(\Delta_n)$.   Let
    $$  Subset[{ \clubsuit}_{n,\alpha}] \ :=  \left\{\vec{a}_{\varepsilon} \ | \  \varepsilon \in [0,1] \right\} $$
        We get   $ \vec{a}_{0} = \vec{x}$, and  
       $ \vec{a}_{1} =:  ( \bar{a}_0, \bar{a}_1,   \ldots   , \bar{a}_n ) \in { \bf BOU}_n $, and we have  \\
        \centerline{  $ 0 = \bar{a}_{\vartheta(0)} = \bar{a}_{\vartheta(1)} =   \ldots   =
                      \bar{a}_{\vartheta({\tt q}-1)} = 0  < \alpha < \bar{a}_{\vartheta({\tt q}+1)} <
                      \bar{a}_{\vartheta({\tt q}+2)} <  \ldots \ldots <  \bar{a}_{\vartheta(n)} $ .  }      \\
         Let   $ \Phi(\vec{a}_{\varepsilon} ) =: \vec{d}_{\varepsilon} =: ( d_o, d_1 , \ldots , d_n ) \in
         { \clubsuit}_{n,\beta}\ ,   \text{for   all}   \ \varepsilon \in [0,1]$.    Because of
         $ \Phi  \in   { \cal COMFORT}(\Delta_n)   $ and because of Lemma  \ref{Lemma drei} and
         Lemma \ref{Lemma fuenf}, with \
         $ \Phi(\vec{a}_{1} ) = \vec{d}_{1} =: ( \bar{d}_o, \bar{d}_1 , \ldots ,   \beta , \ldots ,   \bar{d}_n )
         \in  { \clubsuit}_{n,\beta}\ $ \  we get a single component  $\beta$ in $\vec{d}_{1}$, i.e. there is an index
         $ {\tt j }$ with ${\tt q} \leq {\tt j}$  such that  $ \beta =  \bar{d}_{\vartheta( { \tt j})}$, \ and
     $$ 0 = \bar{d}_{\vartheta(0)} =  \ldots  = \bar{d}_{\vartheta({\tt q}-1)} = 0  <
         \bar{d}_{\vartheta({\tt q})} < \bar{d}_{\vartheta({\tt q}+1)} <  \ldots < \bar{d}_{\vartheta({\tt j}-1)}
         <   \beta  <  \bar{d}_{\vartheta({\tt j}+1)}\ <  \ldots  < \bar{d}_{\vartheta(n)}  \,  .    $$
         Now we use the canonical projections  ${ \sf PROJ_{   \it i } } $,  for   $ i = 0,1, \ldots , n $,
        $$   {\sf PROJ_{   \it i } }:   \Delta_n \rightarrow [0,1] \, , \
                              ( z_0 , z_1 , z_2 , \ldots , z_n ) \mapsto z_i  \,  .  $$
      ${ \sf PROJ_{   \it i } }$ \  is continuous.    \
         By definition, the components of each  $\vec{a}_{\varepsilon} \in Subset[{ \clubsuit}_{n,\alpha}]$
         contain a single \ $ \alpha $ at the $ \vartheta({\tt q})^{th}$  place.  By Lemma \ref{Lemma fuenf}  \
         the image     $ \Phi( \vec{a}_{\varepsilon}) =  d_{\varepsilon}$  contains a single $ \beta$ for all
         $ \varepsilon \in [0,1]$.
         Note that both    $ Subset[{\clubsuit}_{n,\alpha}]$   and its homeomorphic image \
         $ \Phi(Subset[{\clubsuit}_{n,\alpha}])$ \
         are connected subsets of $ \Delta_n$, and also note
         that the set  \ $ \{ \beta \} $  is closed in  $ [0,1]$.   If we restrict   ${ \sf PROJ_{   \it i } } $ \
         for  \ $i = 0,1, \ldots , n$ \ to the set \ 
         $ \Phi(Subset[{\clubsuit}_{n,\alpha}]) \subset {\clubsuit}_{n,\beta}$, 
      $$  { \sf PROJ_{   \it i } }:  \Phi(Subset[{\clubsuit}_{n,\alpha}])  \longrightarrow  [ 0,1 ] \, , $$   
         we can express \ $ \Phi(Subset[{\clubsuit}_{n,\alpha}])$ \ as a disjoint union of domains
         ${ \sf PROJ_{   \it i }}^{-1} (\{ \beta \})$, \ that means
     $$ \Phi(Subset[{\clubsuit}_{n,\alpha}]) \ =  \
         \bigcup  \left\{ {\sf PROJ_{ \it i }}^{-1} (\{ \beta \}) \, | \,  i = 0,1, \ldots , n \right\} \, , $$
         and the union is disjoint.   \\
         Because   $ \Phi(Subset[{\clubsuit}_{n,\alpha}])$   \ is connected and
         ${ \sf PROJ_{   \it i }}^{-1} (\{ \beta \})  $   is closed, we see that
         ${ \sf PROJ_{   \it i }}^{-1} (\{ \beta \})  $   either is empty or the entire set, for each \
         $i \in \{ 0,1, \ldots, n \}$.  We have  
     $$  { \sf PROJ_{\vartheta({\tt k})}}^{-1} (\{ \beta \}) \neq \emptyset
                                  \neq { \sf PROJ_{\vartheta({\tt j})}}^{-1} (\{ \beta \})   \, .  $$
     Hence  it follows  \ $ {\tt k} = {\tt j} $, which contradicts   $ {\tt k} < {\tt q} \leq {\tt j} $ !

     As we already mentioned above, the  case ${\tt k} > {\tt q}$ can be treated  in the same way
          by exchanging the parts of $\alpha$   and   $ \beta$ and considering  $ \Phi^{-1}$ instead of
          $\Phi$. Note that  $ \Phi^{-1}$ induces a homeomorphism
	        ${ \clubsuit}_{n,\beta} \stackrel{\cong}{\longrightarrow} { \clubsuit}_{n,\alpha} $.   \
	        This finishes the proof of  Lemma  \ref{propositionlemma eins}.
  \end{proof}
   From the previous lemma we can deduce an important corollary.
  \begin{corollary}  \label{corolloary zwei}
  	\quad Let   $n \in  \mathbbm{N}$. Let $0 \leq \alpha , \beta < \frac{1}{n+1} $, let
            $ \Phi  \in   { \cal COMFORT}(\Delta_n) $, and assume  that $\Phi $  yields a homeomorphism \
	          $ \Phi |_{ {\clubsuit}_{n,\alpha} }:
	          { \clubsuit}_{n,\alpha}    \stackrel{\cong}{\longrightarrow} {\clubsuit}_{n,\beta}$ .
	          For a point  $ \vec{x} = ( x_0, x_1, \ldots , x_n ) \in  \Delta_n $ let the image  be \
	          $ \Phi(\vec{x}) = \vec{y} = ( y_0, y_1, \ldots , y_n )$.  \\
	          Then \ $ x_i = \alpha$ \ if and only if \ $ y_i = \beta$ \ for all indices $ i = 0,1, \ldots, n $.
  \end{corollary}
  \begin{proof}
  Use the previous Lemma \ref{propositionlemma eins} and note that closure$({\sf Sponge}(\Delta_n)) = \Delta_n$ .
  \end{proof}
    \begin{lemma}  \label{lemma sieben}
           	\quad Let   $n \in  \mathbbm{N}$. Let    $ 0 \leq \alpha , \beta < \frac{1}{n+1} $ ,   and let \
            $ \Phi  \in { \cal COMFORT}(\Delta_n) $  such that $ \Phi $   yields a homeomorphism \quad
	          $ \Phi |_{ {\clubsuit}_{n,\alpha} }   :
	          { \clubsuit}_{n,\alpha}    \stackrel{\cong}{\longrightarrow}   { \clubsuit}_{n,\beta} $.
            Let    $ \vec{x} = ( x_0, x_1,   \ldots   , x_n ) \in  \Delta_n\backslash {\clubsuit}_{n,\alpha}$,
            i.e.  $ \vec{x}$ is not an element of the $ \alpha$-cross.  \
            Trivially, either \  $ \alpha < \min \{ x_0, x_1, \ldots , x_n \}$ \
            or    there is a permutation
            $ \vartheta $   on   $ \{ 0,1, \ldots , n \} $ \ and there are two natural  numbers \
        ${ \tt q} ,  { \tt r}$, with  \\
        \centerline{ $ 0 \leq { \tt q} \leq { \tt r} < n \ ({ \tt q} = {\tt r} \ \text{is possible})$, and }  
    $$ x_{\vartheta(0)} \leq x_{\vartheta(1)} \leq    \ldots    \leq
               x_{\vartheta({\tt q})} <   \alpha   < x_{\vartheta({\tt q}+1)} \leq    \ldots
               \leq x_{\vartheta({\tt r})} \leq  \frac{1}{n+1} <  x_{\vartheta({ \tt r}+1)} \leq   \ldots
               \leq x_{\vartheta(n)} \, .                                                   $$ 
        Let    $ \Phi (\vec{x}) =: \vec{y} =: ( y_0, y_1,   \ldots   , y_n )$.  Since \
        $ \Phi |_{ {\clubsuit}_{n,\alpha} }   :
	          { \clubsuit}_{n,\alpha}    \stackrel{\cong}{\longrightarrow}   { \clubsuit}_{n,\beta} $ \
        we get that \ $\vec{y} \in \Delta_n\backslash {\clubsuit}_{n,\beta}$. \
        Then we claim that either \\
        \centerline { $ \beta < \min \{ y_0, y_1, \ldots , y_n \}$ \ or }
        $$ y_{\vartheta(0)} \leq    \ldots    \leq  y_{\vartheta({\tt q}-1)} \leq
        y_{\vartheta({\tt q})} <   \beta   < y_{\vartheta({\tt q}+1)} \leq    \ldots
        \leq   y_{\vartheta({\tt r})} <  y_{\vartheta({ \tt r}+1)} \leq   \ldots
        \leq y_{\vartheta(n)} \, . $$
   \end{lemma}
   \begin{proof}  \quad The proof follows the line of the previous Lemma  \ref{propositionlemma eins}.
   We have either \ $ \alpha = \beta = 0$ \ or \  $ 0 < \alpha , \beta < \frac{1}{n+1} $.
   The case $ \alpha = \beta = 0$  is trivial since \ ${ \clubsuit}_{n,0} = {\bf BOU}_n $, \ and we get \\
         \centerline{ $ \alpha = 0 < \min \{ x_0, x_1, \ldots , x_n \}$. See Lemma \ref{Lemma zwei}.
         It follows  $ \beta = 0 < \min \{ y_0, y_1, \ldots , y_n \}$. }  \\

   Let \ $ 0 < \alpha , \beta < \frac{1}{n+1} $. \
     Further we assume the second alternative, i.e.  we have a suitable number \
     $ { \tt q} $ \ with \ $ 0 \leq { \tt q} \leq { \tt r} < n $ \ such that
   $$ x_{\vartheta(0)} \leq x_{\vartheta(1)} \leq  \ldots  \leq   x_{\vartheta({\tt q})} <
                  \alpha   < x_{\vartheta({\tt q}+1)} \leq    \ldots   \leq x_{\vartheta({\tt r})} \leq
                  \frac{1}{n+1} <  x_{\vartheta({ \tt r}+1)} \leq   \ldots  \leq x_{\vartheta(n)}  \ . $$
       We have to show that (A): $ y_{\vartheta({\tt q})} < \beta$, \ and (B): $\beta < y_{\vartheta({\tt q}+1)}$.

   (A): \ If \  $x_{\vartheta({\tt q})} = 0$ \  it follows  from Lemma \ref {Lemma drei}  that
       $y_{\vartheta({\tt q})} = 0$, hence   $y_{\vartheta({\tt q})} < \beta$. \\
       Assume that  \ $x_{\vartheta({\tt q})} > 0$.  \  Similar as in the proof of 
       Lemma \ref{propositionlemma eins} we define a subset of $\Delta_n$,  \\
    \centerline{  $ Subset[\Delta_n] \ :=  \left\{\vec{a}_{\varepsilon} \, | \, \varepsilon \in [0,1] \right\}$.  } \\
       For all $ \varepsilon \in [0,1]$ \ let \
       $\vec{a}_{\varepsilon} := (a_0, a_1, \ldots, a_n) \in  Subset[\Delta_n]  $ \ by setting
      \[    a_j   :=
               \begin{cases}   x_j \cdot (1-\varepsilon)   &   \text{ for } \ \
                    j \in  \{ \vartheta(0) , \vartheta(1) , \ldots , \vartheta({\tt q}-1), \vartheta({\tt q})\}  \\
                    x_j    &   \text{ for } \ \  j \in
                    \{\vartheta({\tt q}+1), \vartheta({\tt q}+2), \ldots , \vartheta(n-1) \}  \\
                    x_{\vartheta(n)} +   \varepsilon \cdot \sum_{i=0}^{{\tt q}}   x_{\vartheta(i)}
                    &   \text{ for }  \ \   j =  \vartheta(n) \ .
                \end{cases}
       \]
      We have $ a_{\vartheta({\tt q})} < \alpha < a_{\vartheta({\tt q}+1)}$, hence
      $ \vec{a}_{\varepsilon} \in  \Delta_n \backslash { \clubsuit}_{n,\alpha}$
      for all $ \varepsilon \in [0,1]$.   We get $ \vec{a}_{0} = \vec{x}$,  and
       $ \vec{a}_{1} =:  ( \bar{a}_0, \bar{a}_1,   \ldots   , \bar{a}_n ) \in { \bf BOU}_n $   with
    $$          0 = \bar{a}_{\vartheta(0)} = \bar{a}_{\vartheta(1)} =   \ldots   =
               \bar{a}_{\vartheta({\tt q})} = 0  <  \alpha < \bar{a}_{\vartheta({\tt q}+1)} \leq
                \ldots   \leq \bar{a}_{\vartheta({\tt r})}  <  \bar{a}_{\vartheta({ \tt r}+1)} \leq  \ \ldots
                  \ \leq \bar{a}_{\vartheta(n)}   \ . $$
        Let   $ \Phi(\vec{a}_{\varepsilon} ) =: \vec{d}_{\varepsilon} =: ( d_o, d_1 , \ldots , d_n )$. 
        It is \ $ \vec{d}_{\varepsilon} \notin  {\clubsuit}_{n,\beta}\ , \text{for all} \ \varepsilon \in [0,1]$.
        Because   $ \Phi  \in   { \cal COMFORT}(\Delta_n)   , $    with \
         $ \Phi(\vec{a}_{1} ) = \vec{d}_{1} =: ( \bar{d}_o, \bar{d}_1 , \ldots   \ldots ,   \bar{d}_n )
         \notin { \clubsuit}_{n,\beta}\  $  it  follows that there is an index  $ {\tt j } \geq {\tt q} $
         such that   $  \bar{d}_{\vartheta( { \tt j})}  < \beta <  \bar{d}_{\vartheta( { \tt j}+1)}$,
         see Lemma \ref{Lemma drei},   and
     $$ 0 = \bar{d}_{\vartheta(0)} =   \ldots   = \bar{d}_{\vartheta({\tt q})} = 0  <
         \bar{d}_{\vartheta({\tt q}+1)} \leq    \ldots   < \bar{d}_{\vartheta( { \tt j})}  < \beta <
         \bar{d}_{\vartheta( { \tt j}+1)}<    \ldots    \bar{d}_{\vartheta( { \tt r})}\ \leq  \ldots    \leq
         \bar{d}_{\vartheta(n)}   \ .                                                         $$
         Note that \  $ Subset[\Delta_n] \cap {\clubsuit}_{n,\alpha} = \emptyset$, \ and since \
         $  \Phi |_{ {\clubsuit}_{n,\alpha}} :
	       { \clubsuit}_{n,\alpha}    \stackrel{\cong}{\longrightarrow} {\clubsuit}_{n,\beta} $ \ it also holds that
	   $$ \Phi \left(Subset[\Delta_n]\right) \cap {\clubsuit}_{n,\beta} \ = \ \emptyset  \  . $$
	       Further note that both sets
         $ Subset[\Delta_n]  \text{ and its continous image } \Phi\left(Subset[\Delta_n]\right)$  are connected
         subsets of $ \Delta_n$.  It follows that the projection  onto  the $\vartheta({\tt q})^{th}$  component
         is a connected set, i.e.  \ \
         ${ \sf PROJ_{   \it  \vartheta({\tt q}) }} (\Phi(Subset[\Delta_n])) $   \ is a connected subset of
         $ [0,1]$, i.e. an interval.

      Now we use  $ \Phi (\vec{x}) = \vec{y} \text{ and } \Phi (\vec{a}_{1}) = \vec{d}_{1}$, elements of
         $\Phi\left(Subset[\Delta_n]\right)$. \ They  have  the $ \vartheta({\tt q})^{th}$  component 
         $y_{\vartheta({\tt q})}  \text{ and }    \bar{d}_{\vartheta({\tt q})} $,   respectively.
         We had shown  that  $ \bar{d}_{\vartheta({\tt q})} = 0 $,  hence
     $$      \{ 0 , y_{\vartheta({\tt q})} \} \subset
                     { \sf PROJ_{   \it  \vartheta({\tt q}) }} (\Phi(Subset[\Delta_n])) , $$
         hence the closed interval \  $ [ 0 , y_{\vartheta({\tt q})} ]$ \ is a subset \ of
         $ { \sf PROJ_{   \it  \vartheta({\tt q}) }} (\Phi(Subset[\Delta_n]))$. \  With the empty set \
         $ \Phi(Subset[\Delta_n])  \cap {\clubsuit}_{n,\beta} = \emptyset $ \ it follows that \
         $ y_{\vartheta({\tt q})} <  \beta  $, and (A) is shown.

    (B):  Since  $ \Phi^{-1}  \in  { \cal COMFORT}(\Delta_n)$ we use the same argument as in
         Lemma \ref{propositionlemma eins}  to show \ $ \beta  < y_{\vartheta({\tt q}+1)}$.
         We can exchange the parts of \
         $ \alpha$  and $ \beta$ since \ $ \Phi^{-1} |_{ {\clubsuit}_{n,\beta}}:
	          { \clubsuit}_{n,\beta}  \stackrel{\cong}{\longrightarrow}  { \clubsuit}_{n,\alpha}$. \
	       With  \ $ \Phi^{-1}( y_0, y_1, \ldots , y_n ) = ( x_0, x_1, \ldots , x_n )$ \ we can
	       deduce that \  $ y_{\vartheta({\tt q}+1)} < \beta $ \  means \ $ x_{\vartheta({\tt q}+1)} < \alpha $,
	       and we get a contradiction.   Hence the only possibility is   \ $ y_{\vartheta({\tt q}+1)} > \beta$.
	
	   The first alternative of the lemma is  $ \alpha < \min \{ x_0, x_1, \ldots , x_n \} $. It is treated
	       correspondingly. If we assume that  $ \beta > \min \{ y_0, y_1, \ldots , y_n \} $ we can exchange the
	       parts of  $ \alpha  \text{ and }  \beta $  and consider  $ \Phi^{-1}$ and we would find a contradiction.
	       This ends the proof of Lemma \ref{lemma sieben}.
   \end{proof}

   With the previous rather technical lemmas we are able to discuss some possibilities to extend a map
         $ \varphi $   which is defined for a subset of   $ \Delta_n$   to the entire simplex.
	\begin{proposition}  \label{proposition zwei}
	\quad Let   $n \in  \mathbbm{N}$, and let us assume  $ 0 < \alpha , \beta \leq \frac{1}{n+1} $,  or
	  $ 0 = \alpha = \beta $.  Further we assume the existence of a homeomorphism
	  \ $ \varphi:  { \bf Layer}_{n,\alpha} \stackrel{\cong}{\longrightarrow} { \bf Layer}_{n,\beta}$.
	    Then \ $ \varphi$ can be extended to a homeomorphism on $ \Delta_n$,  i.e. there is a
	  homeomorphism  \  $\Phi:  \Delta_n \stackrel{\cong}{\longrightarrow} \Delta_n$ \  such that \
    $\Phi_{| { \bf Layer}_{n,\alpha} } = \varphi $.  \ The constructed
    $\Phi$   has the property that for any  $ 0 \leq \gamma \leq \frac{1}{n+1} $ there is a \
    $0 \leq \delta \leq \frac{1}{n+1} $  such that we get a homeomorphism  \
    $\Phi_{|{ \bf Layer}_{n,\gamma}}:{ \bf Layer}_{n,\gamma} 
    \stackrel{\cong}{\longrightarrow} { \bf Layer}_{n,\delta}$. \\
    If $ \varphi $  has the properties   ${\bf \widehat{[2]}}$ (respecting permutations) and
    ${ \bf  \widehat{[3]}}$  (keeping the order), then the homeomorphism  $\Phi$ is  an element of
    ${ \cal COMFORT}(\Delta_n)$.
	\end{proposition}	
	\begin{proof} \quad If   $\alpha = \frac{1}{n+1} $,    then is    $\beta = \frac{1}{n+1} $,   and we have
	        $ { \bf Layer}_{n,\alpha} = \{ { \bf Center}_n \} $.  Take \  $\Phi := id(\Delta_n)$.
	        For the other cases we need an auxiliary function $ \sigma $,  it must be any increasing homeomorphism
	        on the interval  $ \left[ 0, \frac{1}{n+1} \right] $   with
        	$ \sigma(\alpha) = \beta $. \	For instance  $ \sigma $  can be the polygon  through three  points
        	$ \left\{ (0,0), ( \alpha,\beta) , \left( \frac{1}{n+1}, \frac{1}{n+1} \right) \right\}$. \
        	Define   $ \Phi({ \bf Center}_n ) := {\bf Center}_n $.  \\
        	Now we assume   \ $ 0 < \alpha , \beta < \frac{1}{n+1} $ . \
          As we noted already on the second page of this section, for any
          $ \vec{x} = ( x_0, x_1 ,  \ldots   x_n ) \in \Delta_n\backslash\{{ \bf Center}_n \}$   we have an unique
	        $\pi(\vec{x}) \in { \bf BOU}_n$ \ and unique number ${\cal A}(\vec{x}) = \min \{ x_0, x_1 , \ldots x_n \}$,
	        $  0 \leq {\cal A}(\vec{x})  < \frac{1}{n+1}$, \ such that \
	        $\vec{x} \in {\bf Layer}_{n,{\cal A}(\vec{x}) }$, and $\vec{x} \in [{\bf Center}_n, \pi(\vec{x})] $, and
	  $$   \vec{x}  =  \pi(\vec{x}) +  {\cal A}(\vec{x}) \cdot (n+1) \cdot ( { \bf Center}_n - \pi(\vec{x})) \ . $$
	        We define  (we use the brackets  `$ \left\langle \, ... \,  \right\rangle$' for a better display)
        	$$ \Phi(\vec{x})   :=   \pi(\varphi(\pi_{\alpha}(\vec{x}))) +  \sigma({\cal A}(\vec{x})) \cdot (n+1) \cdot
	        \left\langle { \bf Center}_n - \pi(\varphi(\pi_{\alpha}(\vec{x}))) \right\rangle , $$
	        and  $\Phi$ has all the demanded properties. (Please note the remarks after 
	        Definition \ref{irgendeine Definition von Abbildungen}).
	        In case of  \ $ \alpha = \beta = 0 $ \ we can
	        take   $ \sigma := id\left( \left[ 0, \frac{1}{n+1} \right]\right)$, and we use the corresponding formula
        	$$ \Phi(\vec{x})   :=   \varphi(\pi(\vec{x})) +  {\cal A}(\vec{x}) \cdot (n+1) \cdot
	        \left\langle { \bf Center}_n - \varphi(\pi(\vec{x})) \right\rangle \ . $$
	        The proof is done.
	\end{proof}	
	Now we consider another possibility to extend a map on  ${ \bf BOU}_n$   to a map on $ \Delta_n$  with
	some properties preserved. This proposition will be very important in the next section.
	\begin{proposition}  \label{proposition drei}
           	\quad Let $n \in  \mathbbm{N}$. Let  $ 0 \leq \alpha , \beta < \frac{1}{n+1} $,  and let \
            $\varphi :  { \bf BOU}_n   \stackrel{\cong}{\rightarrow} { \bf BOU}_n $ be a homeomorphism such that
              $ \varphi $   yields a homeomorphism
        $$ \varphi_{| { \bf BOU}_n \, \cap \, { \clubsuit}_{n,\alpha} } :  \
            { \bf BOU}_n \cap   { \clubsuit}_{n,\alpha}    \stackrel{\cong}{\longrightarrow} \
            { \bf BOU}_n \cap   { \clubsuit}_{n,\beta}  \, .    $$
            Further, let  $ \varphi $  have the properties  ${\bf \widehat{[2]}}$  (respecting permutations) and
            ${\bf \widehat{[3]}}$ (keeping the order), i.e.  $ \varphi  \in { \cal COMFORT}({ \bf BOU}_n )$. \\
              Then    $ \varphi$   can be extended to a  homeomorphism on $ \Delta_n$, more precisely  there
            is a homeomorphism   \ $\Phi:  \Delta_n \stackrel{\cong}{\longrightarrow} \Delta_n$    such that \
            $\Phi_{|   { \bf BOU}_n } = \varphi   , $  the map    $\Phi$   is  an element  of
            ${ \cal COMFORT}(\Delta_n)   , $    and   $\Phi$   has the  property   \
            $ \Phi_{| { \clubsuit}_{n,\alpha} } :
              { \clubsuit}_{n,\alpha}  \stackrel{\cong}{\longrightarrow} { \clubsuit}_{n,\beta} \, .$    \
  \end{proposition}
  \begin{proof}  \quad
        In the case of   $ \alpha = \beta = 0$,  i.e. \
        $  { \bf BOU}_n = { \clubsuit}_{n,\alpha} = { \clubsuit}_{n,\beta}$,
        we have to look  at the previous Proposition \ref{proposition zwei}.   \
        The case    $ 0 < \alpha , \beta < \frac{1}{n+1} $   is treated here.   \\
        Recall  the  almost  disjoint  union   \
       $ \Delta_n  =  \bigcup \:  \left\{ \left[ { \bf Center}_n , \vec{b} \right] \, | \,
       \vec{b} \in {\bf BOU}_n \right\} $, \
	     all the lines intersect only in   ${ \bf Center}_n$. Since $ \varphi$ is a homeomorphism it follows \
	     $ \Delta_n =  \bigcup  \:  \left\{ \left[{\bf Center}_n, \varphi(\vec{b}) \right] \, | \,
	     \vec{b} \in  {\bf BOU}_n \right\} $. \
       We shall construct the homeomorphism  $\Phi$ by mapping every line $ \left[ {\bf Center}_n, \vec{b} \right]$
       homeomorphicly onto the line  $ \left[{ \bf Center}_n  , \varphi(\vec{b}) \right]$   such that if \
       $ \vec{x} = ( x_0, x_1,   \ldots   , x_n ) \in \left[{\bf Center}_n, \vec{b}\right]$  contains a component
       $ \alpha $, then   $\Phi(\vec{x}) \in \left[{\bf Center}_n, \varphi(\vec{b})\right]$ will  contain
       a  component $\beta$ at the same position.

    For  $ \vec{b} = ( b_0, b_1, \ldots , b_n ) \in { \bf BOU}_n $ \ there is a permutation $ \vartheta $  on
       $ \{ 0,1, \ldots , n \}$  and there are two natural  numbers \  ${ \tt q} , { \tt s}$,  with \
       $ 0 \leq { \tt q} \leq  { \tt s} < n $ \  such  that
     $$ 0 = b_{\vartheta(0)} \leq b_{\vartheta(1)} \leq     \ldots    \leq
         b_{\vartheta({\tt q})}  \leq \alpha < b_{\vartheta({\tt q}+1)} \leq    \ldots   \leq  b_{\vartheta({\tt s})}
         < \frac{1}{n+1} \leq  b_{\vartheta({ \tt s}+1)} \leq   \ldots   \leq b_{\vartheta(n)}   \ . $$
       (The case   $ { \tt q} = { \tt s} $ is possible).  \ Now we need two lemmas.  \\
        Let   \ $\varphi(\vec{b}) =: \vec{c} =: ( c_0, c_1,   \ldots   , c_n ) \in { \bf BOU}_n $.   Because 
        $ \varphi  \in { \cal COMFORT}({ \bf BOU}_n ) $  we have
      $$ 0 = c_{\vartheta(0)} \leq c_{\vartheta(1)} \leq   \ldots   \leq
        c_{\vartheta({\tt q})} <  c_{\vartheta({\tt q}+1)} \leq    \ldots   \leq c_{\vartheta({\tt s})}
        <  c_{\vartheta({ \tt s}+1)} \leq   \ \ldots   \ \leq c_{\vartheta(n)} \ . $$
   \begin{lemma} \label{lemma acht}   \  With the assumed properties of \  $ \varphi$ in  
        Proposition \ref{proposition drei}  it holds that \\
        \centerline{   $  c_{\vartheta({\tt q})}  \leq \beta < c_{\vartheta({\tt q}+1)} $,   \ and   \
               $  c_{\vartheta({\tt j})} = \beta $ \  if and only if \ $ b_{\vartheta({\tt j})} = \alpha $, \
               for $ 1 \leq {\tt j} \leq {\tt q} $.       }
          Further, if we have   \ \
          $ b_{\vartheta({\tt j})} < b_{\vartheta({\tt j}+1)} =  b_{\vartheta({\tt j}+2)}\ \ldots   \ldots
          = \, b_{\vartheta({\tt j}+  { \tt t} )} <  b_{\vartheta({\tt j} + { \tt t} +1  )}$   \ for suitable
          $ { \tt j}, { \tt t}$ with  $ 0 \leq { \tt j} \leq n-1 $ and  $ 1 \leq { \tt t} \leq n - { \tt j} $,
          we get \\
          \centerline{  $ c_{\vartheta({\tt j})} < c_{\vartheta({\tt j}+1)}
              =  c_{\vartheta({\tt j}+2)}\ \ldots   \ldots
              = \, c_{\vartheta({\tt j}+  { \tt t} )} <  c_{\vartheta({\tt j} + { \tt t} +1  )}$, \ and vice versa. }
    \end{lemma}
      \begin{proof}   \  We can argue as we did it in the  previous lemmas  \ref{Lemma fuenf} and
       \ref{lemma sieben}. See also Corollary  \ref{corolloary zwei}.
      \end{proof}

    \begin{lemma}  \quad For all \ $ \vec{b} = ( b_0, b_1, \ldots , b_n ) \in  { \bf BOU}_n $ \
        the number of intersections of \ $ \left[ {\bf Center}_n, \vec{b} \right] $ with $ {\clubsuit}_{n,\alpha}$
        is at most $ {\tt q} + 1$,    more precisely
        $$ { \rm cardinality} \left( \left[{ \bf Center}_n, \vec{b} \right] \cap {\clubsuit}_{n,\alpha}\ \right) = \
           {\rm cardinality} \left( \left\{ b_{\vartheta(0)}, b_{\vartheta(1)}, \ldots,
           b_{\vartheta({\tt q})} \right\} \right)  \geq 1 \, . $$
    \end{lemma}
    \begin{proof}   \quad      Recall that \
        for an \ $ \vec{x} = ( x_0, x_1, \ldots , x_n ) \in  \left[ {\bf Center}_n, \vec{b} \right] $ we have
        an unique  number 
        $ t \in [0,1]$ \ such that \ $ \vec{x} = t \cdot { \bf Center}_n + (1-t) \cdot \vec{b} $,  \ i.e.
        for a component  $ x_j$ we have
        $$ x_j = t \cdot \frac{1}{n+1} + (1-t) \cdot b_j = t \cdot \left(\frac{1}{n+1} - b_j \right) + b_j \ .$$
        If we use the  canonical projections   ${ \sf PROJ_{   \it i } } $, for $ i = 0,1, \ldots , n $,
        we have for a fixed index \ $j \in \{ 0, 1,   \ldots ,   { \tt s} \}$, \
        i.e.    $  b_{\vartheta(j)} < \frac{1}{n+1} $ ,   that the map  \\
        \centerline{ $ [0,1] \longrightarrow \left[ {\bf Center}_n, \vec{b} \right] \longrightarrow [0,1]  \ , \
                     t   \mapsto   \vec{x}   \mapsto   { \sf PROJ_{\vartheta(j)}}(\vec{x}) $, \quad   i.e. }  \\
        \centerline{ $ t   \longmapsto   t \cdot  { \bf Center}_n + (1-t) \cdot \vec{b}   \longmapsto
                     t \cdot \left( \frac{1}{n+1} - b_{\vartheta(j)} \right) + b_{\vartheta(j)}$, }    \\
        is strictly monotone increasing,  while for  \ $j \in \{ {\tt s}+1, {\tt s}+2, \ldots , n \} $
        (i.e. $ \frac{1}{n+1} \leq b_{\vartheta(j)}$)  the map \\
        $ [0,1] \longrightarrow [0,1] \ ,
        \ t  \mapsto t \cdot \left( \frac{1}{n+1} - b_{\vartheta(j)} \right) + b_{\vartheta(j)}$ \ is
        monotone decreasing. \ Hence all  components
        $ \left\{ b_{\vartheta(0)} , b_{\vartheta(1)} ,  \ldots   , b_{\vartheta({\tt q})} \right\}$
        meet $ \alpha$  any time while they are increasing to  $\frac{1}{n+1}\ $.
     \end{proof}
        Now we need  for all  $ \vec{b} = ( b_0, b_1, \ldots , b_n ) \in { \bf BOU}_n $  an increasing
        homeomorphism   \ $ \tau[\vec{b}]$   on   $ [0,1]$. \  Let  $ \tau[\vec{b}]$  be the polygon through the set
        of at most $ {\tt q} + 3$   points   \ (some may be identical)
        $$ \left\{ \ (0,0), \ \ldots \ , \ \left( \frac{ \alpha - b_{\vartheta({\tt q}-j)} }
        { \frac{1}{n+1} - b_{\vartheta({\tt q}-j)}}, \frac{ \beta - c_{\vartheta({\tt q}-j)}}
        { \frac{1}{n+1} - c_{\vartheta({\tt q}-j)} } \right), \ \ldots \ , \ (1,1) \right\},  \quad
                     \text{for} \ \ j = 0,1,  \ldots ,  {\tt q} \ .  $$
     \begin{lemma}  \label{Lemma zehn}
         For all   \ $ \vec{b} = ( b_0, b_1,   \ldots   , b_n ) \in { \bf BOU}_n $  the just constructed map
         $ \tau[\vec{b}]$   is a well defined increasing homeomorphism on   $ [0,1]$.
     \end{lemma}
     \begin{proof} \quad Easy.   Note   \
          $ b_{\vartheta(0)} \leq  \ldots \leq b_{\vartheta({\tt q})} \leq \alpha < b_{\vartheta({\tt q}+1)}$ \ and \
          $ c_{\vartheta(0)} \leq  \ldots \leq  c_{\vartheta({\tt q})} \leq \beta  < c_{\vartheta({\tt q}+1)}$.  \
          And for any index $ j \in  \{ 0,1, \ldots , n \} $  we have  $ b_j = \alpha $   if and only if
          $ c_j = \beta $, see the previous Lemma  \ref{lemma acht}.
      \end{proof}
     We have for the number \ 
         $ \widehat{t} := \frac{ \alpha - b_{\vartheta({\tt q}-j)}}{\frac{1}{n+1} - b_{\vartheta({\tt q}-j)}}$, \
         $j = 0,1,  \ldots ,  {\tt q}$, \  that the   $ \vartheta({\tt q}-j)^{th}$ component of an element \
         $ \vec{x} \in  \left[ {\bf Center}_n , \vec{b} \right]$
         is   $\alpha$, \ for \ $\vec{x} = \widehat{t} \cdot { \bf Center}_n + (1-\widehat{t}) \cdot \vec{b}$, \ i.e.
      $$ x_{\vartheta({\tt q}-j)}  =  \widehat{t} \cdot \left( \frac{1}{n+1} - b_{\vartheta({\tt q}-j)} \right) +
                b_{\vartheta({\tt q}-j)}  = \alpha , \ \  \text{ for } \ \ j = 0,1, \ldots, {\tt q}.      $$
      Now  we  define the map $\Phi$. \ For  all $\vec{b} \in { \bf BOU}_n$  we define for all points \
          $ \vec{x} \in \, \left[ {\bf Center}_n, \vec{b} \right]$, i.e.
        $ \vec{x} = t \cdot { \bf Center}_n + (1-t) \cdot \vec{b}$ \ for a suitable \ $ t \in [0,1]$, \,
         with   $ \varphi(\vec{b}) =: \vec{c}$ \ the image \ $\Phi(\vec{x})$ \ by
     $$   \Phi(\vec{x}) :=   \tau[\vec{b}](t) \cdot { \bf Center}_n + \left(1- \tau[\vec{b}](t)\right) \cdot
        \vec{c}  \ \in  \left[ {\bf Center}_n , \vec{c} \right] \, , $$
     and the map $\Phi$ fulfils all the properties which are demanded in  Proposition \ref{proposition drei}.
     The details are left to the reader.   Hence the proof of Proposition \ref{proposition drei}
     is finished.
  \end{proof}	
  We continue our investigations with a further  interesting proposition.
	\begin{proposition}    \label{proposition vier} \quad
	    Let  $n \in  \mathbbm{N}$. \ 	There is an injective group morphism \\
	   $ { \bf \Lambda_n}: \	\left( \left\{ f | f \text{ is an increasing homeomorphism on} \
	   \left[ 0,\frac{1}{n+1} \right] \right\}, \circ \right) \longrightarrow ({ \cal COMFORT}(\Delta_n), \circ ). $
	\end{proposition}
	\begin{remark}   \rm  \label{remark vier}  \quad
	      Of course, we can replace   $  \left[ 0,\frac{1}{n+1} \right]$   by any closed interval.
	      Note that the following construction can be used for   $n=0$, but    ${ \bf \Lambda_0}$   is not injective.
	\end{remark}
	\begin{proof}  \quad
            Let   $f$   be a  map with the above  conditions, i.e.
	          $f$  is an an increasing homeomorphism on   $ \left[0,\frac{1}{n+1}\right]$.
		        It follows that   $f(0) = 0$    and   $f \left(\frac{1}{n+1}\right) = \frac{1}{n+1}$ .
            Let \ $\vec{x} := (x_0,x_1,  \ldots , x_n) \in \Delta_{n}$. \
	          Trivially, there is a permutation  $\vartheta$  of \ $ \{ 0,1,2,   \  \ldots   \ , n \}$  
	          and there is an index   ${ \tt r} \in \{ 0,1,2, \ \ldots , n \}$    such that  	
	   $$ 0 \leq  x_{\vartheta(0)} \leq  x_{\vartheta(1)} \leq      \  \ldots   \ \leq
	       x_{\vartheta({ \tt r})} \leq  \frac{1}{n+1} < x_{\vartheta({ \tt r}+1)} \leq  x_{\vartheta({ \tt r}+2)} \
	       \leq \ \ldots \	\leq  x_{\vartheta(n-1)}  \leq  x_{\vartheta(n)} \leq  1   \ . $$
       	Then we define   \ $ {\bf \Lambda_n}(f) =: F:  \Delta_{n} \longrightarrow \Delta_{n}$. \
      	Let \ $F ( x_0,x_1, \ \ldots , x_n ) =: (y_0,y_1, \ \ldots , y_n)$.  \
      	For all   $ i \in \{0,1,   \ldots   , { \tt r}\}$ \ we set  
      $$         y_{\vartheta(i)}:=f \left(x_{\vartheta(i)}\right)   \,  .  $$
      	Note that in the case  \ ${ \tt r} = n$, \ we  have \
	      $ F(\vec{x} ) = \vec{x}  = \left( \frac{1}{n+1}, \ \ldots \ , \frac{1}{n+1} \right) = { \bf Center}_n $.  \\
	      For  ${ \tt r} < n$ we define two real numbers $ D, \delta$ by
	 $$ \quad D :=  \sum_{i=0}^{\tt r } \left(x_{\vartheta(i)} - y_{\vartheta(i)}\right) \quad  {\rm and } \quad
	    \delta  :=  \frac{D }{ \sum_{i={ \tt r}+1}^{n} \ \left( x_{\vartheta(i)} - \frac{1}{n+1} \right)} \ \ . $$
	    Finally let for all \ $ i \in \{ { \tt r}+1, { \tt r}+2,   \ldots ,   n \}$ \quad
	 $$y_{\vartheta(i)} :=  x_{\vartheta(i)} + \delta \cdot \left( x_{\vartheta(i)} - \frac{1}{n+1} \right)
	     =  x_{\vartheta(i)} \cdot ( 1 + \delta ) - \frac{\delta}{n+1} 	 \, , $$
	    and the definition of  $F$ is complete.
	    To prove Proposition \ref{proposition vier} we have  still  to verify that
	\begin{itemize}
        \item $F$   is  a map to   $\Delta_{n}$
        \item $F$   is injective
        \item $F$   is surjective
        \item $F$   is continuous
        \item $F$ fulfils the  conditions \  $ \widehat{[2]} \text{ and }  \widehat{[3]}
  	              \text{ of  Definition }  \ref{irgendeine Definition von Abbildungen} $.
  \end{itemize}
 \underline{F is  a map to  $\Delta_{n}$ :}  \quad We have for \
                     $ { \bf Center}_n  \neq \vec{x} = (x_0,x_1,  \ldots , x_n),  \text{ i.e. }  { \tt r} < n$:
 \begin{eqnarray*}
      1 - \sum_{i=0}^{n}   y_{\vartheta(i)} &  =  &
      \sum_{i=0}^{n}   (x_{\vartheta(i)} - y_{\vartheta(i)})   \ =   \ \sum_{i=0}^{{ \tt r}}
      (x_{\vartheta(i)} - y_{\vartheta(i)})  + \sum_{i={ \tt r}+1}^{n}   (x_{\vartheta(i)} - y_{\vartheta(i)})  \\
      & =  & \quad D    \ +   \   \sum_{i={ \tt r}+1}^{n}   (x_{\vartheta(i)} - y_{\vartheta(i)})   \\
      & = &  \quad  D   \ -   \    \sum_{i={ \tt r}+1}^{n}   \delta \cdot
                            \left( x_{\vartheta(i)} - \frac{1}{n+1} \right)
      \qquad  (\text{ see the definition of the } \ y_{\vartheta(i)} \text{'s})      \\
      & = & \quad  D   \  -   \  \frac{D }{ \sum_{i={ \tt r}+1}^{n}   \
                                            \left( x_{\vartheta(i)}    -    \frac{1}{n+1} \right) }
         \cdot \sum_{i={ \tt r}+1}^{n}   \left( x_{\vartheta(i)} - \frac{1}{n+1} \right)   \ =    \ 0  \ .
  \end{eqnarray*}
      Hence  \ $\sum_{i=0}^{n}   y_{\vartheta(i)} = 1 $. 
                                                            
       Now we show that \ $y_{\vartheta(j)}  \geq 0$ \ for all  $ j \in \{ 0, \ldots, n \}$.
      This is trivial for $ j \in \{ 0, \ldots , { \tt r} \}$ \ or for \ $D \geq 0 $. \
      Thus let  $j \in \{{ \tt r}+1, { \tt r}+2, \ldots, n \}$ \ and  \ $D <  0$, hence \ $ \delta < 0$.  \
      We show  the even stronger inequality  \ $y_{\vartheta(j)}   >  \frac{1}{n+1}$ .
  \begin{proof}  \quad  We have \ $\frac{1}{n+1} > y_{\vartheta(0)}$ \ and \ $ \frac{1}{n+1} \geq y_{\vartheta(i)} $ \
             for \ $ i = 0, 1, 2, \ldots { \tt r}$.  We get equivalences
  \begin{eqnarray*}
       \frac{{ \tt r}+1}{n+1}   \ >   \  \sum_{i=0}^{{ \tt r}}   y_{\vartheta(i)}
       & \Longleftrightarrow   & 1    -    \frac{n-{ \tt r}}{n+1}   \ >   \ \sum_{i=0}^{{ \tt r}}   y_{\vartheta(i)} \\
       & \Longleftrightarrow &  1    -    \ \frac{n-{ \tt r}}{n+1}   \ -   \sum_{i=0}^{{ \tt r}}   x_{\vartheta(i)}
         \  >   \   \sum_{i=0}^{{ \tt r}}   y_{\vartheta(i)}    \ -    \ \sum_{i=0}^{{ \tt r}}   x_{\vartheta(i)}
         \ =    \  - D  \\
       &  \Longleftrightarrow &    \left( \sum_{i={ \tt r}+1}^{n}   x_{\vartheta(i)} \right)   \   -
                                                       \ \frac{n-{ \tt r}}{n+1}   \ >   \  - D    \\
       & \Longleftrightarrow & 1   \ >   \ \frac{- D}{ \left( \sum_{i={ \tt r}+1}^{n} x_{\vartheta(i)} \right)
                                \ -   \ \frac{n-{ \tt r}}{n+1} }   \ =   \ - \delta  \\
       &  \Longleftrightarrow  &    1   \ >   \ - \delta
       \quad \Longleftrightarrow \quad    x_{\vartheta(j)}    -    \frac{1}{n+1}   \ >   \ (-\delta)
       \cdot   \left( x_{\vartheta(j)}    -    \frac{1}{n+1} \right)  \\
       & \Longleftrightarrow &   y_{\vartheta(j)}    =     x_{\vartheta(j)}   +   \delta   \cdot
       \left( x_{\vartheta(j)}   -   \frac{1}{n+1} \right)   \ >   \ +   \frac{1}{n+1} \ .
   \end{eqnarray*}
   \end{proof}
       Hence we have proved that the image of $F$ is $\Delta_n$, i.e. \ $F : \Delta_{n} \longrightarrow \Delta_{n}$ .
    \begin{lemma}   \label{wohl das zwoelfte lemma}
                 (1)  We can see from the above equivalences that \ $1 > - \delta $,  hence $ 1 + \delta > 0$,
                 and this fact easily leads to the conclusion that $F$ keeps the   order, i.e. the condition
                 ${\bf \widehat{[3]}}$  in Definition \ref{irgendeine Definition von Abbildungen} is fulfiled.  \\
                 (2) $F$ does not  change components $\{ x_0, \ldots , x_n \}$ of the set
    	           $ \left\{ 0,\frac{1}{n+1}, 1 \right\} \cup$  \{the fixed points  of $f$\}.	 \\
    	           (3) $F$ respects permutations on  $\Delta_{n}$,  that means if  \ $ \vartheta $ is any permutation
    	           of  \  $ \{ 0,1,2, \ldots , n \} $,  and if \
    	           $F \left(x_{0}, x_{1}, \ldots , x_{n}\right) = \left(y_{0}, y_{1}, \ldots ,y_{n}\right)$,
	               then we get that \ $F \left(x_{\vartheta(0)} , x_{\vartheta(1)}, \ldots , x_{\vartheta(n)}\right)$
	               $ = \left(y_{\vartheta(0)}, y_{\vartheta(1)}, \ldots , y_{\vartheta(n)}\right)$.  That means  that
	               ${\bf \widehat{[2]}}$  in Definition \ref{irgendeine Definition von Abbildungen} is fulfiled.
	   \end{lemma}
     \begin{proof}  \quad
         (1):   Let     $ 0 \leq  x_{\vartheta(0)} \leq  x_{\vartheta(1)} \leq      \  \ldots   \
	          \leq  x_{\vartheta({ \tt r})} \leq  \frac{1}{n+1} <  x_{\vartheta({ \tt r}+1)} \leq \ \ldots \
           	\leq  x_{\vartheta(n)} \leq  1 $.
            For   $ x_{i} < x_{j} \leq \frac{1}{n+1}$ we have  $ y_{i} < y_{j} \leq \frac{1}{n+1}$,
            because   $f$   is an homeomorphism.  If    $\frac{1}{n+1} < x_{i} <  x_{j}$ \ we defined \
            $y_{i}   =     x_{i} + \delta \cdot \left( x_{i} - \frac{1}{n+1} \right)$  \  and \
            $y_{j}   =    x_{j} + \delta \cdot \left( x_{j} - \frac{1}{n+1} \right) $. \ We can write \
            $ y_{i} =  x_{i} \cdot ( 1 + \delta ) - \frac{\delta}{n+1} <
              y_{j} = x_{j} \cdot ( 1 + \delta ) - \frac{\delta}{n+1}$ .   \quad  It follows  \\
            \centerline{    $ 0 \leq  y_{\vartheta(0)} \leq  y_{\vartheta(1)} \leq      \  \ldots   \
          	                \leq  y_{\vartheta({ \tt r})} \leq  \frac{1}{n+1} <
          	                 y_{\vartheta({ \tt r}+1)} \leq \ \ldots \ \leq y_{\vartheta(n)} \leq 1 $. }  \\
            (2) and (3):    Both facts follow easily from the construction of   $F$.
	     \end{proof}          $ $   \\
	   \underline{$F$ is injective :}   \ This is rather trivial.   \\
     \underline{$F$ is surjective :}   \ This is easy, too. Trivially, for an element \
               $ \vec{y} = (y_0,y_1,  \ldots , y_n) \in \Delta_{n}$ there are a permutation  \
               $\vartheta$ and an index   ${ \tt r} \in \{ 0,1,2, \ \ldots , n \}$  such that  	\\
	  \centerline{   $ 0 \leq  y_{\vartheta(0)} \leq  y_{\vartheta(1)} \leq  \ \ldots  \ \leq
	          y_{\vartheta({ \tt r})} \leq  \frac{1}{n+1} < y_{\vartheta({ \tt r}+1)} \leq
	          y_{\vartheta({ \tt r}+2)} \ \leq \ \ldots \ \leq  y_{\vartheta(n-1)}  \leq  y_{\vartheta(n)} \leq  1 $.} \\
	        Then we define the inverse map $F^{-1}(y_0,y_1,  \ldots , y_n) =: (x_0,x_1,  \ldots , x_n) $ \
	        by \  $ x_{\vartheta(j)} := f^{-1} \left(y_{\vartheta(j)}\right) $ \ for \
	        $ j \in \{ 0,1, \ldots ,  { \tt r} \}$, \ and (in the case of $ { \tt r} < n) $ \ we define \
	        $ D :=  \sum_{i=0}^{\tt r } \left(x_{\vartheta(i)} - y_{\vartheta(i)}\right)$ \ and
    $$ \delta := \frac{D} { \sum_{i={ \tt r}+1}^{n} \ \left(  y_{\vartheta(i)} - \frac{1}{n+1} \right) - D } \ . $$
	       We can define  for the other indices \ $ j \in \{ {\tt r}+1, {\tt r}+2, \ldots , n \}$ \ the components
            $$ x_{\vartheta(j)} := \left( y_{\vartheta(j)} + \frac{\delta}{n+1} \right) \cdot
                                 \left( { 1+ \delta } \right )^{-1}  \, . $$
         Then we get \ $ F(\vec{x}) = \vec{y} $, \ as well as \ $ F \circ F^{-1} = F^{-1} \circ F = id(\Delta_n)$. \\           \underline{$F$ is continuous :}   \ Let   \
        ${ \sf PROJ_{ \it j} } $ \ be the  canonical projection for all  \ $j \in \{ 0, \ldots ,n \}$,  \\
            $  { \sf PROJ_{ \it j} }: \Delta_{n} \longrightarrow [0,1],
              \ ( x_{0} , x_{1} ,   \ldots ,   x_{j},   \ldots ,    x_{n} )   \mapsto  x_{j}$  ,
              \	and let    \  $F_j   :=    { \sf PROJ_{ \it j} }   \circ   F $  .   \\
            Consider  the following commutative  diagram (Figure \ref{picture6}): \\  \\
   \begin{figure}[ht]
  \centering
  \setlength{\unitlength}{1cm}
  \begin{picture}(2.5,1.25)
  \put(-1,1.5){$\Delta_{n}$ }      \put(0,1.6){\vector(1,0){1}}    \put(0.5,1.8){$F$ }  
  \put(1.8,1.5){$\Delta_{n}$ }   \put(3,1.5){$\hookrightarrow$}  \put(4,1.5){$\mathbbm{R}^{n+1}$ }   
   \put(5,1.5){$ \ \ = \ \ \prod \mathbbm{R}$ }  \put(5.85,1.15){${\scriptstyle i=1}$}                                                                                   \put(5.8,1.875){${\scriptstyle n+1}$}
   \put(-0.8,1.3){\vector(1,-1){1}}  \put(2,1.3){\vector(-1,-1){1}}  
   \put(-0.7,0.5){$F_j$ }      \put(1.6,0.5){$ { \sf PROJ_{ \it j} }$ }     \put(0.3,-0.1){[0,1]} 
                                 \put(1.4,-0.1){$\hookrightarrow$}  \put(2.2,-0.1){$\mathbbm{R}$ }    
  \end{picture}    \\   \caption{}\label{picture6}    
  \end{figure}    
     \\  \\
            We see that   \  $F$ is continuous  if and only if  \ $F_j$  is continuous  for all
            $j  \in \{ 0, \ldots  , n \}$.  \\
            The continuity  of the $F_j$'s  are rather trivial; we have to put some attention to the case
            that  components of   $ ( x_{0} , x_{1} , \ldots  ,  x_{n} ) $   have the value   $\frac{1}{n+1}$ .

     This finishes the construction of \  $ { \bf \Lambda_n}(f) := F$, and  $F$ is an  element of
            ${ \cal COMFORT}(\Delta_n)$, since it was already shown in Lemma  \ref{wohl das zwoelfte lemma} 
            that  ${\bf \Lambda_n}(f)$  fulfils \ ${ \bf \widehat{[2]}} \text{ and } 
            { \bf \widehat{[3]}}  \text{ of  Definition } \ref{irgendeine Definition von Abbildungen} $.
     Furthermore, for two increasing homeomorphisms  \ $ f, g \text{ on } \left[0,\frac{1}{n+1}\right]$
            we can confirm with the aid of several calculations that we have the identities    \\
       \centerline  {  ${ \bf \Lambda_n} ( g \circ f ) =  { \bf \Lambda_n} ( g ) \circ { \bf \Lambda_n} ( f )$,
                     and  ${ \bf \Lambda_n} ( f^{-1} ) = ({ \bf \Lambda_n} ( f ))^{-1} , \ \text{and} \
             { \bf \Lambda_n} \left( id\left( [0,\frac{1}{n+1}] \right) \right) =  id(\Delta_n) $.   }
             This proves that  $ {\bf \Lambda_n}$   is a group morphism.  
             The injectivity of   $ { \bf \Lambda_n}$   is trivial. \ 
             Now the proof of Proposition \ref{proposition vier} is complete.
   	\end{proof}
    \begin{remark}	 \rm     \label{remark fuenf}	
          Note that for all  $ 0 \leq \alpha \leq \frac{1}{n+1} $ the homeomorphism $ {\bf \Lambda_n}(f)$ from
          Proposition \ref{proposition vier} maps the $ \alpha$-cross homeomorphicly onto the $ f(\alpha)$-cross,
       $$ {\bf \Lambda_n}(f) _{|{ \clubsuit}_{n,\alpha}}:   \
          { \clubsuit}_{{n,\alpha}} \stackrel{\cong}{\longrightarrow} { \clubsuit}_{{n, f(\alpha)}}    \ . $$
    \end{remark}
    \begin{lemma}
          The map  ${ \bf \Lambda_n}$ from { \rm Proposition \ref{proposition vier}} is not surjective.
    \end{lemma}
    \begin{proof}
         We need an element  which is not in the image of  $ { \bf \Lambda_n}$, for one $n \in  \mathbbm{N}$. \
          We define  $ F: \Delta_2 \rightarrow \Delta_2 $,
        \rm  \[  F_{   | { \bf BOU}_2 \cap { \bf Section}_{2,0} }   (0,x,1-x)  :=
             \begin{cases}
                    \quad \left(0, \frac{1}{2} \cdot x, 1 - \frac{1}{2} \cdot x \right)  &
                    \quad \mbox{for} \quad  x \in \left[0,\frac{1}{4} \right]     \\
                    \quad \left(0, \frac{5}{2} \cdot x - \frac{1}{2}, \frac{3}{2} - \frac{5}{2} \cdot x \right)  &
                    \quad \mbox{for} \quad  x \in \left[\frac{1}{4},\frac{1}{3}\right]     \\
                    \quad \left(0, x, 1 - x \right)  &
                    \quad \mbox{for} \quad  x \in \left[\frac{1}{3},\frac{1}{2}\right]
              \end{cases}           \]
        and   $ F_{   | { \bf BOU}_2 \cap { \bf Section}_{2,0} } $ can uniquely be extended on ${ \bf BOU}_2$
        by the properties ${\bf \widehat{[2]}}$  and  ${\bf \widehat{[3]}}$.
        With Proposition  \ref{proposition zwei} (we have the case \ $ \alpha = \beta = 0 $), \
        $F_{| { \bf BOU}_2}$ can be extended to  a homeomorphism  $ F \in { \cal COMFORT}(\Delta_2)$.
        This map \ $F$ \ has the property that for $0 \leq \mu \leq \frac{1}{3}$  \ we have a homeomorphism  \
        $ F_{|{ \bf Layer}_{2,\mu}}: { \bf Layer}_{2,\mu} \stackrel{\cong}{\rightarrow} { \bf Layer}_{2,\mu}$. \
        That means for an element   $(\mu ,v,w)  \in { \bf Layer}_{2,\mu}$ \ (i.e.  $\mu$  is the smallest 
        element of   $ \{ \mu, v, w \} $)  that  $ F( \mu,v,w ) = (\mu,\tilde{v},\tilde{w})$ for suitable numbers 
        $ \tilde{v},\tilde{w}$. \   Hence if there would exist a homeomorphism
        $ f: \left[0,\frac{1}{3}\right] \stackrel{\cong}{\rightarrow} \left[0,\frac{1}{3}\right]$   with \
        $  { \bf \Lambda_2}(f) = F $,   then   $ f(\mu) = \mu$   for all
        $ \mu \in  \left[0,\frac{1}{3}\right] $, i.e. $f$ would be the identity.
        This contradicts the above definition of  $F$ on ${ \bf BOU}_2 \cap { \bf Section}_{2,0}$ .
   \end{proof}
	 \newpage            
	 \section{Constructions in the Case of   $ L = 1$ }
	 \begin{remark}    \rm     \label{remark sechs}	
	         The reason for introducing  the map $ { \bf \Lambda_n}$ of Proposition  \ref{proposition vier}  
	         is that we shall take   the	searched homeomorphism   $\Theta_{1,n,0}$  from the image of
	         $ { \bf \Lambda_n}$, for each $n$. \  After that the maps  $ \Theta_{1,n,1}$  will  be defined  by 
	         induction on \ $n$ \ to  make the equations  { \tt   EQUATION$_{n,j \leq p,i,k}$ }  true.   
	 \end{remark}
	         We still are considering the case   $L = 1$. In the following we  omit this constant $L=1$.  \\
	           Remember that we have defined (in  Definition  \ref{zweite Definition}) the maps for $n=1$:    \\  \\
             $  \Theta_{1,0}(x,1-x) = (\eta(x),\eta(1-x)) \quad \text{and} \quad
	           \Theta_{1,1}(x,1-x) = (\kappa(x),\kappa(1-x))$ , for   $  (x,1-x) \in \Delta_1$.  \\
	           Note     \qquad  \qquad \qquad  \qquad  \quad
	           $\Theta_{1,0}\left(\frac{1}{4},\frac{3}{4}\right) = \left(\frac{1}{6},\frac{5}{6}\right)$, \ and \
	           $\Theta_{1,1}\left(\frac{1}{4},\frac{3}{4}\right) = \left(\frac{1}{5},\frac{4}{5}\right)$. \\  \\
	      Now we can continue the constructions of the homeomorphisms	\ $ \Theta_{n,0}$   and   $\Theta_{n,1}$ .                  They all will be elements of ${\cal COMFORT}(\Delta_n)$.   Four of them, \
	      $ \Theta_{n-1,0} ,   \Theta_{n-1,1} ,   \Theta_{n,0} ,   \Theta_{n,1}$,
	      will be needed to fulfil the equations     {\tt   EQUATION$_{n,j \leq p,i,k}$} ,
	      for all $i, k \in \{ 0, 1 \}$  and for $j , p \in \{ 0, 1, \ldots , n \} $ \
	      with  $j \leq p $,  for all fixed $n \in  \mathbbm{N}$.
   \begin{definition}   \  \rm  \label{Definition zehn}
         Let  $ n \in  \mathbbm{N}_0$ . First we define homeomorphisms   $ \phi_{n,0}$   on the interval \
         $ \left[0, \frac{1}{n+1} \right]$.
           Let    $\phi_{n,0} $   be the polygon through three points
         $$ \left\{ (0,0), \left(\frac{1}{2 \cdot (n+1)},\frac{1}{2 \cdot (n+2)}\right)   , \
                          \left(\frac{1}{n+1},\frac{1}{n+1}\right)   \right\}   \   . $$
        Now let  $\Theta_{n,0} := {\bf \Lambda_n}(\phi_{n,0})$, and by  Proposition \ref{proposition vier},
        $\Theta_{n,0}$ is an element of ${\cal COMFORT}(\Delta_n)$ for all numbers $n \in \mathbbm{N}$. \hfill $\Box$
   \end{definition}
        Note that  $\phi_{1,0} = \eta_{| \left[0, \frac{1}{2}\right]}$, hence the definition of \
        $\Theta_{1,0} = { \bf \Lambda_1}(\phi_{1,0}) $ \
        corresponds with those we  already   have given in Definition \ref{zweite Definition}, i.e.  \
        $\Theta_{1,0}(x,1-x) = (\eta(x),\eta(1-x))$. \ And note that \\
       \centerline{  $\phi_{n,0} \left(\frac{1}{2 \cdot (n+1)} \right) = \frac{1}{2 \cdot (n+2)}$,   \ \
                                      for all $ n \in \mathbbm{N}$.   }
   \begin{lemma} \label{lemma Dreieinhalb}
        \quad For  all \ $ n \in  \mathbbm{N}$, the  homeomorphisms    $ \Theta_{n,0}$   yield
        homeomorphisms   $$ \Theta_{n,0   | {{ \clubsuit}_{n,\frac{1}{2 \cdot (n+1)}}}}: \ \
          {{ \clubsuit}_{n,\frac{1}{2 \cdot (n+1)}}} \stackrel{\cong}{\longrightarrow}
                                                          {{ \clubsuit}_{n,\frac{1}{2 \cdot (n+2)}}}   \ . $$
   \end{lemma}
   \begin{proof} \quad See  Definition  \ref{vierte definition}   of the $\alpha$-cross   
           $ { \clubsuit}_{n,\alpha}$, \ and see  Proposition   \ref{proposition vier}.
   \end{proof}
        Now we construct the maps      $ \Theta_{n,1}$   for $ n \in  \mathbbm{N}$   by induction to
        fulfil { \tt   EQUATION$_{n,j \leq p,i,k}$}.  Note that  $\Theta_{0,1}$ and $ \Theta_{1,1} $ already exist. \
        Because of $i,k \in \{0,1\}$, we have to look at four possibilities, for \
        $(i,k) \in \{ (0,0), (0,1), (1,0), (1,1) \} $.   The reader should notice that in the case  \
        $i= k= 0$   the  { \tt   EQUATION$_{n,j \leq p,0,0}$} is   trivial. (The maps \
        $ \left\langle  id  \right\rangle_{n, 0, j }$ only add a component $0$ to the components of an element
        $ (x_0, \ldots, x_{n-1})$ of $\Delta_{n-1}$.) \
        We fix  $j := p := i := 0$ and $k := 1$ to construct the maps  $\Theta_{n,1}$  by induction  on $n$.   
                                                        
    At first we define   $\Theta_{n,1}$   on the topological boundary   ${\bf BOU_n}$ to fulfil
        { \tt   EQUATION$_{n,0 \leq 0,i=0,k=1}$}. After that we shall extend the map \
        $\Theta_{n,1}$   on   $\Delta_n$   to fulfil   {\tt EQUATION$_{n,0 \leq 0,i=1,k=1}$}, too.  \\  
                                                         
    To verify  the equation  \ { \tt   EQUATION$_{n,j \leq p,i=0,k=1}$}  we have to show the commutativity
        of Figure \ref{picture7}: 
  \begin{figure}[ht]
    \centering
    \setlength{\unitlength}{1cm}
    \begin{picture}(9,4)
     \put(-3.5,2.0){$\Delta_{n-1}$ }   \put(-2.5,2.9){$\vector(2,1){2}$}   \put(-2.0,2.7){$\Theta_{n-1,0}$}  
     \put(-0.0,3.9){$\Delta_{n-1}$ }   \put(1.0,4.0){$\vector(1,0){2}$} 
                                     \put(1.5,3.5){$\left\langle id \right\rangle_{\: n,\;  0,\; j \; }$}  
       \put(3.8,3.9){$\Delta_{n}$ }  
      \put(5.0,4.0){$\vector(1,0){2}$}   \put(5.4,3.55){$\Theta_{n,1}$}  \put(7.8,3.9){$\Delta_n$ }
      \put(10.0,3.4){$\left\langle id \right\rangle_{\: n+1,\;  1,\; p+1 \; }$}    \put(9.0,3.8){$\vector(2,-1){2}$}    
      \put(-2.5,1.2){$\vector(2,-1){2}$}    \put(-2.0,1.0){$\Theta_{n-1,1}$} 
      \put(-0.0,0.0){$\Delta_{n-1}$ }       \put(1.0,0.0){$\vector(1,0){2}$} 
                                              \put(1.5,0.3){$\left\langle id \right\rangle_{\: n,\;  1,\; p \; }$}  
        \put(3.8,-0.1){$\Delta_n$ }  
            \put(5.0,0.0){$\vector(1,0){2}$}      \put(5.7,0.2){$\Theta_{n,0}$}   \put(7.8,-0.1){$\Delta_n$ } 
           \put(9.0,0.2){$\vector(2,1){2}$}   \put(9.9,0.3){$\left\langle id \right\rangle_{\: n+1,\;  0,\; j \; }$}   
                    \put(11.0,2.0){$\Delta_{n+1} $ } 
    \end{picture}  \\       \caption{}\label{picture7}    
  \end{figure}        
   		\\    \\   \\    \\   \\   \\      
   	We describe  the way of the induction by writing down the case  $n= 2$   explicitly. \
   	We consider     {\tt EQUATION$_{n=2, j=0 \leq p=0,i=0,k=1}$} .   We have to show that
  	  \begin{equation*}
	        \left\langle id\right\rangle_{ 3, 0, 0 } \ \circ \ \Theta_{2,0} \ \circ \
  	              \left\langle  id  \right\rangle_{2, 1, 0 } \ \circ \ \Theta_{1,1} \  = \
  	       \left\langle id\right\rangle_{ 3, 1, 1 } \ \circ \ \Theta_{2,1} \ \circ \
  	              \left\langle  id  \right\rangle_{2, 0, 0 } \ \circ \ \Theta_{1,0} \, ,
  	    \end{equation*}         \\  \\
  	that means that we want the commutativity of  Figure \ref{picture8}:
   	  	  {      }  	\\   \\  
  \begin{figure}[ht]
  \centering
  \setlength{\unitlength}{1cm}
  \begin{picture}(8,4.0)
     \put(-3.25,2.0){$\Delta_{1}$ }   \put(-2.5,2.9){$\vector(2,1){2}$}   \put(-2.0,2.7){$\Theta_{1,0}$}  
     \put(0.0,3.9){$\Delta_{1}$ }   \put(1.0,4.0){$\vector(1,0){2}$} 
                                     \put(1.5,3.5){$\left\langle id \right\rangle_{\: 2,\;  0,\; 0 \; }$}  
       \put(3.8,3.9){$\Delta_{2}$ }  
      \put(5.0,4.0){$\vector(1,0){2}$}   \put(5.7,3.55){$\Theta_{2,1}$}  \put(7.8,3.9){$\Delta_2$ }
      \put(10.0,3.4){$\left\langle id \right\rangle_{\: 3,\;  1,\; 1 \; }$}    \put(9.0,3.8){$\vector(2,-1){2}$}    
            \put(-2.5,1.2){$\vector(2,-1){2}$}    \put(-2.0,1.0){$\Theta_{1,1}$} 
      \put(0.0,0.0){$\Delta_{1}$ }       \put(1.0,0.0){$\vector(1,0){2}$} 
                                              \put(1.5,0.3){$\left\langle id \right\rangle_{\: 2,\;  1,\; 0 \; }$}  
        \put(3.8,-0.1){$\Delta_2$ }  
            \put(5.0,0.0){$\vector(1,0){2}$}      \put(5.7,0.2){$\Theta_{2,0}$}   \put(7.8,-0.1){$\Delta_2$ } 
           \put(9.0,0.2){$\vector(2,1){2}$}   \put(9.9,0.3){$\left\langle id \right\rangle_{\: 3,\;  0,\; 0 \; }$}   
                    \put(11.0,2.0){$\Delta_{3} $ } 
  \end{picture}                  \\    \caption{}\label{picture8} 
  \end{figure} 
        { \ \   }  	\\   \\  
        Take an arbitrary element   \ $ ( x , 1-x ) \in \Delta_1 $,   $ x \in [0,1]$.
        Note that    $\phi_{2,0}\left(\frac{1}{6}\right) = \frac{1}{8}$,   and hence
        $ \Theta_{2,0}$   yields a  homeomorphism    $ \Theta_{2,0}|_{{{ \clubsuit}_{2,\frac{1}{6}}}}:
        {{ \clubsuit}_{2,\frac{1}{6}}} \stackrel{\cong}{\longrightarrow}  {{ \clubsuit}_{2,\frac{1}{8}}}$.
        We have \ $ \Theta_{1,1} \ \left( x, 1-x \right) =  \left( \kappa(x) , \kappa(1-x) \right) $,  
        and \ $\left\langle id \right\rangle_{\: 2,\;  1,\; 0 \; } \left( \kappa(x) , \kappa(1-x) \right) =
        \left(\frac{1}{6},\frac{5}{6}\cdot \kappa(x),\frac{5}{6}\cdot \kappa(1-x) \right) $. \
        This element will be mapped by \ $\Theta_{2,0}$,  \ 
        we write \ $\Theta_{2,0} \left(\frac{1}{6},\frac{5}{6}\cdot \kappa(x),\frac{5}{6}\cdot \kappa(1-x) \right)
        =: (\frac{1}{8},y,z)$ \ for  suitable real numbers \ $ y, z$. Shortly, an element  $( x , 1-x ) \in \Delta_1 $
        will  be mapped in the bottom line of the above Figure \ref{picture8} in the following way :
    $$   \left\langle id \right\rangle_{\: 3,\;  0,\; 0 \; } \circ  \Theta_{2,0}  \circ  \left\langle id
                  \right\rangle_{\: 2,\;  1,\; 0 \; } \circ  \Theta_{1,1} \ \left( x, 1-x \right) \
                  = \ \left( 0,\frac{1}{8},y,z \right)  \ . $$  
        Since  Figure  \ref{picture8} shall commute, the map \
        $\left\langle id \right\rangle_{\: 3,\;  1,\; 1 \; } $ \ from the upper line must take the element 
        $\left( 0, \frac{8}{7}\cdot y, \frac{8}{7}\cdot z \right)$ to map it to $\left( 0,\frac{1}{8},y,z \right)$. 
                                   
     We must define  $\Theta_{2,1}$ in a way that the following  diagram (Figure \ref{picture9}) commutes:
           \\   \\  \\ 
  \begin{figure}[ht]
  \centering
  \setlength{\unitlength}{1cm}
  \begin{picture}(9,4)
     \put(-3.5,2.0){$ \left(x , 1-x \right)$ }   \put(-2.5,2.9){$\vector(2,1){1.75}$}   \put(-2.0,2.7){$\Theta_{1,0}$}  
     \put(-1.5,4.0){$ \left( \eta(x),\eta(1-x) \right)  $ }   \put(1.4,4.1){$\vector(1,0){1.2}$} 
                                     \put(1.5,3.5){$\left\langle id \right\rangle_{\: 2,\;  0,\; 0 \; }$}  
       \put(2.8,4.0){$ \left( 0,\eta(x),\eta(1-x) \right)  $ }  
      \put(6.2,4.1){$\vector(1,0){1.3}$}   \put(6.6,3.55){$\Theta_{2,1}$}  
                             \put(7.8,4.0){$ \left(0,\frac{8}{7}\cdot y,\frac{8}{7}\cdot z \right) $ }
      \put(10.5,3.2){$\left\langle id \right\rangle_{\: 3,\;  1,\; 1 \; }$}    \put(9.5,3.6){$\vector(2,-1){2}$}    
            \put(-2.5,1.2){$\vector(2,-1){2}$}    \put(-2.0,1.0){$\Theta_{1,1}$} 
      \put(-0.9,-0.2){$ \left( \kappa(x),\kappa(1-x) \right)  $ }       \put(2.0,0.0){$\vector(1,0){1.5}$} 
                                              \put(2.0,0.3){$\left\langle id \right\rangle_{\: 2,\;  1,\; 0 \; }$}  
        \put(3.6,-0.1){$ \left( \frac{1}{6},\frac{5}{6}\cdot \kappa(x),\frac{5}{6}\cdot \kappa(1-x) \right)  $ }  
            \put(7.75,0.0){$\vector(1,0){0.7}$}      \put(7.8,0.3){$\Theta_{2,0}$}   
                             \put(8.6,-0.1){$(\frac{1}{8},y,z)$ } 
           \put(9.6,0.5){$\vector(2,1){2}$}   \put(10.5,0.6){$\left\langle id \right\rangle_{\: 3,\;  0,\; 0 \; }$}   
                    \put(11.0,1.9){$ \left(0,\frac{1}{8},y,z \right) $ } 
  \end{picture}            \\    \caption{}\label{picture9} 
  \end{figure}        
      $ { } $  	\\
    Hence we must define \ 
    $\Theta_{2,1} \left( 0, \eta(x), \eta(1-x) \right) := \left( 0, \frac{8}{7}\cdot y, \frac{8}{7}\cdot z \right)$.
    This  will be explained  now in details. 
                                                
  We are able  to  change  the  directions of  the maps  $\Theta_{1,0}, \
    \left\langle id \right\rangle_{\: 2,\;  0,\; 0 \; }$, and $\left\langle id \right\rangle_{\: 3,\;  1,\; 1 \; }$,
    respectively.  We  call these maps
    $\Theta_{1,0}^{-1}, \ \left\langle id \right\rangle_{\: 2,\; 0,\; 0 \; }^{-1}$, and
    $\left\langle id \right\rangle_{\: 3,\;  1,\; 1 \; }^{-1}$,  respectively. Since $\Theta_{1,0}$ is a
    homeomorphism, the meaning of   $\Theta_{1,0}^{-1}$ is clear. And
    $ \left\langle id \right\rangle_{\: 2,\;  0,\; 0 \; }^{-1}$ acts on the 3-tuple
    $ \left(0,\eta(x),\eta(1-x) \right) $ by deleting the $0$ at the first place,
    $$ \left\langle id \right\rangle_{\: 2,\;  0,\; 0 \; }^{-1} \left(0,\eta(x),\eta(1-x)\right)   :=
                                          \left(\eta(x),\eta(1-x)\right) \, , $$
      $\left\langle id \right\rangle_{\: 3,\;  1,\; 1 \; }^{-1}$   acts on \
       $\left(0,\frac{1}{8},y,z\right)$   by deleting  the fraction   $\frac{1}{8}$   at the second place and 
        multiplying the other components with    $\frac{8}{7}$ ,
      $$\left\langle id \right\rangle_{\: 3,\;  1,\; 1 \; }^{-1} \left(0,\frac{1}{8},y,z \right)   :=
                           \left(0,\frac{8}{7} \cdot y, \frac{8}{7} \cdot z \right)   .  $$
        Note that the maps    $ \left\langle id \right\rangle_{\: 2,\;  0,\; 0 \; }^{-1}$   and
        $\left\langle id \right\rangle_{\: 3,\;  1,\; 1 \; }^{-1}$   are left inverse maps, that means it holds
        $$ \left\langle id \right\rangle_{\: 2,\;  0,\; 0 \; }^{-1} \circ
         \left\langle id \right\rangle_{\: 2,\;  0,\; 0 \; }   =   id( \Delta_1) \qquad \text{as well as}  \qquad
        \left\langle id \right\rangle_{\: 3,\;  1,\; 1 \; }^{-1} \circ
        \left\langle id \right\rangle_{\: 3,\;  1,\; 1 \; }   =   id( \Delta_2)   \  .  $$
        And note that we also have the identity \
        $ \left\langle id \right\rangle_{\: 3,\;  1,\; 1 \; } \circ  \left\langle id
                  \right\rangle_{\: 3,\;  1,\; 1 \; }^{-1} \left(0,\frac{1}{8},y,z\right)
                  =  \left(0,\frac{1}{8},y,z\right)$. \

   If we turn around the directions of the maps  $\Theta_{1,0}$,
        $\left\langle id \right\rangle_{\: 2,\;  0,\; 0 \; }$, $\left\langle id \right\rangle_{\: 3,\;  1,\; 1 \; }$,
        respectively,  in the described way,  and if  we demand the  commutativity  of  the above  diagram,
        the map $\Theta_{2,1}$ is uniquely defined on  the face \ ${ \bf BOU}_2 \cap { \bf Section}_{2,0}$.
                                     
    That means that we define   $\Theta_{2,1}$  for every   3-tuple   
        $(0,\eta(x),\eta(1-x)) \in  { \bf BOU}_2 \cap { \bf Section}_{2,0}$, \ i.e. we define
        $$\Theta_{2,1}|_{ { \bf BOU}_2 \cap { \bf Section}_{2,0} }
              :=    \left\langle id \right\rangle_{\: 3,\;  1,\; 1 \; }^{-1}   \circ
            \left\langle id \right\rangle_{\: 3,\;  0,\; 0 \; }   \circ   \Theta_{2,0}   \circ
            \left\langle id \right\rangle_{\: 2,\;  1,\; 0 \; }   \circ
            \Theta_{1,1}   \circ   \Theta_{1,0}^{-1}   \circ    \left\langle id \right\rangle_{\: 2,\;  0,\; 0 \; }^{-1}
              \,  ,  $$
         $$  \Theta_{2,1}|_{ { \bf BOU}_2 \cap { \bf Section}_{2,0}  }    (0,\eta(x),\eta(1-x))    \ :=   \
             \left( 0,\frac{8}{7}\cdot y,\frac{8}{7}\cdot z \right)  \, ,   $$
        see  the following  commutative diagram  (Figure \ref{picture10}). \
        There we start at $ \left(0,\eta(x),\eta(1-x) \right)$:  \\  \\  \\  \\  \\  \\
  \begin{figure}[ht]
     \centering
     \setlength{\unitlength}{1cm}
     \begin{picture}(9,4)
          \put(-4.0,2.5){$ (x , 1-x)$ }  \put(-1.05,4.25){$\vector(-2,-1){1.75}$} \put(-2.0,3.2){$\Theta_{1,0}^{-1}$}
          \put(-1.5,4.5){$ \left(  \eta(x),\eta(1-x) \right)  $ }   \put(2.6,4.6){$\vector(-1,0){1.2}$}
                                     \put(1.5,4.0){$\left\langle id \right\rangle_{\: 2,\;  0,\; 0 \; }^{-1}$}
          \put(2.8,4.5){$ \left(  0,\eta(x),\eta(1-x) \right)  $ }
          \put(6.2,4.6){$\vector(1,0){1.3}$}   \put(6.6,4.05){$\Theta_{2,1}$}
                             \put(7.8,4.5){$\left(0,\frac{8}{7}\cdot y,\frac{8}{7}\cdot z\right)$ }
          \put(10.5,3.8){$\left\langle id \right\rangle_{\: 3,\;  1,\; 1 \; }^{-1}$}
          \put(11.5,3.1){$\vector(-2,1){2}$}
          \put(-3.0,1.7){$\vector(2,-1){2}$}    \put(-2.0,1.5){$\Theta_{1,1}$}
          \put(-0.9,0.3){$( \kappa(x),\kappa(1-x)) $ }       \put(2.0,0.5){$\vector(1,0){1.5}$}
          \put(2.0,0.8){$\left\langle id \right\rangle_{\: 2,\;  1,\; 0 \; }$}
          \put(3.6,0.4){$\left( \frac{1}{6},\frac{5}{6}\cdot \kappa(x),\frac{5}{6}\cdot \kappa(1-x)\right) $ }
          \put(7.75,0.5){$\vector(1,0){0.7}$}      \put(7.8,0.8){$\Theta_{2,0}$}
                                                   \put(8.6,0.4){$\left( \frac{1}{8},y,z \right) $ }
          \put(9.6,1.0){$\vector(2,1){2}$}   \put(10.5,1.1){$\left\langle id \right\rangle_{\: 3,\;  0,\; 0 \; }$}
                                            \put(11.0,2.5){$\left(0,\frac{1}{8},y,z\right)$ }
        \end{picture}    \\       \caption{}\label{picture10}
     \end{figure}
   {      }   	\\  
   Up to now the map  $\Theta_{2,1}$ is a homeomorphism on ${ \bf BOU}_2 \cap { \bf Section}_{2,0} $, i.e.
   it maps triples   $ (0,a,b) \in \Delta_2$.
   We had fixed the positions $j = p = 0$. If we vary $j,p \in \{ 0,1,2 \} $ with $ j \leq p $, \ we get
   $\Theta_{2,1}$ \ defined on  the other faces
   \ ${ \bf BOU}_2 \cap { \bf Section}_{2,1}$ and ${\bf BOU}_2 \cap { \bf Section}_{2,2}$,
   \ respectively. The seven maps \ $\Theta_{1,0}, \ \Theta_{1,1}, \ \left\langle id \right\rangle_{\: 2,\;  0,\; j}, \
   \left\langle id \right\rangle_{\: 2,\;  1,\; p }, \ \Theta_{2,0}, \
   \left\langle id \right\rangle_{\: 3,\;  0,\; j } ,  \left\langle id \right\rangle_{\: 3,\;  1,\; p+1 }$  \
   respect permutations, see Lemma   \ref{lemma eins}, hence there are no contradictions in the definition of
   $\Theta_{2,1}$   on   ${ \bf BOU}_2   $,   the boundary of $ \Delta_2$ .
  \begin{lemma}   \quad
      The just constructed map  $\Theta_{2,1}|_{ { \bf BOU}_2} $ is a homeomorphism on  $ { \bf BOU}_2 $.
      Further, it satisfies the  conditions \  $ {\bf \widehat{[2]}}$ (respecting permutations) and 
      $ {\bf \widehat{[3]}}$ (keeping the order). 
  \end{lemma}
  \begin{proof}   { $ $  }
     	\begin{itemize}
         \item Continuity:   $\Theta_{2,1}|_{\ { \bf BOU}_2 \cap { \bf Section}_{2,j} } $   is a product
               of seven continuous maps  \\
               $  \left\langle id \right\rangle_{\: 2,\;  0,\; j \; }^{-1}, \, \Theta_{1,0}^{-1}, \, \Theta_{1,1}, \,
               \left\langle id \right\rangle_{\: 2,\;  1,\; p \; }, \,  \Theta_{2,0}, \,
               \left\langle id \right\rangle_{\: 3,\;  0,\; j \; },  \,
               \left\langle id \right\rangle_{\: 3,\;  1,\; p+1 \; }^{-1}$,  \ for $ j,p \in \{ 0,1,2 \}$, $j\leq p$.
         \item    $\Theta_{2,1}|_{ { \bf BOU}_2} $ is injective, because all seven maps are injective.
         \item    $\Theta_{2,1}|_{ { \bf BOU}_2} $ is surjective on   $ { \bf BOU}_2 $.
         \item    The properties  $ {\bf \widehat{[2]}}$ and $  {\bf \widehat{[3]}}$ are easy to verify 
                  because of the construction of  $\Theta_{2,1}|_{ { \bf BOU}_2} $.   
                  See  Figure \ref{picture10}, and note that all maps \ $\Theta_{1,0}, \Theta_{1,1}, \Theta_{2,0}$ \
                  are from ${ \cal COMFORT}(\Delta_1)$ or ${ \cal COMFORT}(\Delta_2)$, respectively.
                  Note Lemma  \ref{lemma eins}.
      \end{itemize}
  \end{proof}  
    If the reader starts with the element \
      $ \left( \frac{1}{4}, \frac{3}{4} \right) \in  { \clubsuit}_{1,\frac{1}{4}} $ in the above diagram (Figure \ref{picture9}),
      she or he  can see  that
      $$ \Theta_{2,1}{   _{ |   { \bf BOU}_2   \cap   { \clubsuit}_{2,\frac{1}{6}} } } :      \
           { \bf BOU}_2   \cap   { \clubsuit}_{2, \frac{1}{6}}    \stackrel{\cong}{\longrightarrow} \
           { \bf BOU}_2   \cap   { \clubsuit}_{2, \frac{1}{7}} , \    \text{ e.g. }
           \left( 0,\frac{1}{6},\frac{5}{6} \right) \longmapsto \left( 0,\frac{1}{7},\frac{6}{7} \right) \, . $$
      Now the reader should remember Proposition \ref{proposition drei}.  By this proposition the map
      $\Theta_{2,1}|_{ { \bf BOU}_2} $  can be extended to a
            homeomorphism    $\Theta_{2,1}$   on
            $ \Delta_2 $ , even      $\Theta_{2,1} \in { \cal COMFORT}(\Delta_2)$, such that
      $$ \Theta_{2,1}{   _{ |   { \clubsuit}_{2,\frac{1}{6}} } } :      \
             { \clubsuit}_{2, \frac{1}{6}}    \stackrel{\cong}{\longrightarrow} \
             { \clubsuit}_{2, \frac{1}{7}}     \   .  $$

     By the construction of  $\Theta_{2,1}|_{ { \bf BOU}_2}$, the  {\tt EQUATION$_{2,j \leq p,i,k}$}
     is satisfied for the two pairs $(i,k) \in \{ (0,1),(1,0) \} $.  For  $i=k=0$,
     { \tt   EQUATION$_{2,j \leq p,0,0}$ }   is trivial. (The maps \ $\left\langle  id  \right\rangle_{2, 0, j }$
     only add a third component $0$ to the components of an  element  $ (y,1-y)$ of $\Delta_{1}$). 

     It remains to consider the case $i=k=1$, i.e. we must show  { \tt   EQUATION$_{2,j \leq p,1,1}$}.
     With the property   \  $ \Theta_{2,1}{   _{ |   { \clubsuit}_{2,\frac{1}{6}}}} :  \
                { \clubsuit}_{2, \frac{1}{6}} \stackrel{\cong}{\longrightarrow} \ {\clubsuit}_{2, \frac{1}{7}} $ \
                  it is easy to see that the following diagram (Figure \ref{picture11})  \\  \\  \\ 
  \begin{figure}[ht]
    \centering
    \setlength{\unitlength}{1cm}
    \begin{picture}(9.5,6)
     \put(-3.5,2.4){$ (x , 1-x)$ }   \put(-2.5,3.4){$\vector(2,1){1.75}$}   \put(-2.0,3.2){$\Theta_{1,1}$}  
     \put(-0.9,4.5){$ \left( \kappa(x),\kappa(1-x) \right)  $ }             \put(2.0,4.6){$\vector(1,0){1.2}$} 
                                     \put(2.0,4.0){$\left\langle id \right\rangle_{\: 2,\; 1,\; 0 \; }$}  
       \put(3.6,4.5){$ \left( \frac{1}{6},\frac{5}{6}\cdot \kappa(x),\frac{5}{6}\cdot \kappa(1-x)  \right)  $ }  
      \put(7.76,4.6){$\vector(1,0){0.7}$}                                   \put(7.8,4.05){$\Theta_{2,1}$}  
                                     \put(8.6,4.5){$ \left( \frac{1}{7},v,w \right) $ }
      \put(10.5,3.7){$\left\langle id \right\rangle_{\: 3,\;  1,\; 1 \; }$} \put(9.5,4.1){$\vector(2,-1){2}$}    
                                     \put(-2.5,1.7){$\vector(2,-1){2}$}     \put(-2.0,1.5){$\Theta_{1,1}$} 
      \put(-0.9,0.3){$ \left( \kappa(x),\kappa(1-x)  \right)  $ }           \put(2.0,0.5){$\vector(1,0){1.5}$} 
                                     \put(2.0,0.8){$\left\langle id \right\rangle_{\: 2,\;  1,\; 0 \; }$}  
      \put(3.6,0.4){$ \left( \frac{1}{6},\frac{5}{6}\cdot \kappa(x),\frac{5}{6}\cdot \kappa(1-x)  \right)  $ }  
                                     \put(7.75,0.5){$\vector(1,0){0.7}$}    \put(7.8,0.8){$\Theta_{2,1}$}   
                             \put(8.6,0.4){$ \left( \frac{1}{7},v,w \right) $ } 
      \put(9.6,1.0){$\vector(2,1){2}$}  
      \put(10.5,1.1){$\left\langle id \right\rangle_{\: 3,\; 1,\; 0 \; }$}   
      \put(10.5,2.4){$ \left( \frac{1}{8},\frac{1}{8},\frac{7}{8} \cdot v, \frac{7}{8} \cdot w  \right) $ } 
     \end{picture}          \\         \caption{}\label{picture11} 
  \end{figure}        
      $  \\   $  
   commutes, with suitable numbers $ v, w $, and  we also have confirmed that 
      {\tt EQUATION$_{2,0 \leq 0,i=1,k=1}$}   holds.  We get similar commutative diagrams if we take other pairs 
      $ (j,p) $, for  $ j,p \in \{ 0,1,2 \}$   with   $ j \leq p $ .   \                   
   Hence, for  $n=2$  all $ 24 $ cases of  { \tt   EQUATION$_{2,j \leq p,i,k}$}   are proved, for
       $ i,k \in \{ 0,1 \} , $   and   $ j,p \in \{ 0,1,2 \}$   with   $ j \leq p $ .  \\  \\

    At  this point we make a summary of the results that we have got so far.
      We have two trivial  homeomorphisms   $\Theta_{0,0},  \Theta_{0,1}$  on  $ \Delta_0 = \{1\} $,
      and we have defined  two    homeomorphisms   $\Theta_{1,0},  \Theta_{1,1}$ on  $\Delta_1$   and  two
      homeomorphisms  $\Theta_{2,0}, \Theta_{2,1}$  on  $\Delta_2$, respectively.   All four constructed
      homeomorphisms are even from  ${ \cal COMFORT}(\Delta_1)$ or  ${ \cal COMFORT}(\Delta_2)$,  respectively.  \
      Furthermore, the four homeomorphisms   $\Theta_{0,0}, \Theta_{0,1}, \Theta_{1,0}, \Theta_{1,1}$  \
      satisfy  the equations  \ {\tt EQUATION$_{n=1,j \leq p,i,k}$}, for $j,p \in \{0,1\}$  with  $j \leq p$,
      and  the four homeomorphisms  $\Theta_{1,0},  \Theta_{1,1}, \Theta_{2,0}$ and  $\Theta_{2,1}$  \
      satisfy  the equations   \ {\tt EQUATION$_{n=2,j \leq p,i,k}$}, \
      for  $j,p \in \{0,1,2\}$ with  $j \leq p$, and always for all \ $i,k \in \{ 0,1 \}$. \\
   \begin{theorem}  \label{zweites theorem}
          We are able to construct a homeomorphism   $\Theta_{n,1} \in { \cal COMFORT}(\Delta_n) $
          for all natural numbers $ n $,  such that, together with the already defined \   
          $\Theta_{n,0} = { \bf \Lambda_{n}}(\phi_{n,0}) \in { \cal COMFORT}(\Delta_{n}) $,  the four  maps \
          $\Theta_{n-1,0} , \Theta_{n-1,1}, \Theta_{n,0}, \Theta_{n,1}$   satisfy  the equation
          {\tt   EQUATION$_{n,j \leq p,i,k}$},  for integers \
          $  j,p \in \{0,1,2,   \ldots , n \} $ with $ j \leq p$,  and   $ i,k \in \{ 0,1 \}$.  \\
          Moreover, if we restrict   $\Theta_{n,1}$ to the    $ \frac{1}{2 \cdot (n+1)}$-cross of \
          $ \Delta_n$,  then  $\Theta_{n,1}$  maps $ { \clubsuit}_{{n,\frac{1}{2 \cdot (n+1)}}}$   homeomorphicly
          onto the   $ \frac{1}{2 \cdot (n+1)+1}$-cross,    \
      $$     \Theta_{n,1} | { \clubsuit}_{{n,\frac{1}{2 \cdot (n+1)}}}: \ \
             { \clubsuit}_{{n,\frac{1}{2 \cdot (n+1)}}} \stackrel{\cong}{\longrightarrow}
             { \clubsuit}_{{n,\frac{1}{2 \cdot (n+1) + 1 }}}   \ .  $$
   \end{theorem}   { $  $ }  \\  \\
   With  Theorem  \ref{zweites theorem} we would get, because of { \tt EQUATION$_{n,j \leq p,i,k}$} \
        for \ $  i := 0 $ and $ k := 1 $,  that the following diagram (Figure \ref{picture12})  commutes, for \
        $  j,p \in \{0,1,2, \ldots, n \} $ with $ j \leq p$ :  \\ \\ 
    \begin{figure}[ht]
  \centering
  \setlength{\unitlength}{1cm}
  \begin{picture}(9,4)
     \put(-3.5,2.0){$\Delta_{n-1}$ }   \put(-2.5,2.9){$\vector(2,1){2}$}   \put(-2.0,2.7){$\Theta_{n-1,0}$}  
     \put(-0.2,3.9){$\Delta_{n-1}$ }   \put(1.0,4.0){$\vector(1,0){2}$} 
                                     \put(1.5,3.5){$\left\langle id \right\rangle_{\: n,\;  0,\; j \; }$}  
       \put(3.8,3.9){$\Delta_{n}$ }  
      \put(5.0,4.0){$\vector(1,0){2}$}   \put(5.4,3.55){$\Theta_{n,1}$}  \put(7.8,3.9){$\Delta_n$ }
      \put(10.0,3.4){$\left\langle id \right\rangle_{\: n+1,\;  1,\; p+1 \; }$}    \put(9.0,3.8){$\vector(2,-1){2}$}    
      \put(-2.5,1.2){$\vector(2,-1){2}$}    \put(-2.0,1.0){$\Theta_{n-1,1}$} 
      \put(-0.2,-0.1){$\Delta_{n-1}$ }       \put(1.0,0.0){$\vector(1,0){2}$} 
                                              \put(1.5,0.3){$\left\langle id \right\rangle_{\: n,\;  1,\; p \; }$}  
        \put(3.8,-0.1){$\Delta_n$ }  
            \put(5.0,0.0){$\vector(1,0){2}$}      \put(5.7,0.2){$\Theta_{n,0}$}   \put(7.8,-0.1){$\Delta_n$ } 
           \put(9.0,0.2){$\vector(2,1){2}$}   \put(10.0,0.3){$\left\langle id \right\rangle_{\: n+1,\;  0,\; j \; }$}   
                    \put(11.0,2.0){$\Delta_{n+1} $ } 
  \end{picture}  \\       \caption{}\label{picture12}    
  \end{figure}        
      {      }   	\\    \\ \\   	\\    \\   
   \begin{proof} (of Theorem \ref{zweites theorem}).   \ \ 
   We repeat the beginning of the induction from the third section. \\
    \underline{Start of the induction:} We had defined the \ $ \frac{1}{4}$-cross of  $ \Delta_1$,
       $ {\clubsuit}_{{1,\frac{1}{4}}} =
       \left\{ \left( \frac{1}{4}, \frac{3}{4} \right), \left( \frac{3}{4}, \frac{1}{4} \right) \right\}$.
     The constructed two homeomorphisms  $\Theta_{1,0}, \Theta_{1,1}$    are elements  of  \
      ${ \cal COMFORT}(\Delta_1)$, and together with  $\Theta_{0,0},  \Theta_{0,1}$
    they fulfil the equations   \  { \tt   EQUATION$_{n=1,j \leq p,i,k}$}   for
      $ j,p \in \{ 0,1 \} $   with   $ j \leq p$,   and   $ i,k \in \{ 0,1 \}$.  \
       Moreover, if we restrict  $\Theta_{1,0}$  and  $\Theta_{1,1}$ to ${\clubsuit}_{{1,\frac{1}{4}}}$, then
       $\Theta_{1,0}$  maps \ $ { \clubsuit}_{{1,\frac{1}{4}}}$  onto the
       $ \frac{1}{6}$-cross $ { \clubsuit}_{{1,\frac{1}{6}}}$, and  $\Theta_{1,1}$  maps
       $ { \clubsuit}_{{1,\frac{1}{4}}}$   onto $ { \clubsuit}_{{1,\frac{1}{5}}}$. \  (E.g. we have
       $\Theta_{1,1} \left( \frac{1}{4}, \frac{3}{4} \right) = \left( \frac{1}{5}, \frac{4}{5} \right) $). \
       Shortly, we have maps
    $$ \Theta_{1,0}|  { \clubsuit}_{{1,\frac{1}{4}}}: {\clubsuit}_{{1,\frac{1}{4}}} \stackrel{\cong}{\longrightarrow}
        { \clubsuit}_{{1,\frac{1}{6}}} \ , \qquad \text{and} \qquad  \Theta_{1,1}| { \clubsuit}_{{1,\frac{1}{4}}}:
        \ {\clubsuit}_{{1,\frac{1}{4}}} \stackrel{\cong}{\longrightarrow}  {\clubsuit}_{{1,\frac{1}{5}}} \, , $$
        which are trivially homeomorphisms. \\
   \underline{The induction step from $n$ to $n+1$:} \
       Let for  an   $ n \in  \mathbbm{N}$   for all
       $ q \in \{ 0,1,2, \ldots , n \}$   the homeomorphisms   $\Theta_{q,1}$   on   $\Delta_q$
       be constructed,   $\Theta_{q,1}$   and   $\Theta_{q,0}$  be elements of  ${ \cal COMFORT}(\Delta_q)$,
       and four at a time be used to  satisfy the equations  {\tt   EQUATION$_{q,j \leq p,i,k}$},  for  all
       $ j,p \in \{ 0,1,   \ldots, q \}$, with $j \leq p$,  and  $ i,k \in \{ 0,1 \}$.   \
       Furthermore, if we restrict  $\Theta_{q,0}$  and  $\Theta_{q,1}$ to the  $\frac{1}{2 \cdot (q+1)}$-cross
       of   $ \Delta_q$, then we get a homeomorphism   \
           $ \Theta_{q,1}:    { \clubsuit}_{{q,\frac{1}{2 \cdot (q+1)}}} \stackrel{\cong}{\longrightarrow}
           { \clubsuit}_{{q,\frac{1}{2 \cdot (q+1) + 1 }}}, $ \ (by the assumption of the induction),   \
           and also a homeomorphism  \
           $\Theta_{q,0}: {\clubsuit}_{{q,\frac{1}{2 \cdot (q+1)}}} \stackrel{\cong}{\longrightarrow}
          { \clubsuit}_{{q,\frac{1}{2 \cdot (q+2)}}} $ \
           (by the construction of  $\Theta_{q,0}$, see Lemma  \ref {lemma Dreieinhalb}).

    That means for  $q := n$  and for $j := p := i := 0$  and  $ k := 1$
      we assume by  the induction hypothesis that the following equation {\tt   EQUATION$_{n,j=0 \leq p=0,i=0,k=1}$}
      holds, i.e. we assume that the following diagram (Figure \ref{picture13})  commutes,
        \begin{equation*}
	            \left\langle id \right\rangle_{\: n+1,\;  0,\; 0 \; } \ \circ \ \Theta_{n,0}  \ \circ \
  	            \left\langle id \right\rangle_{\: n,\;  1,\; 0 \; } \ \circ \ \Theta_{n-1,1}  \quad  =  \quad 	
  	         \left\langle id \right\rangle_{\: n+1,\;  1,\; 1 \; } \ \circ \  \Theta_{n,1} \ \circ \
  	          \left\langle  id  \right\rangle_{\: n,\;  0,\; 0 \; } \ \circ  \  \Theta_{n-1,0}   \   ,
       \end{equation*}  	
           \\  \\    \\   \\    \\   \\ 
    \begin{figure}[ht]
  \centering
  \setlength{\unitlength}{1cm}
  \begin{picture}(8,4)
     \put(-3.5,2.5){$\Delta_{n-1}$ }   \put(-2.5,3.4){$\vector(2,1){2}$}   \put(-2.0,3.2){$\Theta_{n-1,0}$}  
     \put(-0.0,4.4){$\Delta_{n-1}$ }   \put(1.2,4.5){$\vector(1,0){2}$} 
                                     \put(1.5,4.0){$\left\langle id \right\rangle_{\: n,\;  0,\; 0 \; }$}  
       \put(3.8,4.4){$\Delta_{n}$ }  
      \put(5.0,4.5){$\vector(1,0){2}$}   \put(5.4,4.05){$\Theta_{n,1}$}  \put(7.8,4.4){$\Delta_n$ }
      \put(10.0,3.9){$\left\langle id \right\rangle_{\: n+1,\;  1,\; 1 \; }$}    \put(9.0,4.3){$\vector(2,-1){2}$}    
            \put(-2.5,1.7){$\vector(2,-1){2}$}    \put(-2.0,1.5){$\Theta_{n-1,1}$} 
      \put(-0.0,0.5){$\Delta_{n-1}$ }       \put(1.2,0.5){$\vector(1,0){2}$} 
                                              \put(1.5,0.8){$\left\langle id \right\rangle_{\: n,\;  1,\; 0 \; }$}  
        \put(3.8,0.4){$\Delta_n$ }  
            \put(5.0,0.5){$\vector(1,0){2}$}      \put(5.7,0.7){$\Theta_{n,0}$}   \put(7.8,0.4){$\Delta_n$ } 
           \put(9.0,0.7){$\vector(2,1){2}$}   \put(10.0,0.8){$\left\langle id \right\rangle_{\: n+1,\;  0,\; 0 \; }$}   
                    \put(11.0,2.5){$\Delta_{n+1} $ } 
    \end{picture}          \caption{}\label{picture13}   
  \end{figure}          
      \\   \\  
    Now we want to make  the induction step \ $ n \mapsto n+1$. We fix $j := p := i := 0$ and  $ k := 1$, 
        i.e. we want to show { \tt EQUATION$_{n+1,j=0 \leq p=0,i=0,k=1}$}. That means that we shall show the
        commutativity of    Figure \ref{picture14}. \\
   \begin{figure}[ht]
  \centering
  \setlength{\unitlength}{1cm}
  \begin{picture}(8,5)
     \put(-3.5,2.5){$\Delta_{n}$ }   \put(-2.5,3.4){$\vector(2,1){2}$}   \put(-2.0,3.2){$\Theta_{n,0}$}  
     \put(-0.0,4.4){$\Delta_{n}$ }   \put(1.0,4.5){$\vector(1,0){2}$} 
                                     \put(1.5,4.0){$\left\langle id \right\rangle_{\: n+1,\;  0,\; 0 \; }$}  
       \put(3.8,4.4){$\Delta_{n+1}$ }  
      \put(5.0,4.5){$\vector(1,0){2}$}   \put(5.4,4.05){$\Theta_{n+1,1}$}  \put(7.8,4.4){$\Delta_{n+1}$ }
      \put(10.0,3.9){$\left\langle id \right\rangle_{\: n+2,\;  1,\; 1 \; }$}    \put(9.0,4.3){$\vector(2,-1){2}$}    
            \put(-2.5,1.7){$\vector(2,-1){2}$}    \put(-2.0,1.5){$\Theta_{n,1}$} 
      \put(-0.0,0.5){$\Delta_{n}$ }       \put(1.0,0.5){$\vector(1,0){2}$} 
                                              \put(1.5,0.8){$\left\langle id \right\rangle_{\: n+1,\;  1,\; 0 \; }$}  
        \put(3.8,0.4){$\Delta_{n+1}$ }  
            \put(5.0,0.5){$\vector(1,0){2}$}      \put(5.7,0.7){$\Theta_{n+1,0}$}   \put(7.8,0.4){$\Delta_{n+1}$ } 
           \put(9.0,0.7){$\vector(2,1){2}$}   \put(10.0,0.8){$\left\langle id \right\rangle_{\: n+2,\;  0,\; 0 \; }$}   
                    \put(11.0,2.5){$\Delta_{n+2} $ } 
  \end{picture}      \caption{}\label{picture14}      
  \end{figure}        
    \\  \\ \\   
   The homeomorphism  $ \Theta_{n,1} \in { \cal COMFORT}(\Delta_n)$ already exists by the
   assumption of the induction. And we already had defined  $\Theta_{n,0}$ and
   $\Theta_{n+1,0} = { \bf \Lambda_{n+1}}(\phi_{n+1,0}) \in { \cal COMFORT}(\Delta_{n+1})$.
   Hence it lacks   $\Theta_{n+1,1}$.  Note that  $\Theta_{n+1,0}$   yields a
    homeomorphism  (see Definition \ref{Definition zehn}),
      $$ \Theta_{n+1,0   | {{ \clubsuit}_{n+1,\frac{1}{2 \cdot (n+2)}}}}: \
                {{ \clubsuit}_{n+1,\frac{1}{2 \cdot (n+2)}}} \stackrel{\cong}{\longrightarrow}
                    {{ \clubsuit}_{n+1,\frac{1}{2 \cdot (n+3)}}} \ , $$
   as it was mentioned in  Lemma \ref{lemma Dreieinhalb}.  \\
   We show the the construction of   \ $\Theta_{n+1,1}$ in the way we have described it for   $n = 2 $.
   Let  $ \vec{x}$  be an arbitrary element of $\Delta_n$.  We name  the images of  $ \vec{x}$ by \
   $\Theta_{n,0}(\vec{x}) =: \vec{y}$, and  $\Theta_{n,1}(\vec{x}) =:  \vec{z}$.  Hence
     $ \left\langle id \right\rangle_{\: n+1,\;  0,\; 0 \; }(\vec{y}) = ( 0, \vec{y} )$
   and        $ \left\langle id \right\rangle_{\: n+1,\;  1,\; 0 \; }(\vec{z}) =
   \left( \frac{1}{ 2 \cdot (n+2) }, \left( 1 - \frac{1}{ 2 \cdot (n+2) } \right) \cdot \vec{z} \right)
   \in   \Delta_{n+1}$ .    We call    \ $\Theta_{n+1,1}  ( 0, \vec{y} ) =:  ( 0, \vec{v} )$   and
      $\Theta_{n+1,0} \left(  \frac{1}{ 2 \cdot (n+2) }, \left( 1 - \frac{1}{ 2 \cdot (n+2) } \right)
    \cdot \vec{z} \right) =:  \left( \frac{1}{ 2 \cdot (n+3)} , \vec{w} \right) $.  \
    Hence the bottom line of the above diagram (Figure \ref{picture14}) \ is
    $$ \left\langle id \right\rangle_{\: n+2,\;  0,\; 0 \; } \circ  \Theta_{n+1,0} \circ
    \left\langle id \right\rangle_{\: n+1,\;  1,\; 0 \; }  \circ \Theta_{n,1} \ (\vec{x})
       =   \left( 0, \frac{1}{ 2 \cdot (n+3)} , \vec{w} \right)   \  \in   \Delta_{n+2}  \ . $$
   The upper line of the above diagram is
   \begin{eqnarray*}
       &  \left\langle id \right\rangle_{ \: n+2,\; 1,\; 1 \; } \circ  \Theta_{n+1,1} \circ
          \left\langle id \right\rangle_{\: n+1,\;  0,\; 0 \; }  \circ \Theta_{n,0} \  (\vec{x})
             =      \left\langle id \right\rangle_{\: n+2,\; 1,\; 1 \; } ( 0, \vec{v} ) \\
       &  =      \left(  0, \frac{1}{ 2 \cdot (n+3)} ,
          \left( 1 - \frac{1}{ 2 \cdot (n+3)} \right) \cdot \vec{v} \right)     \in   \Delta_{n+2}   \  .
   \end{eqnarray*}
    Because we want the commutativity of the above diagram, we   must define    $\Theta_{n+1,1}$   such that
      $  \vec{w} =   \left( 1 - \frac{1}{ 2 \cdot (n+3)} \right) \cdot \vec{v}$. \
    We want the commutativity of  Figure  \ref{picture15}, \\  
    \begin{figure}[ht]
       \centering
       \setlength{\unitlength}{1cm}
       \begin{picture}(9.6,5)
          \put(-3.2,2.5){$\vec{x}$ }   \put(-2.5,3.4){$\vector(2,1){2}$}   \put(-2.0,3.2){$\Theta_{n,0}$}  
          \put(-0.0,4.4){$\vec{y}$ }   \put(1.0,4.5){$\vector(1,0){2}$} 
                                     \put(1.5,4.0){$\left\langle id \right\rangle_{\: n+1,\;  0,\; 0 \; }$}  
          \put(3.8,4.4){$( 0, \vec{y} )$ }  
          \put(5.0,4.5){$\vector(1,0){2}$}   \put(5.4,4.05){$\Theta_{n+1,1}$}  \put(7.8,4.4){$( 0, \vec{v} )$ }
          \put(9.4,4.3){$\left\langle id \right\rangle_{\: n+2,\;  1,\; 1 \; }$}    \put(9.0,4.3){$\vector(2,-1){1.5}$}    
           \put(-2.5,1.7){$\vector(2,-1){2}$}    \put(-2.0,1.5){$\Theta_{n,1}$} 
           \put(-0.2,0.4){$\vec{z}$ }     
                                \put(0.2,1.25){$\left\langle id \right\rangle_{\: n+1,\;  1,\; 0 \; }$}  
           \put(1.4,0.4){$ \left(  \frac{1}{ 2 \cdot (n+2) }, \left( 1 - \frac{1}{ 2 \cdot (n+2) } \right)
                                            \cdot \vec{z} \right)   $ } 
           \put(6.0,0.5){$\vector(1,0){0.75}$}    \put(0.25,0.5){$\vector(1,0){0.75}$}                              
           \put(6.0,1.25){$\Theta_{n+1,0}$}    \put(6.8,0.4){$\left( \frac{1}{ 2 \cdot (n+3)} , \vec{w} \right)$ } 
           \put(9.0,0.7){$\vector(2,1){1.5}$}   \put(9.4,0.5){$\left\langle id \right\rangle_{\: n+2,\;  0,\; 0 \; }$}                      \put(11.0,2.4){$ \|  $ } 
                    \put(10.0,1.8){$\left( 0, \frac{1}{ 2 \cdot (n+3)} , \vec{w} \right) $ } 
                    \put(8.25,3.0){$ \left(  0, \frac{1}{ 2 \cdot (n+3)} , 
                              \left( 1 - \frac{1}{ 2 \cdot (n+3)} \right) \cdot \vec{v} \right)   $ } 
         \end{picture}    \caption{}\label{picture15} 
    \end{figure}     
    \\  \\  \\                     
      To force the commutativity,  we   reverse the directions of
      $    \Theta_{n,0} , \left\langle id \right\rangle_{\: n+1,\;  0,\; 0 \; }$,   and
      $  \left\langle id \right\rangle_{\: n+2,\; 1,\; 1 \; } $, respectively, in a way we described it before in
      the case $ n = 2 $.    
                  
   Since  $ \Theta_{n,0}$ is a homeomorphism the meaning of $ \Theta_{n,0}^{-1}$ is clear. We have   
      \  $ \Theta_{n,0}^{-1}(\vec{y}) = (\vec{x})$.  
      We define the left inverse maps  $ \left\langle id \right\rangle_{\: n+1,\;  0,\; 0 }^{-1}$   and
      $ \left\langle id \right\rangle_{\: n+2,\; 1,\; 1 \; }^{-1}$.  \
      Let \  $ \left\langle id \right\rangle_{\: n+1,\;  0,\; 0 }^{-1} ( 0, \vec{y} ) := (\vec{y})$. \
      For any $(n+1)$-tuple   $ \vec{w} $   such that
        $ \left(  0, \frac{1}{ 2 \cdot (n+3)},\vec{w} \right)   \in  \Delta_{n+2} $ \ we must define
    $$ \left\langle id \right\rangle_{\: n+2,\; 1,\; 1 \; }^{-1}
            \left(  0, \frac{1}{ 2 \cdot (n+3)},\vec{w}  \right)
            :=   \left(  0, \left( 1 - \frac{1}{ 2 \cdot (n+3)} \right)^{-1} \cdot \vec{w} \right)
            =    \left(  0,  \frac{2 \cdot n + 6}{2 \cdot n + 5}  \cdot \vec{w} \right)     \in  \Delta_{n+1}   \ . $$
    Hence we get the identity
   $$  \left\langle id \right\rangle_{\: n+2,\;1,\;1\;} \circ \left\langle id \right\rangle_{\: n+2,\; 1,\; 1 \; }^{-1}
       \left( 0, \frac{1}{ 2 \cdot (n+3)},\vec{w}  \right) =  \left( 0, \frac{1}{ 2 \cdot (n+3)},\vec{w} \right) \, . $$
     At first we define  $\Theta_{n+1,1}$ on the face $ { \bf BOU}_{n+1} \cap { \bf Section}_{n+1,0} $ :
    $$\Theta_{n+1,1}  :=   \left\langle id \right\rangle_{\: n+2,\; 1,\; 1 \; }^{-1}  \circ
         \left\langle id \right\rangle_{\: n+2,\;  0,\; 0 \; }   \circ   \Theta_{n+1,0}   \circ
         \left\langle id \right\rangle_{\: n+1,\;  1,\; 0 \; }   \circ   \Theta_{n,1}     \circ  
         \Theta_{n,0}^{-1} \circ  \left\langle id \right\rangle_{\: n+1,\;  0,\; 0 \; }^{-1}   \  ,  $$
    $$  i.e. \qquad   \Theta_{n+1,1} ( 0, \vec{y} )  
         \ := \ \left(  0, \left( 1 - \frac{1}{ 2 \cdot (n+3)} \right)^{-1} \cdot \vec{w} \right) \ ,   $$ 
     and   $\Theta_{n+1,1}$   is uniquely defined on  the face $ { \bf BOU}_{n+1} \cap { \bf Section}_{n+1,0} $.
           We get the commutativity of Figure \ref{picture16}, see the following diagram. We start at  $(0,\vec{y})$:
    \begin{figure}[ht]
       \centering
       \setlength{\unitlength}{1cm}
       \begin{picture}(9,5)
         \put(-3.2,2.5){$\vec{x}$ }   \put(-0.5,4.4){$\vector(-2,-1){2}$}   \put(-2.0,3.2){$\Theta_{n,0}^{-1}$}  
         \put(-0.0,4.4){$\vec{y}$ }   \put(3.0,4.5){$\vector(-1,0){2}$} 
                                     \put(1.5,4.0){$\left\langle id \right\rangle_{\: n+1,\;  0,\; 0 \; }^{-1}$}  
         \put(3.8,4.4){$( 0, \vec{y} )$ }  
         \put(5.0,4.5){$\vector(1,0){2}$}   \put(5.4,4.05){$\Theta_{n+1,1}$} 
         \put(7.8,4.4){$ \left( 0, \left( 1 - \frac{1}{ 2 \cdot (n+3)} \right)^{-1} \cdot \vec{w} \right)$ }
         \put(10.1,3.5){$\left\langle id \right\rangle_{\: n+2,\;  1,\; 1 \; }^{-1}$}                                                                                                    \put(10.6,2.9){$\vector(-3,2){1.5}$}    
            \put(-2.5,1.7){$\vector(2,-1){2}$}    \put(-2.0,1.5){$\Theta_{n,1}$} 
         \put(-0.3,0.5){$\vec{z}$ }   
         \put(0.2,1.25){$\left\langle id \right\rangle_{\: n+1,\;  1,\; 0 \; }$}  
         \put(1.4,0.4){$ \left( \frac{1}{ 2 \cdot (n+2) }, \left( 1 - \frac{1}{ 2 \cdot (n+2) } \right)
                                            \cdot \vec{z} \right) $ }                                           
         \put(6.0,0.5){$\vector(1,0){0.75}$}    \put(0.25,0.5){$\vector(1,0){0.75}$} 
         \put(6.0,1.25){$\Theta_{n+1,0}$} 
         \put(6.8,0.4){$\left( \frac{1}{ 2 \cdot (n+3)} , \vec{w} \right)$ } 
         \put(9.0,0.7){$\vector(2,1){1.5}$}   \put(10.1,0.75){$\left\langle id \right\rangle_{\: n+2,\;  0,\; 0 \; }$}            \put(10.0,2.2){$\left( 0, \frac{1}{ 2 \cdot (n+3)} , \vec{w} \right) $ } 
       \end{picture}         \caption{}\label{picture16} 
     \end{figure} 
                  \\  \\  \\    \\  \\         
   \begin{lemma}    \label{lemma sechzehn}
     The just constructed map   $\Theta_{n+1,1}|_{ { \bf BOU}_{n+1} \cap { \bf Section}_{n+1,0} } $ \
     is a homeomorphism on     $ { \bf BOU}_{n+1} \cap { \bf Section}_{n+1,0} $.
     Further, it satisfies the  conditions \  $ {\bf \widehat{[2]}}$ (respecting permutations) and 
      $ {\bf \widehat{[3]}}$ (keeping the order)  of  Definition  \ref{irgendeine Definition von Abbildungen}.  
  \end{lemma}     
  \begin{proof}   { $ $  }
     	\begin{itemize}
           \item Continuity:  It is a product of seven continuous maps.
           \item  $\Theta_{n+1,1}|_{ { \bf BOU}_{n+1} \cap { \bf Section}_{n+1,0} }  $ is injective,
                   because all seven maps are injective.
           \item  $\Theta_{n+1,1}|_{ { \bf BOU}_{n+1} \cap { \bf Section}_{n+1,0} } $
                            is surjective on   $ { \bf BOU}_{n+1} \cap { \bf Section}_{n+1,0} $ .
           \item  $ {\bf \widehat{[2]}}$  and  $ {\bf \widehat{[3]}}$  are easy to verify.  
                  See  Figure \ref{picture16}, and note that all maps \ $\Theta_{n,0}, \Theta_{n,1}, \Theta_{n+1,0}$ \
                  are from \ ${ \cal COMFORT}(\Delta_n)$ \ or \ ${ \cal COMFORT}(\Delta_{n+1})$, respectively.                            Note Lemma  \ref{lemma eins}.
       \end{itemize}
   \end{proof}

   As yet the map    $\Theta_{n+1,1}$   is defined on   ${ \bf BOU}_{n+1} \cap { \bf Section}_{n+1,0} $ .
       We had fixed the positions $ j = p = 0 $. If we vary $ j,p$ and take them from the set \
       $ \{ 0,1, \ldots , n, n+1 \} $ \ with $ j \leq p $,  we get  $\Theta_{n+1,1}$ \ defined
       on  the other faces ${ \bf BOU}_{n+1} \cap { \bf Section}_{n+1,j}$,   for    
        {$ 1 \leq j \leq n+1$}.                                               
       All   maps  \ $\Theta_{n,0}, \, \Theta_{n,1}, \, \left\langle id \right\rangle_{\: n+1,\;  0,\; j \; }, \,
       \left\langle id \right\rangle_{\: n+1,\; 1,\; p }, \Theta_{n+1,0}, \  \left\langle id
       \right\rangle_{\: n+2,\;  0,\; j }, \, \left\langle id \right\rangle_{\: n+2,\;  1,\; p+1 }$  \
   respect permutations, see Lemma   \ref{lemma eins}, hence there are no contradictions in the definition of
   $\Theta_{n+1,1}$   on   ${ \bf BOU}_{n+1}   $,   the boundary of $ \Delta_{n+1}$.

   As in the case  \ $ n=2$ \ we use Proposition \ref{proposition drei}. \
        The reader should remember this proposition again. Before we use   Proposition \ref{proposition drei}
        we have to check whether all the conditions of  this   proposition are fulfiled.
  \begin{lemma}  \label{lemma siebzehn}  \quad
           The constructed map   $\Theta_{n+1,1}|_{ { \bf BOU}_{n+1}} $   is  a homeomorphism on 
             ${ \bf BOU}_{n+1}$.   Further, it satisfies the  conditions \  $ {\bf \widehat{[2]}}$ 
             (respecting permutations) and  $ {\bf \widehat{[3]}}$ (keeping the order).   \
           If we restrict \ $\Theta_{n+1,1}|_{ { \bf BOU}_{n+1}}$ to \
           ${ \clubsuit}_{n+1,\frac{1}{2 \cdot (n+2)}}$, \ we get a homeomorphism
           $$ \Theta_{n+1,1}|_{ { \bf BOU}_{n+1}}   :       \
                 { \bf BOU}_{n+1} \cap {{ \clubsuit}_{n+1,\frac{1}{2 \cdot (n+2) } } }    \
                 \stackrel{\cong}{\longrightarrow}    \
                 { \bf BOU}_{n+1} \cap {{ \clubsuit}_{n+1, \frac{1}{2 \cdot (n+2) + 1 } } }   \   . $$                     \end{lemma}
   \begin{proof}     \quad
            The fact that   $ \Theta_{n+1,1}|_{ { \bf BOU}_{n+1}}:  \
            { \bf BOU}_{n+1} \stackrel{\cong}{\longrightarrow} { \bf BOU}_{n+1}$ and 
            the properties  \  $ {\bf \widehat{[2]}}$ and   $ {\bf \widehat{[3]}}$  \            
            have been discussed  in and after Lemma  \ref{lemma sechzehn}. 
                      
       The third claim means that
            if  \ $\Theta_{n+1,1} (y_0,y_1, \ldots, y_n, y_{n+1} )$ $= (z_0,z_1, \ldots , z_n, z_{n+1})$,  for \
            $(y_0,y_1, \ldots , y_{n+1}) \in {\bf BOU}_{n+1}$, and if \
            $y_j = \frac{1}{2 \cdot (n+2)}$ \ for any \ $j \in \{ 0, 1, \ldots, n,n+1 \}$, \ then it follows  that \
            $z_j = \frac{1}{2 \cdot (n+2)+1}$. \ This can be shown by using the definition of
            $\Theta_{n+1,1} \text{ on }  \ { { \bf BOU}_{n+1}}$.   The reader should also look at
           Figure \ref{picture16}.   We have to note that, if we restrict the following three maps to the corresponding
           `crosses', we have three homeomorphisms \\
    \centerline{  $ \Theta_{n,0}^{-1}: {{ \clubsuit}_{n,\frac{1}{2 \cdot (n+2)}}} \stackrel{\cong}{\longrightarrow}
                  {{ \clubsuit}_{n,\frac{1}{2 \cdot (n+1)}}} , $   \qquad
                  $ \Theta_{n,1}: {{ \clubsuit}_{n,\frac{1}{2 \cdot (n+1)}}} \stackrel{\cong}{\longrightarrow}
                  {{ \clubsuit}_{n,\frac{1}{2 \cdot (n+1) + 1}}}, $ \qquad and    }   \\
    \centerline{  $ \Theta_{n+1,0}: {{ \clubsuit}_{n+1,\frac{1}{2 \cdot (n+2)}}} \stackrel{\cong}{\longrightarrow}
                  {{ \clubsuit}_{n+1,\frac{1}{2 \cdot (n+3)}}} $ .    }
   \end{proof}
        Now we are prepared to use Proposition \ref{proposition drei}.  
        Let \  $ \alpha := \frac{1}{2 \cdot (n+2)}   \ \text{and}   \ \beta := \frac{1}{2 \cdot (n+2)+1 } $ .
        Hence by Proposition \ref{proposition drei}, the constructed map
     $\Theta_{n+1,1}|_{ { \bf BOU}_{n+1}}$    can be extended to a homeomorphism
     $\Theta_{n+1,1}$   on  $ \Delta_{n+1} $, even  $\Theta_{n+1,1} \in { \cal COMFORT}(\Delta_{n+1})$, \
     with the property   
   $$ \Theta_{n+1,1}|_{{\clubsuit}_{n+1,\frac{1}{2 \cdot (n+2)}}}: \
             {\clubsuit}_{n+1,\frac{1}{2 \cdot (n+2) }}
             \stackrel{\cong}{\longrightarrow} {{\clubsuit}_{n+1, \frac{1}{2 \cdot (n+2) + 1 }}} \ . $$ 
    
      By the construction of  $\Theta_{n+1,1}|_{ { \bf BOU}_{n+1}}$,  the equations
      { \tt EQUATION$_{n+1,j \leq p,i,k}$}  are satisfied for $(i,k) \in \{ (0,1),(1,0)\}$.  For  $ i=k=0$,
      { \tt EQUATION$_{n+1,j \leq p,0,0}$} is trivial. It remains to consider the case
      $ i=k=1 $, i.e. we want to  show    { \tt   EQUATION$_{n+1,j \leq p,1,1}$}.
       With  the homeomorphism     $ \Theta_{n+1,1}: {{\clubsuit}_{n+1,\frac{1}{2 \cdot (n+2) }}}
               \stackrel{\cong}{\longrightarrow}  { { \clubsuit}_{n+1, \frac{1}{2 \cdot (n+2) + 1 }}}$ \
       it is easy to see that the following  Figure \ref{picture17} \\  \\ \\  \\ \\   \\  
    \begin{figure}[ht]
      \centering
      \setlength{\unitlength}{1cm}
      \begin{picture}(9.5,3)
         \put(-3.0,2.5){$\vec{x}$ }   \put(-2.5,3.4){$\vector(2,1){2}$}   \put(-2.0,3.2){$\Theta_{n,1}$}  
         \put(-0.3,4.4){$\vec{z}$ }   \put(0.25,4.5){$\vector(1,0){0.75}$} 
                                  \put(0.2,3.5){$\left\langle id \right\rangle_{\: n+1,\; 1,\; 0 \; }$}  
         \put(1.4,4.4){$ \left( \frac{1}{ 2 \cdot (n+2) }, \left( 1 - \frac{1}{ 2 \cdot (n+2) } \right)
                                            \cdot \vec{z} \right) $ }  
         \put(6.0,4.5){$\vector(1,0){0.75}$}   \put(6.0,3.65){$\Theta_{n+1,1}$}  
                                  \put(6.8,4.4){$\left( \frac{1}{ 2 \cdot (n+2)+1} , \vec{w} \right)$}
         \put(10.0,4.0){$\left\langle id \right\rangle_{\: n+2,\;  1,\; 1 \; }$}  \put(9.2,4.2){$\vector(2,-1){1.3}$}    
         \put(-2.5,1.7){$\vector(2,-1){2}$}    \put(-2.0,1.5){$\Theta_{n,1}$} 
         \put(-0.3,0.5){$\vec{z}$ }    \put(0.2,1.25){$\left\langle id \right\rangle_{\: n+1,\;  1,\; 0 \; }$}  
         \put(1.4,0.4){$ \left( \frac{1}{ 2 \cdot (n+2) }, \left( 1 - \frac{1}{ 2 \cdot (n+2) } \right)
                                            \cdot \vec{z} \right) $ } 
         \put(6.0,0.5){$\vector(1,0){0.75}$}    \put(0.25,0.5){$\vector(1,0){0.75}$}                              
         \put(6.0,1.25){$\Theta_{n+1,1}$}    \put(6.8,0.4){$\left( \frac{1}{ 2 \cdot (n+2)+1} , \vec{w} \right)$ } 
         \put(9.2,0.8){$\vector(2,1){1.3}$}   \put(10.0,0.80){$\left\langle id \right\rangle_{\: n+2,\; 1,\; 0 \; }$}            \put(7.25,2.5){$ \left( \frac{1}{ 2 \cdot (n+3)} , \frac{1}{ 2 \cdot (n+3)} , 
                                    \left( 1 - \frac{1}{ 2 \cdot (n+3)} \right) \cdot \vec{w} \right) $ } 
      \end{picture}     \caption{}\label{picture17} 
   \end{figure}     
            \quad  \\ \\
     commutes. Hence we have confirmed that also  {\tt EQUATION$_{n+1,j=0 \leq p=0,i=1,k=1}$} holds. 
     The other cases of   {\tt EQUATION$_{n+1,j \leq p,1,1}$} work correspondingly for  
     $ j,p \in \{ 0,1,2, \ldots, n, n+1  \}$  with   $ j \leq p $. 
                      
   Now  all cases of the equations  { \tt   EQUATION$_{n+1,j \leq p,i,k}$}  are proved, for
     $ i,k \in \{ 0,1 \} , $   and   $ j,p \in \{ 0,1,2, \ldots, n, n+1  \}$   with   $ j \leq p $.   \\  \\
   As a brief summary  we fixed  $ i=0, k=1 $, and we varied $j,p \in \{0,1,2, \ldots ,   n , n+1 \} $ 
     with   $ j \leq p$. Thereby we were able to define  the  homeomorphism   $\Theta_{n+1,1}$ on the  boundary  
     ${ \bf BOU}_{n+1}$ without   contradictions.   
     After that we could extend $\Theta_{n+1,1}$ on the entire $\Delta_{n+1}$ by 
     Proposition \ref{proposition drei}.  The map \ $\Theta_{n+1,1}$ has all  properties demanded in 
     Theorem   \ref{zweites theorem},  and we confirmed all equations  { \tt   EQUATION$_{n+1,j \leq p,i,k}$}. 
     
  Thus we have done the induction step from   $ n$   to   $n+1$,   and the proof  of 
     Theorem \ref{zweites theorem} is finished.
    \end{proof}
  With Theorem  \ref{zweites theorem} we have completed the proof of Theorem \ref{theorem eins}, i.e. 
      we have proved that  \  $_{\vec{m}}\partial_{n} \; \circ \; _{\vec{m}}\partial_{n+1} (T)  =  0 $  \ 
      holds for an arbitrary  $T \in {\cal SS}_{n+1}(X)$ \ in the case of    { \bf $ L = 1$ }! 
    \newpage
  \section{The Homology Modules of a Point}
        The reader should recall the definitions of   ${\cal SS}_{n}(X)$ and  of
        ${\cal F}\left({\cal R}\right)_{n}(X) $  from the introduction of this paper.
	 	    Let us take an arbitrary  chain {\tt u } $\in {\cal F}\left({\cal R}\right)_{n}(X)$,
	      i.e. {\tt u } is a linear combination of some  $T$'s from the set  ${\cal SS}_{n}(X)$.
	      Now let us assume that for fixed $L \in \;\mathbbm{N}_{0}$ \ and $ \vec{m} \in  {\cal R}^{L+1}$ \
	      we proved $_{\vec{m}}\partial_{n-1} \;  \circ \; 	 _{\vec{m}}\partial_{n} (T)  =  0$ \ for all basis
	      elements  $T \in {\cal SS}_{n}(X)$,  as we just have done in the case $L = 1$.
	      We can extend the boundary operators   $_{\vec{m}}\partial_{n}$ \
	      from   ${\cal SS}_{n}(X) $   to   ${\cal F}\left({\cal R}\right)_{n}(X)$  \
	      by linearity, hence we obtain for each chain  {\tt u }  the equation \
	      $_{\vec{m}}\partial_{n-1} \;  \circ \; _{\vec{m}}\partial_{n} ({\tt u }) = 0$, \ and for the chain complex
     \[ \ \cdots  \            \labto{_{\vec{m}}\partial_{n+1}} {\cal F}({\cal R})_{n}(X)
                               \labto{_{\vec{m}}\partial_{n}}   {\cal F}({\cal R})_{n-1}(X)
                               \labto{_{\vec{m}}\partial_{n-1}}  \ \cdots \
                               \labto{_{\vec{m}}\partial_2}    {\cal F}({\cal R})_{1}(X)
                               \labto{_{\vec{m}}\partial_1}    {\cal F}({\cal R})_{0}(X)
                               \labto{_{\vec{m}}\partial_0} \{ 0 \}
     \]
          we can deduce that \  $ \text{image}(	_{\vec{m}}\partial_{n+1 })$ \ is a submodule of \
          $ \text{kernel}(	_{\vec{m}}\partial_{n})$, \ hence the ${\cal R}$-module
       $$  {\cal H}_{n}(X) :=  \frac {kernel(	_{\vec{m}}\partial_{n })}{image(	_{\vec{m}}\partial_{n+1} )}  $$
          is  well defined for all fixed  $L \in  \mathbbm{N}_0$ and  fixed  tuples
          $\vec{m}  = \left( m_{0}, m_{1}, m_{2}, \dots  ,m_{L} \right)  \in \: {\cal R}^{L+1}$,
	        for every topological space   $X$   and all   $n \in  \;\mathbbm{N}_{0}$.  \\
	    Example :     \quad
                    Let   ${\cal R}$  := $\mathbbm{Z}$ .  \
                    For the one-point space  $\{p\}$   and for  $n \in  \mathbbm{N}_0$   there is only one
                    $T: \Delta_n \rightarrow \{p\}$,    thus it holds \
                    $ {\cal F}\left(  \mathbbm{Z}\right)_{n}(p) \cong  \mathbbm{Z}$. \
                     And for the generated chain complex    
     \[ \ \cdots  \            \labto{_{\vec{m}}\partial_{n+1}} {\cal F}(\mathbbm{Z})_{n}(p)
                               \labto{_{\vec{m}}\partial_{n}}   {\cal F}(\mathbbm{Z})_{n-1}(p)
                               \labto{_{\vec{m}}\partial_{n-1}}  \ \cdots \
                               \labto{_{\vec{m}}\partial_2}    {\cal F}(\mathbbm{Z})_{1}(p)
                               \labto{_{\vec{m}}\partial_1}    {\cal F}(\mathbbm{Z})_{0}(p)
                               \labto{_{\vec{m}}\partial_0} \{ 0 \}
     \]
         we get
     \[ \ \cdots \cdots \      \labto{_{\vec{m}}\partial_{n+1}} \mathbbm{Z}
                               \labto{_{\vec{m}}\partial_{n}}   \mathbbm{Z}
                               \labto{_{\vec{m}}\partial_{n-1}}  \ \cdots \cdots
                               \labto{_{\vec{m}}\partial_3}     \mathbbm{Z}
                               \labto{_{\vec{m}}\partial_2}     \mathbbm{Z}
                               \labto{_{\vec{m}}\partial_1}     \mathbbm{Z}
                               \labto{_{\vec{m}}\partial_0} \{ 0 \}  .
     \]
               We abbreviate  \ $ \sigma := \sum_{i=0}^{L}  m_{i}  $.
               If we define the map    $\times \sigma: \mathbbm{Z}  \rightarrow  \mathbbm{Z} ,   x \mapsto
                \sigma \cdot x$,  we can describe the  boundary operators  by
        \[       _{\vec{m}}\partial_{n}
                   \cong \begin{cases} 0 & \text{  if  }   n   \text{ is odd, or }   \ n = 0  \\
                                       \times \sigma   & \text{  if  }   n   \text{ is even, but } n \neq 0 \, .
                         \end{cases}
        \]

        Explanation: In the definition of   \
        $  _{\vec{m}}\partial_{n}	 (T) =  \sum_{j =  0}^{n}   (-1)^{j}  \cdot \sum_{i = 0}^{L}
        m_{i} \cdot \left( \left\langle T\right\rangle_{\: L,\: n,\;  i,\; j \;} \circ  \Theta_{L,n-1,i}\right)$  \
        we have that \  $\sigma$ of the unique map from $\Delta_{n-1}$ to $\{p\}$ \
        cancel pairwise because of the alternating signs.  \   That means for   $\sigma \neq  0$:  
           \[  _{\vec{m}}{\cal H}_{n}(p)
                      \cong \begin{cases}  \mathbbm{Z} & \text{  if  } n  = 0   \\
                                        \mathbbm{Z}_{/(\sigma \cdot \mathbbm{Z})} & \text{  if } n  \text{ is odd } \\
                                        0  & \text{  if  }   n   \text{ is even and }   n \neq 0   ,
                      \end{cases}
           \]
      and  for    $\sigma = 0$    we get
            \  $ _{\vec{m}}{\cal H}_{n}(p) \cong \mathbbm{Z}$ \ for all  $n \in  \mathbbm{N}_0$ .  \\ \\
  	
   \end{document}